\documentclass{article} 

\usepackage[paperheight=11in, 
    paperwidth=8.5in,
    outer=1.25in,
    inner=1.25in,
    bottom=1.4in,
    top=1.4in]{geometry} 

\usepackage{amsmath, amsxtra, amsthm, amsfonts, amssymb}
\DeclareMathOperator{\gen}{\mathsf{gen}}
\DeclareMathOperator{\idx}{\mathsf{idx}}
\usepackage{dsfont}
\usepackage{stmaryrd}
\usepackage{mathrsfs}
\usepackage{graphicx}
\usepackage{xurl}
\usepackage{color}
\usepackage{quiver}
\usepackage[dvipsnames]{xcolor}
\usepackage[absolute,overlay]{textpos}
\usepackage{bbm}
\usepackage{tikz-cd}
\usepackage[utf8x]{inputenc}
\usepackage{csquotes}
\usetikzlibrary{decorations.pathmorphing} 
\usetikzlibrary{decorations.pathreplacing}
\usetikzlibrary{braids}
\usepackage{stmaryrd}
\usepackage[T1]{fontenc}
\usepackage{lmodern}

\usepackage{upgreek}
\usepackage{cprotect}
\usepackage{booktabs}
\usepackage{longtable}

\definecolor{DeepViolet}{RGB}{85,32,128}
\definecolor{BluePeriwinkle}{RGB}{105,115,205}

\usepackage{hyperref}
\hypersetup{
    colorlinks=true,
    linkcolor=DeepViolet,
    citecolor=BluePeriwinkle,
    urlcolor=magenta
}
\usepackage{marginnote}

\NewDocumentCommand{\leanlink}{v}{%
  \marginnote{%
    \href{#1}{%
      \includegraphics[height=1em]{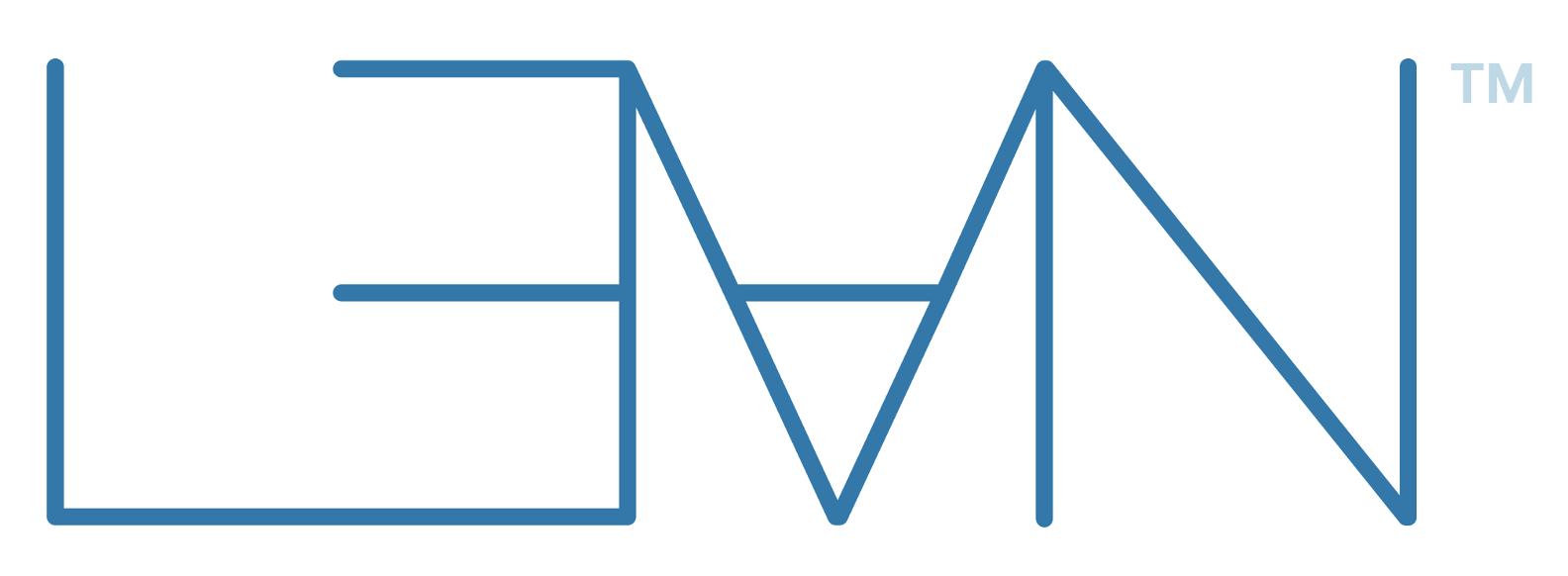}%
    }%
  }%
}

\NewDocumentCommand{\leanlinkfootnote}{m}{%
  \marginnote{%
    \href{#1}{%
      \includegraphics[height=0.85em]{lean-logo-teal.png}%
    }%
  }%
}
\definecolor{keywordcolor}{rgb}{0.7, 0.1, 0.4}   
\definecolor{commentcolor}{rgb}{0.4, 0.4, 0.4}   
\definecolor{symbolcolor}{rgb}{0.0, 0.1, 0.6}    
\definecolor{sortcolor}{rgb}{0.2, 0.5, 0.1}      
\definecolor{errorcolor}{rgb}{1, 0, 0}           
\definecolor{stringcolor}{rgb}{0.5, 0.3, 0.2}    

\usepackage{listings}

\usepackage[backend=biber, style=numeric]{biblatex} 
\usetikzlibrary{calc}

\AtEveryBibitem{
    \clearfield{urlyear}
    \clearfield{urlmonth}
}
\theoremstyle{definition}
\newtheorem{theorem}{Theorem}[section]

\newtheorem{definition}[theorem]{Definition}
\newtheorem{lemma}[theorem]{Lemma}

\newtheorem{conjecture}[theorem]{Conjecture}
\newtheorem{fallacy}[theorem]{Fallacy}
\newtheorem{fact}[theorem]{Fact}
\theoremstyle{remark}
\newtheorem*{remark}{Remark}
\title{Braid Groups in Lean : Formalizing Dehornoy's Subword Reversing Algorithm}
\author{
    Hannah Fechtner\\
    Carnegie Mellon University}
\date{}
\newcommand{\PML}{\mathsf{PML}}
\begin{document}
\maketitle
\begin{abstract}
We formalize the group-theoretic definition of infinite and finite Artin braid groups in Lean. We implement a polynomial-time procedure for determining equivalence in the braid group known as Dehornoy's subword reversing algorithm. Termination and correctness of this algorithm are shown using rectangular grids, which are a specific form of van Kampen diagram. The literature only hints at a proof linking rectangular grids to braid equivalence. We present a novel, three-stage proof of this result, introducing partial grids and associated size metrics.
\end{abstract}
\tableofcontents
\newpage
\section{Introduction} 

We begin by defining the basic objects of study, including the algebraic definition of braids, the Lean theorem prover, and the subword reversing algorithm. At the end of this chapter, we give an overview of the entire project.

\subsection{Mathematical Braids}

Braids were given their first full definition and exposition as a mathematical object by Emil Artin in 1925 \cite{artin_theorie_1925}. He gave both a topological and algebraic definition of these physical objects. A physical braid is a set of strings in three-dimensional space, affixed at top and bottom to the ``floor'' and ``ceiling.'' These strings can be twisted together or crossed, but must always flow downwards. No doubling-back is allowed. This is the sort of image to keep in mind:
\begin{center}
\begin{tikzpicture}
\pic[
braid/.cd,
every strand/.style={ultra thick},
gap=0.2,
strand 1/.style={Rhodamine},
strand 2/.style={YellowGreen},
strand 3/.style={Periwinkle},
strand 4/.style={Mulberry},
] {braid={s_1^{-1} s_3 s_2^{-1} s_1 s_2^{-1} s_1  s_3 s_2}};
\end{tikzpicture}
\end{center}

The topological definition lies outside the scope of this project. We work solely with Artin's algebraic definition, occasionally using images of the corresponding physical objects to aid in exposition.

\begin{definition} [Presented Group]
    Consider a set $S$ and a binary relation $R$ on elements of the free group over $S$, $F_S$. The \textbf{presented group} $\langle S\,| R\rangle$ is the quotient of $F_S$ by the congruence generated over $R$ (the reflexive, symmetric, transitive, multiplicative closure of $R$). 
\end{definition}

\begin{definition}
    Artin's braid group on $n$ strands is given by the presentation

\[ B_n :=
\left\langle 
\sigma_1, \ldots, \sigma_{n-1} \;\middle|\;
\begin{array}{ll}
\sigma_i \sigma_j = \sigma_j \sigma_i & \text{for } |i-j| \ge 2,\\[4pt]
\sigma_i \sigma_j \sigma_i = \sigma_j \sigma_i \sigma_j & \text{for } |i-j| = 1
\end{array}
\right\rangle
\]
\noindent The infinite braid group $B_\infty$ can be defined as the direct limit $B_\infty := \bigcup\limits_{n \in \mathbb{N}} B_n$ or equivalently  by the presentation

\[ B_\infty:=
\left\langle 
\sigma_1, \sigma_2,\sigma_3 \ldots \;\middle|\;
\begin{array}{ll}
\sigma_i \sigma_j = \sigma_j \sigma_i & \text{for } |i-j| \ge 2,\\[4pt]
\sigma_i \sigma_j \sigma_i = \sigma_j \sigma_i \sigma_j & \text{for } |i-j| = 1
\end{array}
\right\rangle
\]

\end{definition}

The generator $\sigma_i$ represents the strand in position $i$ ($1\leq i < n$) crossing over the strand to its right. 

\begin{center}
\begin{tikzpicture}
\pic[
braid/.cd,
every strand/.style={ultra thick},
gap=0.2,
strand 1/.style={Rhodamine},
strand 2/.style={YellowGreen},
] {braid={s_1}};
\end{tikzpicture}
\end{center}

\noindent ${\sigma_i}^{-1}$ represents the strand in position $i$ crossing under the strand to its right.

\begin{center}
\begin{tikzpicture}
\pic[
braid/.cd,
every strand/.style={ultra thick},
gap=0.2,
strand 1/.style={Rhodamine},
strand 2/.style={YellowGreen},
] {braid={s_1^{-1}}};
\end{tikzpicture}
\end{center}

For example, the braid at the opening of this chapter is coded ${\sigma_1}^{-1}\sigma_3{\sigma_2}^{-1}\sigma_1{\sigma_2}^{-1}\sigma_1\sigma_3\sigma_2$.
\vspace{2em}

\noindent The physical representation shows how inverses untangle : here is $\sigma_i^{-1}\sigma_i$ :

\begin{center}
    \begin{tikzpicture}
\pic[
braid/.cd,
every strand/.style={ultra thick},
gap=0.2,
strand 1/.style={Rhodamine},
strand 2/.style={YellowGreen},
] {braid={s_1^{-1} s_1}};
\end{tikzpicture}
\hspace{1cm} is equivalent to \hspace{1cm}
\begin{tikzpicture}
\pic[
braid/.cd,
every strand/.style={ultra thick},
gap=0.2,
strand 1/.style={Rhodamine},
strand 2/.style={YellowGreen},
strand 3/.style={White},
strand 4/.style={White},
] {braid={s_3 s_3}};
\end{tikzpicture}
\end{center}

\noindent The defining relations of the braid group can be visualized as follows :

Crossings far apart can be slid up and down (if $|i-j|\geq2,$ $ \sigma_i\sigma_j$ is equivalent to $\sigma_j\sigma_i$) : 
\medskip

\begin{tikzpicture}
\pic[
braid/.cd,
gap=0.2,
strand 1/.style={Rhodamine, {ultra thick}},
strand 2/.style={YellowGreen, ultra thick},
strand 3/.style={Gray, dashed},
strand 4/.style={Gray, dashed},
strand 5/.style={Periwinkle, ultra thick},
strand 6/.style={Mulberry, ultra thick},
] {braid={s_1 s_5}};
\end{tikzpicture}
\hspace{.5cm}  is equivalent to \hspace{.5cm} 
\begin{tikzpicture}
\pic[
braid/.cd,
every strand/.style={ultra thick},
gap=0.2,
strand 1/.style={Rhodamine},
strand 2/.style={YellowGreen},
strand 3/.style={Gray, dashed},
strand 4/.style={Gray, dashed},
strand 5/.style={Periwinkle},
strand 6/.style={Mulberry},
] {braid={s_5 s_1}};
\end{tikzpicture}
\medskip

An adjacent strand can be pulled through a crossing ($\sigma_i\sigma_{i+1}\sigma_i$ is equivalent to $\sigma_{i+1}\sigma_i\sigma_{i+1}$) :

\begin{center}
    \begin{tikzpicture}
\pic[
braid/.cd,
every strand/.style={ultra thick},
gap=0.2,
strand 1/.style={Rhodamine},
strand 2/.style={YellowGreen},
strand 3/.style={Periwinkle},
] {braid={s_1 s_2 s_1}};
\end{tikzpicture}
\hspace{1cm}  is equivalent to \hspace{1cm} 
\begin{tikzpicture}
\pic[
braid/.cd,
every strand/.style={ultra thick},
gap=0.2,
strand 1/.style={Rhodamine},
strand 2/.style={YellowGreen},
strand 3/.style={Periwinkle},
] {braid={s_2 s_1 s_2}};
\end{tikzpicture}
\end{center}

Of course, to rigorously show that these relations are sufficient to capture equivalence in the physical objects, one must give a topological definition and the proof that the algebraic and topological definitions coincide. The interested reader can find a modern exposition in Patrick Dehornoy's textbook ``The Calculus of Braids'', chapters one and two \cite{dehornoy_calculus_2022}.

\subsection{The Lean Theorem Prover}
Lean is a functional programming language and interactive theorem prover based on dependent type theory. 
Leonardo de Moura launched Lean in 2013; the language is now on version 4 \cite{platzer_lean_2021}. Mathematicians worldwide use Lean to formally verify proofs across mathematics, from undergraduate basics to the cutting edge of research. 
In Lean, one can implement algorithms and execute them using the \lstinline{#eval} command along with stating and proving theorems about their properties, all in the same file. While algorithms of course must be computable in order to evaluate, Lean is not limited to computable functions or constructive proofs. One can choose to use Lean with constructive or classical logic. When proving correctness of an algorithm, for example, the algorithm must be computable but the proofs of correctness and termination need not be.

A large collaborative project to formalize mathematics in one central repository, named Mathlib, is ongoing, and currently counts 772 contributors \cite{the_mathlib_community_lean_2020}.\footnote{\href{https://leanprover-community.github.io/mathlib_stats.html}{Mathlib statistics, as of October 12, 2026}} As far as we can tell, there has been no formalization of Artin's braid group in Mathlib or elsewhere in Lean, nor in any other theorem prover. Coxeter groups, another type of presented group which is similar in structure to the braid groups, have recently been formalized in Lean.\footnote{See the project page at \url{https://www.majiajun.org/formalizing-coxeter-group-hecke-algebra-and-kazhdan-lusztig-theory-in-lean/} and the Mathlib file at \url{https://leanprover-community.github.io/mathlib4_docs/Mathlib/GroupTheory/Coxeter/Basic.html}}. There has been work in formalized category theory which could be related to braid groups \cite{shah_vicar_2024, pous_string_2026, carlier_unbiasing_2026}, including basic definitions in Lean\footnote{\href{https://leanprover-community.github.io/mathlib4_docs/Mathlib/CategoryTheory/Monoidal/Braided/Basic.html}{Braid Monoidal Categories Mathlib page}}, but the major required theorem\footnote{From Joyal and Street \cite{joyal_braided_1993} : Letting $X$ be the single generating object in the free braided monoidal category, the set of endomorphisms on the tensor product of $n$ copies of $X$ is isomorphic to the braid group on $n$ strands.} is not yet formalized, and the word problem in such categorical structures is as hard as that of the unknot problem \cite{delpeuch_word_2022}, compared to polynomial-time complexity for braid groups \cite{dehornoy_efficient_2008}. Thus, a connection to category theory does not simplify the project, and lies beyond the scope of this paper.
 Knots and a few invariants thereof were formalized in Isabelle/HOL by Prathamesh in 2015 \cite{prathamesh_formalising_2015}. However, the mathematics used to study braids (presented groups, localizations, orderings) differs from that used for knots (polynomial invariants). There are also a number of computational packages for working with braids, for example in MATLAB \cite{thiffeault_braidlab_2019}, SAGE\footnote{\url{https://doc.sagemath.org/html/en/reference/groups/sage/groups/braid.html}}, and GAP \cite{michel_development_2013}. Lastly, two theorems about braid groups (unrelated to this project) very recently appeared on Prove2Me,\footnote{\href{https://prove2.me/missions/Braids\%2C\%20Links\%20and\%20Mapping\%20Class\%20Groups\%20I\%3A\%20Artin\%27s\%20Presentation\%20of\%20the\%20Braid\%20Group}{Prove2Me link}} an AI-based website where users employ their agents to solve mathematical ``missions'' proposed by humans. Upon a brief perusal, we were unable to find the Lean code on the Prove2Me site, and so are unsure of what definition of the braid groups was used.\footnote{This project's Lean definition of braid groups has been available on the internet since late 2024, as a \href{https://github.com/leanprover-community/mathlib4/pull/18626}{pull request to Mathlib}.}

As to the formalization of a word problem decision procedure, existing similar work is limited. A recent paper on decidability in presented monoids does offer a (slow) Rocq algorithm (\lstinline{homog_dec}) which can be used for any presented monoid whose relations preserve length \cite{cirpons_certifying_2026}. This could be applied to determine equivalence in braid \textit{monoids} but cannot be extended to braid groups. 

\subsection{Word Problem Algorithm}

The word problem in a group is to determine whether two representatives represent the same equivalence class. In general, the word problem of a presented group is unsolvable \cite{boone_word_1959, novikov_algorithmic_1955}. For the braid group in particular, this problem is solvable in polynomial time. A large number of algorithms exist, in part prompted by the early-2000s hope that braid group problems could be used for post-quantum security protocols \cite{lee_braid_2004,myasnikov_braid-based_2004}. Early algorithms with exponential complexity included Artin combing \cite{artin_theory_1947} and Garside's method of universal denominators \cite{garside_braid_1969}, although both methods have subsequently received updates which run in polynomial time \cite{gonzalez-meneses_polynomial_2018,epstein_braid_1992}.
A more modern solution by Birman, Ko, and Lee introduced an alternative presentation of the braid group, yielding a normal form with computational advantages \cite{birman_new_1998}. Other methods such as Dynnikov coordinates \cite{dynnikov_yang-baxter_2002, dehornoy_efficient_2008} and the Lawrence-Krammer representation \cite{krammer_braid_2002, bigelow_braid_2000} give exceptionally simple algorithms which work quite efficiently, but the proofs of correctness involve advanced topology. Dehornoy's most recent algorithm is handle reduction, which is experimentally quite fast (but there is no proof of polynomial-time complexity) \cite{dehornoy_efficient_2008}.

For this project, we have implemented Dehornoy's earlier method: subword reversing (as introduced in a journal article \cite{dehornoy_deux_1992} and elucidated in his textbook \cite{dehornoy_calculus_2022}). While this may not be the fastest method, it runs in quadratic time, relies on a manageable amount of theory, and could be generalized to larger classes of groups \cite{dehornoy_subword_2011}. Moreover, this solution builds infrastructure that could be reused in a future implementation of other algorithms, such as handle reduction.

The primary contribution lies in the treatment of this algorithm’s termination. Dehornoy wrote about this algorithm over a period of approximately twenty-five years. Most of Dehornoy's work presents a compelling sketch of why subword reversing terminates via associating rectangular grid diagrams to individual instances of reversing, but does not give a proof for this method. Bridging this gap required the development of a new proof of termination, done in three stages. It does appear that in one of Dehornoy's earliest papers he gives a different proof which may be fully rigorous as written. This proof is quite delicate (we will discuss this further in Chapter \ref{chap:grids_rw}). It seems that Dehornoy preferred the visual intuition and simplicity of the grid-based method, so we are pleased to have been able to fill in this method.

The heart of the algorithm is a reversing procedure on braid words. To underscore that we are working with words and not braid group elements, we use $\overline{\sigma_i}$ in place of $\sigma_i^{-1}$. The reversing procedure consists of repeatedly rewriting negative-positive pairs $\overline{\sigma_i}\sigma_j$ as follows:

\hspace{3em} If $|i - j|=0$, $\overline{\sigma_i}\sigma_i$ is deleted (more precisely, replaced by the empty word $\varepsilon$).

\hspace{3em} If $|i-j|>1$, $\overline{\sigma_i}\sigma_j$ is replaced by $\sigma_j\overline{\sigma_i}$

\hspace{3em} If $|i-j|=1$,  $\overline{\sigma_i}\sigma_j$ is replaced by $\sigma_j\sigma_i\overline{\sigma_j}\,\overline{\sigma_i}$.

Note that these rewritings preserve equivalence under the braid relations. The proofs are straightforward algebraic manipulations; here is an instance of the third case, with numerical indices for readability:
\begin{align*}
    \overline{\sigma_2}\sigma_3 &\equiv\overline{\sigma_2}\sigma_3
    \sigma_2\overline{\sigma_2} &(\text{as } \varepsilon \equiv \sigma_2\overline{\sigma_2})\\
    &\equiv\overline{\sigma_2}\sigma_3 
    \sigma_2\sigma_3\overline{\sigma_3}\,\overline{\sigma_2} &(\text{as } \varepsilon \equiv\sigma_3\overline{\sigma_3})\\
    &\equiv \overline{\sigma_2}\sigma_2
    \sigma_3\sigma_2\overline{\sigma_3}\,\overline{\sigma_2}&(\text{as } \sigma_3\sigma_2\sigma_3\equiv \sigma_2\sigma_3\sigma_2)\\
    & \equiv \sigma_3\sigma_2\overline{\sigma_3}\,\overline{\sigma_2} & (\text{as } \overline{\sigma_2}\sigma_2\equiv \varepsilon)\\
\end{align*}

\noindent Now, let us show reversing in action. We begin with $\overline{\sigma_1}\,\sigma_0 \overline {\sigma_2}\sigma_1$; here is the sequence of rewritings:
\begin{align*}
    &\overline{\sigma_1}\,\sigma_0 \overline {\sigma_2}\sigma_1 \\
& \sigma_0\sigma_1 \overline{\sigma_0}\,\overline{\sigma_1}\, \overline {\sigma_2} \sigma_1 \\
& \sigma_0\sigma_1 \overline{\sigma_0}\,\overline{\sigma_1} \sigma_1 \sigma_2\overline {\sigma_1} \,\overline{\sigma_2}\\
& \sigma_0\sigma_1 \overline{\sigma_0}\sigma_2\overline {\sigma_1}\, \overline{\sigma_2} \\
& \sigma_0\sigma_1 \sigma_2 \overline{\sigma_0} \,\overline {\sigma_1} \,\overline{\sigma_2}
\end{align*}

So, reversing $\overline{\sigma_1}\,\sigma_0 \overline {\sigma_2}\sigma_1$ gives $\sigma_0\sigma_1 \sigma_2 \overline{\sigma_0} \,\overline {\sigma_1} \,\overline{\sigma_2}$.  Note that the final word is in positive-negative form, with no negative-positive pairs appearing. 

To determine if two braid group words are equivalent, we may simply check if $ab^{-1}$ is equivalent to the empty word.  The full subword reversing algorithm will reverse its input, exchange the positive and negative parts, and reverse again, returning true if that final reversing produces the empty string, and false otherwise. 

\subsection{Outline}

\leanlink{https://github.com/hannahfechtner/braids_better/tree/098309403cda787eaf1fa954e07ff0ca3b2927ee} Throughout the text, clickable Lean icons in the righthand margin link to the relevant lines in the GitHub formalization.  The repository, frozen at the time of this paper's writing, can be found at \url{https://github.com/hannahfechtner/braids_better/tree/098309403cda787eaf1fa954e07ff0ca3b2927ee}. Information about the Lean and Mathlib versions needed to build the code are available in the \texttt{README} file. 

Chapters 2--5 discuss formalization of known results, following Patrick Dehornoy's textbook \textit{The Calculus of Braids} \cite{dehornoy_calculus_2022}. Chapters 6--8 expand on Dehornoy's work and contain new material necessary to implement his subword reversing algorithm in a theorem prover. 

Chapter 2 develops the algebraic framework needed for the reversing algorithm, including the braid group, braid monoids, and Semi-Thue systems. Before we implement the algorithm itself, we spend the next five chapters preparing theorems which will be used to show its termination and correctness. 

Chapters 3--5 prove the localization result that connects equivalence in the braid monoid with equivalence in the braid group, a connection that will be essential to the proof of correctness. Chapter 3 develops the necessary theory of Ore localization. Chapters 4 and 5 then show that the braid monoid satisfies the hypotheses required for this localization: Chapter 4 formalizes Dehornoy's proof of the existence of common multiples, while Chapter 5 proves cancellativity. The proof of cancellativity introduces grid diagrams, which encode braid monoid equivalences and will play a central role throughout the remainder of the project.

Chapters 6 and 7 develop the structures and metrics needed to prove termination of the reversing algorithm. Chapter 6 extends the interpretation of grid diagrams by introducing partial grids and grid-style rewriting, and establishes the correspondence between these new structures and reversing rewrites. Chapter 7 defines metrics on these structures. The material in these two chapters will yield a well-founded measure that decreases as reversing proceeds, providing the basis for the termination proof. The constructions and results developed in these two chapters appear, to the best of my knowledge, to be novel.

Chapter 8 at last brings all these strands together in the implementation of the subword reversing algorithm, first for braid monoid words and then for braid group words. The results of the preceding chapters are combined to prove correctness and termination. The procedure is also extended to the finite braid group. Finally, Chapter 9 concludes with directions for future work.

\newpage
\section{Defining Braids}

This chapter develops the algebraic structures that will be used throughout the remainder of the paper. We begin by defining Artin--Tits groups and specializing this construction to Artin's infinite braid group. We then turn to two structures necessary for subword reversing. We begin with presented monoids, and we define the infinite braid monoid and establish its basic properties. We next discuss Semi-Thue systems, which we will use to model string rewriting. Together, these definitions and results provide the algebraic foundation for the subword reversing algorithm. Lastly, we return to all the structures discussed, this time in the finite setting. We introduce the finite braid groups and monoids and develop mappings between the infinite and finite braid group/monoids. 

\subsection{Artin-Tits Groups}

Braid groups are a special case of Artin-Tits groups.
\begin{definition} \leanlink{https://github.com/hannahfechtner/braids_better/blob/098309403cda787eaf1fa954e07ff0ca3b2927ee/BraidProject/ArtinTits.lean#L28-L30}
Given a semigroup $G$ and $a,b\in G$,\textbf{ Artin-Tits relations} are those of the form $aba\, ...=bab\, ...$, where each side of the equality has the same length.
    
\end{definition}
\begin{definition}\leanlink{https://github.com/hannahfechtner/braids_better/blob/098309403cda787eaf1fa954e07ff0ca3b2927ee/BraidProject/ArtinTits.lean#L48-L49}
    \textbf{Artin-Tits groups} are presented groups $\langle S|R\rangle$ in which all relations R are Artin-Tits, with at most one relation for each pair of generators. We may succinctly package the information for the relations into a function $m : (S\times S) \to \mathbb{N}$, whose graph can be arranged into a symmetric matrix. Then all relations are of the form 
    \[
\underbrace{a b a\cdots}_{\text{length } m(a,b)}
\;=\;
\underbrace{b a b \cdots}_{\text{length } m(a,b)}
\]
for generators $a, b\in S$, viewed as elements of the free group over $S$. If no Artin-Tits relation exists between generators $a$ and $b$, we will set $m(a,b) := 0$.\footnote{Often in this case $m(a,b) := \infty$, which is not suitable for formalization.} 
\end{definition}

Artin-Tits groups are very similar to Coxeter groups, which have the additional property that for every generator $a$, $a^2=1$. We are thus able to lightly edit much of the work done by Newell Jensen, Mitchell Lee, and Óscar Álvarez, who formalized the definition and many properties of Coxeter groups.\footnote{\url{https://github.com/leanprover-community/mathlib4/blob/971f540b1b5f10f11d8f73b8d43986759b9f6dcd/Mathlib/GroupTheory/Coxeter/Basic.lean}}

In Lean, following Jensen, Lee, and Álvarez, we give the set of generators as the type \lstinline{α}, and the matrix $m$ by the  \lstinline{ArtinTitsMatrix} structure. The \lstinline{off_diagonal} property ensures that different generators are distinct.

\leanlink{https://github.com/hannahfechtner/braids_better/blob/098309403cda787eaf1fa954e07ff0ca3b2927ee/BraidProject/ArtinTits.lean#L10-L17}\begin{lstlisting} 
structure ArtinTitsMatrix (α : Type*) where
  M : Matrix α α ℕ
  isSymm : M.IsSymm := by decide
  off_diagonal i i' : i ≠ i' → M i i' ≠ 1 := by decide
\end{lstlisting}

\noindent A simple recursive definition produces an alternating string of group elements:\footnote{This is generalized to monoid elements in the code}

\leanlink{https://github.com/hannahfechtner/braids_better/blob/098309403cda787eaf1fa954e07ff0ca3b2927ee/BraidProject/Additions/Monoid.lean#L6-L8}
\begin{lstlisting}
def Group.alternate [Monoid G] (a b : G) (n : ℕ) : G := match n with
  | 0 => 1
  | n' + 1 => a * alternate b a n'
\end{lstlisting}

We generate the Artin-Tits relation for any two generators with \lstinline{ArtinTits.relation}; instead of giving an equality, we give a free group element equivalent to the identity:    \[
(\underbrace{a b a\cdots}_{\text{length } m(a,b)})\cdot 
(\underbrace{b a b \cdots}_{\text{length } m(a,b)})^{-1}
\]

\noindent The occurrences of \lstinline{FreeGroup.of} in the below definition lift an element of a set to the corresponding free group element.

\leanlink{https://github.com/hannahfechtner/braids_better/blob/098309403cda787eaf1fa954e07ff0ca3b2927ee/BraidProject/ArtinTits.lean#L28-L29}
\begin{lstlisting}
def ArtinTits.relation (M : ArtinTitsMatrix α) (a b : α) : FreeGroup (α) :=
  Group.alternate (FreeGroup.of a) (FreeGroup.of b) (M a b) * (Group.alternate (FreeGroup.of b) (FreeGroup.of a) (M a b))⁻¹
\end{lstlisting}
We gather all the relations into a single set:

\leanlink{https://github.com/hannahfechtner/braids_better/blob/098309403cda787eaf1fa954e07ff0ca3b2927ee/BraidProject/ArtinTits.lean#L31-L32}
\begin{lstlisting}
def ArtinTits.relation_set : Set (FreeGroup α) :=
  Set.range (Function.uncurry (relation M))
\end{lstlisting}
Finally, using the already-existing \lstinline{PresentedGroup} structure from Mathlib, we define \lstinline{ArtinTitsGroup} which will build the corresponding Artin-Tits group given a matrix $m$.

\leanlink{https://github.com/hannahfechtner/braids_better/blob/098309403cda787eaf1fa954e07ff0ca3b2927ee/BraidProject/ArtinTits.lean#L48}
\begin{lstlisting}
def ArtinTitsGroup (M : ArtinTitsMatrix α) := PresentedGroup (relation_set M)
\end{lstlisting}

At this point, we must decide whether or not to duplicate the \lstinline{PresentedGroup} API, specializing it to Artin-Tits groups. If the construction were not canonical, we would duplicate the API so that their implementation would be opaque to users. Here, the definition as a presented group is canonical, and one can expect the user to know that an Artin-Tits group is a presented group. Often in this case no custom API is given, especially when the code is intended to go into Mathlib, for this leads to unnecessary extra code to maintain. There is a bit of leeway for non-Mathlib code, and also when the specialization leads to simpler hypotheses for certain theorems. The latter is true for the universal property for Artin-Tits groups. We recall the universal property for presented groups in general : 

\begin{theorem}\label{thm:presented_group_universal_property}
(Universal Property for Presented Groups)
Given a set of generators $S$ and a set of relations $R$ on $S^*$, the universal property for the presented group
$\langle S\,|R\rangle$ states that for any group $G$ and function
$f : S \to G$ respecting the relations $R$, there exists a unique
group homomorphism $g : \langle S\,|R\rangle\to G$ such that
$ g(\llbracket s\rrbracket_R)=f(s)$
for every $s\in S$.
\end{theorem}

Here, we decide to implement a custom universal property for the Artin-Tits groups, as we can simplify the hypothesis that the function to be lifted respects the relations. 

\begin{theorem} \leanlink{https://github.com/hannahfechtner/braids_better/blob/098309403cda787eaf1fa954e07ff0ca3b2927ee/BraidProject/ArtinTits.lean#L65}
    (Universal Property for Artin-Tits Groups)
    Given a set of generators $S$, a set of Artin-Tits relations $R$, and a function $m$ which generates $R$, the universal property for the Artin-Tits group $\langle S\,|R\rangle$ states that for any group $G$ and function $f : S \to G$ such that for all $a,b\in S$, 
        \[
\underbrace{f(a) \,f(b) \,f(a)\cdots}_{\text{length } m(a,b)}
\;=\;
\underbrace{f(b) \,f(a) \,f(b) \cdots}_{\text{length } m(a,b)}
\]
    
    there exists a unique group homomorphism $g : \langle S\,|R\rangle\to G$ such that
$ g(\llbracket s\rrbracket_R)=f(s)$
for every $s\in S$. 
\end{theorem}

In Lean, we can phrase this property in terms of the matrix $m$ and \lstinline{Group.alternate}, without needing to reference \lstinline{ArtinTits.relation} and \lstinline{ArtinTits.relation_set}; we simply require that the function $f$ satisfies

\leanlink{https://github.com/hannahfechtner/braids_better/blob/098309403cda787eaf1fa954e07ff0ca3b2927ee/BraidProject/ArtinTits.lean#L53-L55}
\begin{lstlisting}
def ArtinTits.IsLiftable (M : ArtinTitsMatrix α) [Group G] (f : α → G) : Prop :=
  ∀ i j, Group.alternate (f i) (f j) (M i j) = Group.alternate (f j) (f i) (M i j)
\end{lstlisting}

Then we can prove that \lstinline{IsLiftable} functions satisfy the liftability requirement for the presented group universal property. We do so by unpacking the definition of \lstinline{ArtinTits.relation_set} and \lstinline{ArtinTits.relation}; that way, they may remain opaque to the user. Then our custom universal property needs only an instance of \lstinline{IsLiftable} :

\leanlink{https://github.com/hannahfechtner/braids_better/blob/098309403cda787eaf1fa954e07ff0ca3b2927ee/BraidProject/ArtinTits.lean#L63-L74}
\begin{lstlisting}
def ArtinTits.toGroup  {G : Type*} [Group G] {f : α → G} (M : ArtinTitsMatrix α) (hf : IsLiftable M f) : ArtinTitsGroup M →* G 
\end{lstlisting}

\noindent By contrast, we do not duplicate the rest of the \lstinline{PresentedGroup} API (such as extensionality and induction principles) as doing so would not simplify theorems for the user. 

\subsection{The Braid Group}

We now implement the specific case of the ``infinite'' braid group, hereafter simply referred to as the braid group. The set of generators is $S=\{\sigma_i\, | \, i \in \mathbb{N}\}$ (we will say $i$ is the \textit{index} of the generator $\sigma_i$). This set is represented by the natural numbers in Lean. This does require 0-indexing the generators, in contrast to Artin's original definition. If this is bothersome, it is simple at the end of the project to write a wrapper shifting all indices up by one. Recall the braid relation:

\begin{definition} \leanlink{https://github.com/hannahfechtner/braids_better/blob/098309403cda787eaf1fa954e07ff0ca3b2927ee/BraidProject/BraidGroup.lean#L24-L34}
    The braid relation $R_B$ is defined as the smallest relation satisfying the following two clauses. Let $i,j\in \mathbb{N}$. Then
\begin{itemize}
\item If  $|i-j|=1$, then 
 $R_B (\sigma_i\sigma_j\sigma_i, \,\sigma_{j}\sigma_i\sigma_{j})$  (the braid relation, Artin-Tits of length 3)

\item If $|i-j|\geq2,$  then $R_B(\sigma_i\sigma_j, \, \sigma_j\sigma_i)$ (the commutative relation, Artin-Tits of length 2)
\end{itemize}
\end{definition}
That is, if $i=j$, the relation has length 0 (it is unneeded); if $|i-j|=1$ the relation has length 3; otherwise, the relation has length 2. In Lean, we write the corresponding Artin-Tits matrix as follows, and obtain the braid group using the \lstinline{ArtinTits} API :

\leanlink{https://github.com/hannahfechtner/braids_better/blob/098309403cda787eaf1fa954e07ff0ca3b2927ee/BraidProject/BraidGroup.lean#L24-L33}
\begin{lstlisting}
def BraidMatrixInf : ArtinTitsMatrix ℕ where
  M := Matrix.of fun i j : ℕ ↦
    if i = j then 0 else (if i.dist j = 1 then 3 else 2)
  isSymm := by
    grind [Matrix.IsSymm, Matrix.transpose, Matrix.of_apply,  
      EmbeddingLike.apply_eq_iff_eq, Nat.dist]
  off_diagonal := by aesop
\end{lstlisting}

\leanlink{https://github.com/hannahfechtner/braids_better/blob/098309403cda787eaf1fa954e07ff0ca3b2927ee/BraidProject/BraidGroup.lean#L35}
\begin{lstlisting}
def BraidGroupInf := ArtinTitsGroup BraidMatrixInf
\end{lstlisting}

In our mathematical writeup, we have let $\sigma_i$ represent an element of the set of generators $S$, as well as an element of the free group over $S$ when defining $R$. To give this as an element of the braid group, we would write $\llbracket \sigma_i\rrbracket_{R_B}$. In Lean, we begin with a natural number \lstinline{(i : ℕ)} as a generator. We write \lstinline{FreeGroup.of i} for the free group element. Finally, \lstinline{(PresentedGroup.of i : BraidGroupInf)} will give us the braid group element. We abbreviate  \lstinline{(PresentedGroup.of i : BraidGroupInf)} by \lstinline{σ i} as so :

\leanlink{https://github.com/hannahfechtner/braids_better/blob/098309403cda787eaf1fa954e07ff0ca3b2927ee/BraidProject/BraidGroup.lean#L53}
\begin{lstlisting}
def σ (k : ℕ) : BraidGroupInf := PresentedGroup.of k
\end{lstlisting}

While convenient in Lean, \lstinline{σ i} is unfortunately nearly-overlapping notation with $\sigma_i$, though they represent quite different mathematical objects. To belabor the point,  the former (\lstinline{σ i}) is an element of the braid group, while the latter ($\sigma_i$) is an element of the generating set.

For the braid group, we give a specialized API. Again, the construction of the braid group as a presented group is not meant to be hidden. The custom API will make the rest of the project move more smoothly --- no need to constantly unfold the definition of \lstinline{BraidGroupInf}. For example, for the universal property, we may omit mention of \lstinline{BraidMatrixInf} and simplify the required hypothesis to 

\leanlink{https://github.com/hannahfechtner/braids_better/blob/098309403cda787eaf1fa954e07ff0ca3b2927ee/BraidProject/BraidGroup.lean#L154-L157}
\begin{lstlisting}
def BraidGroupInf.IsLiftable {G : Type*} [Group G] (f : ℕ → G) : Prop :=
  (∀ i j : ℕ, i.dist j = 1 → f i * f j * f i = f j * f i * f j) ∧
  (∀ i j : ℕ, 2 ≤ i.dist j → f i * f j = f j * f i)
\end{lstlisting}

Beyond specializations of the universal property, induction principles, et cetera,
we moreover prove useful lemmas such as $\llbracket \sigma_{i+1}\sigma_{i}\sigma_{i+1}\rrbracket =\llbracket \sigma_i\sigma_{i+1}\sigma_i\rrbracket$. In Lean, this is written as follows :

\leanlink{https://github.com/hannahfechtner/braids_better/blob/098309403cda787eaf1fa954e07ff0ca3b2927ee/BraidProject/BraidGroup.lean#L81-L91}
\begin{lstlisting}
theorem BraidGroupInf.braid {i : ℕ} :
    σ i * σ (i+1) * σ i = σ (i+1) * σ i * σ (i+1)
\end{lstlisting} 

Such lemmas are mathematically immediate, but in Lean require careful traversal of the structures at play (quotient groups, congruences, setoids, etc).

\subsection{Presented Monoids}

We next consider a structure akin to presented groups, but lacking inverses : presented monoids. These may be defined as a free monoid quotiented by the congruence generated by a given relation. That is the implementation chosen here, as we are able to build upon the \lstinline{Congruence} structure which already is in Mathlib as well as the rich API surrounding quotients. 

\begin{definition}[Presented Monoid] \leanlink{https://github.com/leanprover-community/mathlib4/blob/6015dbbfc399f4548bafd46f9560ff61595fd328/Mathlib/Algebra/PresentedMonoid/Basic.lean#L35-L41} 
    We let $\langle S\, | R\rangle^+$ denote the presented monoid on a set $S$ and relation $R$. We implement this by generating the free monoid on $S$, then quotienting by the congruence generated on $R$ (the reflexive, symmetric, transitive, multiplicative closure of $R$). Given $x\in S^*$, we write the equivalence class of $x$ as $\llbracket x\rrbracket^+_R$. We drop the subscript $R$ when it can be inferred from context. Given $x, y\in S^*$, we say $x \equiv_R^+y$ if they are equivalent under the congruence generated by $R$, once again dropping the subscript when possible.
\end{definition}

In Lean, we first must give a formal definition of $S$-strings.  We can select from lists, free monoids, strings, etc. However, \lstinline{Congruence} must be applied to a type endowed with a \lstinline{Mul} structure. \lstinline{FreeMonoid} is the only type in Mathlib that has been given a \lstinline{Mul} structure. So, we represent the $S$-strings by \lstinline{FreeMonoid} instead of \lstinline{List}. 

Then, given a set of relations \lstinline{rels} over a type \lstinline{α}, \lstinline{PresentedMonoid} constructs the monoid with
generators \lstinline{x : α} and relations \lstinline{rels} as a quotient of a congruence relation over \lstinline{rels}. Here \lstinline{conGen} generates a congruence.

\leanlink{https://github.com/leanprover-community/mathlib4/blob/cfc55fc58714a59dc06e6e79edc6d1d95140d3c2/Mathlib/Algebra/PresentedMonoid/Basic.lean#L40-L41}
\begin{lstlisting}
def PresentedMonoid (rels : FreeMonoid α → FreeMonoid α → Prop):= (conGen rels).Quotient
\end{lstlisting}

We provide basic properties of the presented monoid, such as extensionality, the universal property, and lifts of a function $f$ from $S^*$ to $\langle S\, | R\rangle^+$ using various formulations of the required conditions on $f$ (respecting the congruence on $R$).

Presented monoids are also called abstract rewriting systems, string rewriting systems, and Thue systems. They model the idea of rewriting parts of a string according to a given symmetric relation. (Such derivations will be discussed more fully in the following section, so the following definition is not yet linked to Lean code).

\begin{definition}\label{derivation} [$R$-derivation]
    Given a set $S$ and a binary relation $R$ on $S^*$, an $R$-derivation is a finite list of $S$-strings $R_0, R_1, ..., R_n$ such that for every $k$, $0\leq k<n$, there exist $(x,y)\in R$ such that $R_k=axb$ and $R_{k+1}=ayb$.
\end{definition}

\begin{theorem}\label{thm:deriv_equiv_presented_monoid} \leanlink{https://github.com/hannahfechtner/braids_better/blob/098309403cda787eaf1fa954e07ff0ca3b2927ee/BraidProject/SemiThue.lean#L88-L102}
    For any \textit{symmetric} relation $R$, there exists an $R$-derivation linking $R_0$ to $R_n$ if and only if they are equivalent in the corresponding presented monoid : $\llbracket R_0\rrbracket^+ _{R}=\llbracket R_n\rrbracket^+ _{R}$
\end{theorem}

A presented monoid is thus a symmetric rewriting system under another name. We give a robust API for a user looking to prove two words are equivalent in a presented monoid, with over a dozen lemmas. For example, we show:
\begin{lemma}\label{lem:presented_monoid_rewrite_example}\leanlink{https://github.com/hannahfechtner/braids_better/blob/098309403cda787eaf1fa954e07ff0ca3b2927ee/BraidProject/Additions/PresentedMonoid.lean#L72-L73}
    For $x, y, z, w\in S^*$, if $(x, y)\in R$ and $z\equiv^+_Rw$, then $xz\equiv^+yw$.
\end{lemma}

Beyond inspiration for future use cases, theorem \ref{thm:deriv_equiv_presented_monoid} gives us an alternate induction principle for presented monoids: 
\begin{theorem}\label{thm:presented_monoid_better_induction} \leanlink{https://github.com/hannahfechtner/braids_better/blob/098309403cda787eaf1fa954e07ff0ca3b2927ee/BraidProject/Additions/PresentedMonoid.lean#L187-L199}
    Consider $a, b\in S^*$ such that $a\equiv^+_Rb$. Let $P$ be a property on pairs of elements of $S^*$.  $P(a,b)$ holds if the following conditions are met:
    \begin{itemize}
        \item \textit{(reflexivity)} $P(a,a)$ holds for all $a\in S^*$.
        \item \textit{(one-step reduction)} If $(x,y)\in R$, then for any $c,d\in S^*$, $P(cxd, cyd)$.
        \item \textit{(symmetric one-step reduction)} If $(x,y)\in R$, then for any $c,d\in S^*$, $P(cyd, cxd)$.
        \item \textit{(transitivity)} For all $a, b, c\in S^*$, if $P(a,b)$ and $P(b,c)$, then $P(a,c)$.
    \end{itemize}
\end{theorem}

\subsection{Semi-Thue Systems}
We return to the definition of $R$-derivations in Definition \ref{derivation}, now discussing their implementation in Lean. We could give a definition almost exactly mirroring Definition \ref{derivation}, for example requiring there exist some natural number $n$ and a list of $n+1$ $A$-strings, with the property that each consecutive pair of $A$-strings in the list is an instance of a one-step reduction. However, it is much cleaner to give an inductive definition. We begin by defining the $\xrightarrow[R]{}$ notation for the one-step reduction operator  : 

\begin{definition} \leanlink{https://github.com/hannahfechtner/braids_better/blob/098309403cda787eaf1fa954e07ff0ca3b2927ee/BraidProject/SemiThue.lean#L10}
    Consider a set $A$ and a binary relation $R$ on $A$-strings. For all $A$-strings $c$ and $d$, if $R(a, b)$ holds, then we have $cad \xrightarrow[R]{} cbd$.
\end{definition}

Next, we take its reflexive, transitive closure. 
\begin{definition} \leanlink{https://github.com/hannahfechtner/braids_better/blob/098309403cda787eaf1fa954e07ff0ca3b2927ee/BraidProject/SemiThue.lean#L8-L11}

Consider a set $A$ and a binary relation $R$ on $A$-strings. Then $\xrightarrow[R]{*} $ is the smallest relation satisfying : 
\begin{itemize}

\item reflexivity : For all $a\in A^*$, $a \xrightarrow[R]{*} a$

\item closure : For all $a, b\in A^*$, $a\xrightarrow[R]{} b$ implies  $a\xrightarrow[R]{*} b$

\item transitivity : For all $a, b, c \in A^*$, if  $a\xrightarrow[R]{*} b$ and  $b\xrightarrow[R]{*} c$, then  $a \xrightarrow[R]{*} c$
\end{itemize}
\end{definition}

This relation is used in the traditional definition of a Semi-Thue system.
\begin{definition}
    A Semi-Thue system consists of a tuple $(A, R)$ (where $A$ is a set and $R$ a binary relation over $A$-strings) along with the relation $\xrightarrow[R]{*}$. 
\end{definition}

In Lean, we call this \lstinline{SemiThue}. For this project, it is convenient to provide a second definition of Semi-Thue systems. First, we define a new relation $\underset{R}{\rightarrowtriangle}$ : 
\begin{definition} \leanlink{https://github.com/hannahfechtner/braids_better/blob/098309403cda787eaf1fa954e07ff0ca3b2927ee/BraidProject/SemiThue.lean#L15-L18}
    Given a set $A$ and a binary relation $R$ over $A$-strings, let $\underset{R}{\rightarrowtriangle}$ be the smallest relation satisfying the below two properties :
\begin{itemize}
\item reflexivity : For all $a\in A^*$, $a \underset{R}{\rightarrowtriangle} a$

\item one-step transitivity :
For all $a, b_1,x,b_2,y \in A^*$, if  $a\underset{R}{\rightarrowtriangle} b_1xb_2$ and  $xRy$, then  $a\underset{R}{\rightarrowtriangle} b_1yb_2$
\end{itemize}
\end{definition} 

Note that the second condition is equivalent to the following, which makes the pen-and-paper notation for recording derivations cleaner:
\begin{itemize}
    \item \textit{alternative one-step transitivity} : For all $a, b_1,x,b_2,y \in A^*$, if  $a\underset{R}{\rightarrowtriangle} b_1xb_2$ and $b_1xb_2 \xrightarrow[R]{} b_1yb_2$, then  $a\underset{R}{\rightarrowtriangle} b_1yb_2$
\end{itemize}

Whichever formulation of one-step transitivity we make, we note the initial word is always preserved throughout an $\underset{R}{\rightarrowtriangle}$ --derivation. This is the main difference between our two definitions. Their behavior otherwise is nearly identical : 

\begin{theorem}\label{thm:SemiThueDefsEq} \leanlink{https://github.com/hannahfechtner/braids_better/blob/098309403cda787eaf1fa954e07ff0ca3b2927ee/BraidProject/SemiThue.lean#L28-L39}
    For all $a,b\in A^{*}$,
    $a \underset{R}{\rightarrowtriangle} b$ 
    iff $a \xrightarrow[R]{*} b$.
\end{theorem} 

Theorem \ref{thm:SemiThueDefsEq} allows us to replace $\xrightarrow[R]{*}$ with $\underset{R}{\rightarrowtriangle}$ in the definition of Semi-Thue systems, to obtain a variant which we call  \lstinline{SemiThueDerivation} in Lean.

The definition for both $\xrightarrow[R]{*}$ and $\underset{R}{\rightarrowtriangle}$ is straightforward in Lean (for the latter, using the original one-step transitivity condition). As to the choice of type for $A$-strings, we make a different decision than that for presented monoids. Here, we use lists, called \lstinline{List} in Lean. The focus for Semi-Thue systems is on rewriting, which has a more computational flair. The API surrounding lists, with computationally-minded functions such as \lstinline{getLast}, is far easier to use. The downside is that any work discussing both presented monoids and Semi-Thue systems requires repeated conversions between lists and free monoids. The occasional hassle of dealing
with \lstinline{FreeMonoid.toList} and \lstinline{List.toFreeMonoid} is outweighed by the convenience of using \lstinline{List} for Semi-Thue systems.

\subsection{The ``Braid'' Monoid}

With both presented monoids and Semi-Thue derivations in place, we can now specialize these general constructions to the braid generators and relations. At this point, we have referenced a proof that the braid \textit{group} aligns with the topological definition of braids. Although it is tempting to think of a monoid in which the braid relations hold as a ``braid monoid'' (perhaps braids with only overcrossings), until we can connect such a monoid to the braid group mathematically, we have no geometric interpretation of the monoid. It is, for the moment, an arbitrarily defined mathematical structure with no connection to the braid group. 

\begin{definition} \leanlink{https://github.com/hannahfechtner/braids_better/blob/098309403cda787eaf1fa954e07ff0ca3b2927ee/BraidProject/BraidMonoidInf.lean#L39}
The ``infinite'' braid monoid $B^+$, from now on simply referred to as the braid monoid, is given by the following presentation:

\[ B^+:=
\left\langle 
\sigma_0, \sigma_1 ,\sigma_2 \ldots \;\middle|\;
\begin{array}{ll}
\sigma_i \sigma_j = \sigma_j \sigma_i & \text{for } |i-j| \ge 2,\\[4pt]
\sigma_i \sigma_j \sigma_i = \sigma_j \sigma_i \sigma_j & \text{for } |i-j| = 1
\end{array}
\right\rangle^+
\]
\end{definition}

\begin{definition}
    The set of \textbf{braid monoid words} is 
    $$BW^+:=\left(\underset{i\in \mathbb{N}}{\cup}\sigma_i\right)^*$$ 
    
    We denote the \textbf{relation} used in the definition of the braid monoid as $R_{B^+}$.
\end{definition}

The absolute-value condition in the relations leads to a doubling of cases in many proofs. Our Lean implementation of the braid monoids is direct; we do not generalize to some sort of ``Artin-Tits monoid'', and thus there is no requirement for the relation to be symmetric.
To simplify future proofs, we define a different relation for our formalization.

\begin{definition} \leanlink{https://github.com/hannahfechtner/braids_better/blob/098309403cda787eaf1fa954e07ff0ca3b2927ee/BraidProject/BraidMonoidInf.lean#L13-L15}
    The simplified relation $R^s_{B^+}$ is defined as follows:
    \[
R^s_{B^+} :=
\begin{cases}
\sigma_i \sigma_{j} = \sigma_{j} \sigma_i & \text{for } j - i \geq 2\\
\sigma_i \sigma_{i+1} \sigma_i = \sigma_{i+1} \sigma_i \sigma_{i+1} & \text{for } j - i = 1
\end{cases}
\]
\end{definition}

\begin{fact}
    $\langle \sigma_0, \sigma_1, \sigma_2, ...|R^s_{B^+}\rangle=B^+$. We have the analogous property for $R^s_{B_n^+}$ and $B_n^+$.
\end{fact}
\begin{proof} (Sketch)
 Presented monoids are generated by taking the congruence over the given relation (hence the symmetric closure), and $R_{B^+}$ is the symmetric closure of $R^s_{B^+}$.
\end{proof}

For our Lean definition of the relations for $B^+$, note that in the second case $j=i+1$. We thus define the relation for the braid monoid as follows: 

\leanlink{https://github.com/hannahfechtner/braids_better/blob/098309403cda787eaf1fa954e07ff0ca3b2927ee/BraidProject/BraidMonoidInf.lean#L13-L15}
\begin{lstlisting}
inductive braid_monoid_rels_inf : FreeMonoid ℕ → FreeMonoid ℕ → Prop
  | adjacent (i : ℕ): braid_monoid_rels_inf (of i * of (i+1) * of i) (of (i+1) * of i * of (i+1))
  | separated (i j : ℕ) (h : i + 2 ≤ j) : braid_monoid_rels_inf (of i * of j) (of j * of i)
\end{lstlisting}

\noindent This can be plugged in to our previously-defined \lstinline{PresentedMonoid} definition, neatly giving us the braid monoid. 

We next define a few properties of braid monoid elements. 

\begin{definition} \label{def:braid_length} \leanlink{https://github.com/hannahfechtner/braids_better/blob/098309403cda787eaf1fa954e07ff0ca3b2927ee/BraidProject/BraidMonoidInf.lean#L87-L93}
    The \textbf{length} of a braid monoid element $b$ is written $|b|$. We define $|\llbracket a\rrbracket| $ as the length of $a$ (number of letters; the empty word $\varepsilon$ has length 0). This is well-defined, as a quick inspection of the relations shows length is preserved under the braid monoid equivalence. 
\end{definition}

\begin{definition} \leanlink{https://github.com/hannahfechtner/braids_better/blob/098309403cda787eaf1fa954e07ff0ca3b2927ee/BraidProject/BraidMonoidInf.lean#L113-L124}
    The \textbf{set of generators} of a braid monoid element is given by the function $\gen$. $\gen(\llbracket a\rrbracket) $ is defined as the set of generators appearing in $a$. This is similarly well-defined. 
\end{definition}

\begin{definition}
    The \textbf{index} of a braid monoid word $b$ is the index of the largest generator occurring in $b$. We denote this $\idx$, defining $\idx(b):= \max {\{i \, |\, \sigma_i \in \gen (b)\}}$. 
\end{definition}

\begin{lemma} \label{lem:reverse_equiv} \leanlink{https://github.com/hannahfechtner/braids_better/blob/098309403cda787eaf1fa954e07ff0ca3b2927ee/BraidProject/BraidMonoidInf.lean#L148-L204}
        The reversing function on free monoid words, $r$, reverses the order of the letters in a free monoid word. $a\equiv^+b$ if and only if $r(a)\equiv^+r(b)$. This is proved in each direction by structural induction on the equivalence. It is a useful tool to derive left common multiples and right cancellativity from right common multiples and left cancellativity respectively (or vice versa).
\end{lemma}

We will thus need to define length, set of symbols in, and reversing for free monoid elements. Let us return for a moment to the choice made earlier about the type for $S$-strings. We were forced to select free monoids instead of lists. We now see one downside to this choice : \lstinline{List} is much richer than \lstinline{FreeMonoid}. \lstinline{List} already has functions \lstinline{length}, \lstinline{toFinset}, \lstinline{reverse}, etc.
Although it is simple to implement these for \lstinline{FreeMonoid}, we must tread carefully. While free monoids are in fact defined as lists in Lean, they are meant to perform different roles. Free monoids generally bend more to mathematical theorems in abstract algebra, whereas lists are more used for computation. Even when the two share similarly-minded functionality, it is from a different viewpoint and with different implementation details. For example, if the idea is to apply a function to every element of our string and then combine them, we have \lstinline{FreeMonoid.lift} but \lstinline{List.foldr}. It is a balancing act to keep the two structures distinct, but also usable. Each has lemmas the other lacks.

We have decided to add definitions and lemmas to \lstinline{FreeMonoid} if they will be able to be defined on the eventual \textit{presented} monoid, the braid monoid. Such definitions seem reasonably ``algebraic'' since they can be defined for even more complex algebraic structures. Thus, we add in \lstinline{FreeMonoid.length}, \lstinline{FreeMonoid.reverse}, and \lstinline{FreeMonoid.symbols} (the set of generators appearing in an element), because these play nicely with the braid monoid structure. 

We take a different approach for \lstinline{List} definitions that do not remain invariant in the braid structure, such as \lstinline{getLast} (this returns the last element in a list). Equivalent braids need not have the same last generator --- $\sigma_1\sigma_3$ and $\sigma_3\sigma_1$ certainly do not! In this case, we have decided to use \lstinline{FreeMonoid.toList} when needed. This is a bit cumbersome, but thankfully quite rare for this project. 

We amass simple results of equivalence classes for words:
\begin{lemma}\label{lem:generator_equiv} \leanlink{https://github.com/hannahfechtner/braids_better/blob/098309403cda787eaf1fa954e07ff0ca3b2927ee/BraidProject/BraidMonoidInf.lean#L206-L213}
If $\sigma_i\equiv^+a$, then $a=\sigma_i$. Thus, $\llbracket \sigma_i \rrbracket=\{\sigma_i\}$
\end{lemma}
\begin{lemma}\label{lem:equiv_relation} \leanlink{https://github.com/hannahfechtner/braids_better/blob/098309403cda787eaf1fa954e07ff0ca3b2927ee/BraidProject/BraidMonoidInf.lean#L215-L298}
    Let $(w,x)\in R_{B^+}$. If $w\equiv^+a$, then $a=w$ or $a=x$. Similar reasoning for $x$ shows us $\llbracket w\rrbracket = \{w, x\} =\llbracket x\rrbracket$
\end{lemma}
\begin{proof}
    For the case of the commutative relation, $(\sigma_i\sigma_j, \sigma_j\sigma_i)\in R_{B^+}$ with $|i-j| \geq 2$, since we know that both length and set of generators is preserved under equivalence, our conclusion is immediate. For the case of the braid relation, $(\sigma_i\sigma_j\sigma_i, \sigma_j\sigma_i\sigma_j)\in R_{B^+}$ with $|i-j| = 1$, the information about length and set of generators is not enough as those allow for, say, $a = \sigma_j\sigma_i\sigma_i$. Instead, we induct on the equivalence $w\equiv^+a$ following Theorem \ref{thm:presented_monoid_better_induction}.
\end{proof}

For long elements of $BW^+$, it is frequently difficult to prove equivalence by the rewriting conditions for presented monoids (i.e. Lemma \ref{lem:presented_monoid_rewrite_example}), as parentheses need to be moved around so that Lean's rewrite tactic, \lstinline{rw}, can exactly pattern-match the terms. Once we move to a specific presented monoid, we have concrete relations to work with. We can associate all parentheses to the left, and then use one of the following lemmas:

\begin{lemma}\label{lem:comm_rel_rw} \leanlink{https://github.com/hannahfechtner/braids_better/blob/098309403cda787eaf1fa954e07ff0ca3b2927ee/BraidProject/BraidMonoidInf.lean#L366-L368}
    For all $x\in BW^+$, if $|i-j|>1$, then 
    $$\left(x\cdot\llbracket \sigma_i \rrbracket^+\right)\cdot\llbracket \sigma_j \rrbracket^+=
    \left(x\cdot\llbracket \sigma_j \rrbracket^+\right)\cdot\llbracket \sigma_i \rrbracket^+$$ 
\end{lemma}
\begin{lemma}\label{lem:braid_rel_rw} \leanlink{https://github.com/hannahfechtner/braids_better/blob/098309403cda787eaf1fa954e07ff0ca3b2927ee/BraidProject/BraidMonoidInf.lean#L375-L379}
For all $x\in BW^+$, if $|i-j|=1$, then 
    $$\left( \left(x\cdot\llbracket \sigma_i\rrbracket^+\right) \cdot\llbracket \sigma_j \rrbracket^+\right) \cdot\llbracket \sigma_i \rrbracket^+=\left( \left( x\cdot\llbracket \sigma_j\rrbracket^+\right) \cdot\llbracket \sigma_i \rrbracket^+\right) \cdot\llbracket \sigma_j \rrbracket^+$$
\end{lemma}

\noindent For example, say we are trying to transform the element $\llbracket\sigma_2\sigma_7\sigma_3\sigma_4\sigma_5\sigma_4\sigma_0\sigma_1\rrbracket$ into $\llbracket\sigma_2\sigma_7\sigma_5\sigma_3\sigma_4\sigma_5\sigma_0\sigma_1\rrbracket$. First, we separate all of the generators :
$\llbracket\sigma_2\rrbracket\cdot\llbracket\sigma_7\rrbracket\cdot\llbracket\sigma_3\rrbracket\cdot\llbracket\sigma_4\rrbracket\cdot\llbracket\sigma_5\rrbracket\cdot\llbracket\sigma_4\rrbracket\cdot\llbracket\sigma_0\rrbracket\cdot\llbracket\sigma_1\rrbracket$. Next, we associate all the parentheses to the left :

\newcommand{\pA}[1]{\textcolor{YellowGreen}{\raisebox{-0.04em}{\scalebox{1.1}{$($}}}#1\textcolor{YellowGreen}{\raisebox{-0.04em}{\scalebox{1.1}{$)$}}}}

\newcommand{\pB}[1]{\textcolor{Emerald}{\raisebox{-0.08em}{\scalebox{1.25}{$($}}}#1\textcolor{Emerald}{\raisebox{-0.08em}{\scalebox{1.25}{$)$}}}}

\newcommand{\pC}[1]{\textcolor{cyan}{\raisebox{-0.12em}{\scalebox{1.4}{$($}}}#1\textcolor{cyan}{\raisebox{-0.12em}{\scalebox{1.4}{$)$}}}}

\newcommand{\pD}[1]{\textcolor{Blue}{\raisebox{-0.16em}{\scalebox{1.6}{$($}}}#1\textcolor{Blue}{\raisebox{-0.16em}{\scalebox{1.6}{$)$}}}}

\newcommand{\pE}[1]{\textcolor{Mulberry}{\raisebox{-0.2em}{\scalebox{1.8}{$($}}}#1\textcolor{Mulberry}{\raisebox{-0.2em}{\scalebox{1.8}{$)$}}}}

\newcommand{\pF}[1]{\textcolor{WildStrawberry}{\raisebox{-0.24em}{\scalebox{2.0}{$($}}}#1\textcolor{WildStrawberry}{\raisebox{-0.24em}{\scalebox{2.0}{$)$}}}}
\[
    \pF{
      \pE{
        \pD{
          \pC{
            \pB{
              \pA{
                \llbracket\sigma_2\rrbracket\cdot\llbracket\sigma_7\rrbracket
              }\cdot\llbracket\sigma_3\rrbracket
            }\cdot\llbracket\sigma_4\rrbracket
          }\cdot\llbracket\sigma_5\rrbracket
        }\cdot\llbracket\sigma_4\rrbracket
      }\cdot\llbracket\sigma_0\rrbracket }
\cdot\llbracket\sigma_1\rrbracket
\]

\noindent Then, Lean's \lstinline{rw} tactic used with lemma \ref{lem:braid_rel_rw} will assign $\pB{
              \pA{
                \llbracket\sigma_2\rrbracket\cdot\llbracket\sigma_7\rrbracket
              }\cdot\llbracket\sigma_3\rrbracket}$ as $x$, and rewrite
              $      \color{Mulberry}\Bigl( \color{black}
        \left(
          \left(
            \left(  \left(
                \llbracket\sigma_2\rrbracket\cdot\llbracket\sigma_7\rrbracket
              \right)\cdot\llbracket\sigma_3\rrbracket
            \right) 
            \cdot\llbracket\sigma_4\rrbracket \right)\cdot\llbracket\sigma_5\rrbracket
        \right)\cdot\llbracket\sigma_4\rrbracket  
        \color{Mulberry}\Bigr)$ to 
$      \color{Mulberry}\Bigl( \color{black}
        \left(
          \left(
            \left( 
             \left(
                \llbracket\sigma_2\rrbracket\cdot\llbracket\sigma_7\rrbracket
              \right)\cdot\llbracket\sigma_3\rrbracket
            \right) 
            \cdot\llbracket\sigma_5\rrbracket \right)\cdot\llbracket\sigma_4\rrbracket
        \right)\cdot\llbracket\sigma_5\rrbracket \color{Mulberry}\Bigr)$, giving us \[
    \pF{
      \pE{
        \pD{
          \pC{
            \pB{
              \pA{
                \llbracket\sigma_2\rrbracket\cdot\llbracket\sigma_7\rrbracket
              }\cdot\llbracket\sigma_3\rrbracket
            }\cdot\llbracket\sigma_5\rrbracket
          }\cdot\llbracket\sigma_4\rrbracket
        }\cdot\llbracket\sigma_5\rrbracket
      }\cdot\llbracket\sigma_0\rrbracket }
\cdot\llbracket\sigma_1\rrbracket
\]

\noindent We then use \lstinline{rw} with lemma \ref{lem:comm_rel_rw}, explicitly setting $x$ as $\pA{
                \llbracket\sigma_2\rrbracket\cdot\llbracket\sigma_7\rrbracket
              }$ to obtain
    \[
    \pF{
      \pE{
        \pD{
          \pC{
            \pB{
              \pA{
                \llbracket\sigma_2\rrbracket\cdot\llbracket\sigma_7\rrbracket
              }\cdot\llbracket\sigma_5\rrbracket
            }\cdot\llbracket\sigma_3\rrbracket
          }\cdot\llbracket\sigma_4\rrbracket
        }\cdot\llbracket\sigma_5\rrbracket
      }\cdot\llbracket\sigma_0\rrbracket }
\cdot\llbracket\sigma_1\rrbracket
\]          

\noindent which can be simplified to $\llbracket\sigma_2\sigma_7\sigma_5\sigma_3\sigma_4\sigma_5\sigma_0\sigma_1\rrbracket$, as desired. The once-off hassle of separating generators and associating parentheses to the left makes multi-step calculations significantly simpler, without requiring tedious, targeted applications of associativity of multiplication.

\subsection {Finite Braid Groups and Finite Braid Monoids} \label{sec:finite_def}

Above, we have discussed only infinite braid groups and infinite braid monoids and the equivalences on them, $\equiv$ and $\equiv^+$. We may also consider the finite braid groups/monoids. Note that a braid on $n$ strands has $n-1$ generators. Those generators are numbered $\sigma_0, \sigma_1, ... \sigma_{n-2}$.  

\begin{definition} \leanlink{https://github.com/hannahfechtner/braids_better/blob/098309403cda787eaf1fa954e07ff0ca3b2927ee/BraidProject/BraidGroup.lean#L37}
The finite braid group $B_n$, representing braids on $n$ strands, is given by 
    \[ B_n :=
\left\langle 
\sigma_0, \ldots, \sigma_{n-2} \;\middle|\;
\begin{array}{ll}
\sigma_i \sigma_j = \sigma_j \sigma_i & \text{for } |i-j| \ge 2,\\[4pt]
\sigma_i \sigma_j \sigma_i = \sigma_j \sigma_i \sigma_j & \text{for } |i-j| = 1
\end{array}
\right\rangle
\]

We denote the defining relation (on the righthand side) as $R_{B_n}$ and the equivalence on elements of $B_n$ by $\equiv_n$.
\end{definition}
\noindent Note that $B_0$ and $B_1$ are both the trivial group. The finite braid monoids are given using the same generators and relations as for the finite braid groups :

\begin{definition} \leanlink{https://github.com/hannahfechtner/braids_better/blob/098309403cda787eaf1fa954e07ff0ca3b2927ee/BraidProject/BraidMonoidFin.lean#L8-L16}
\[ B_n^+ :=
\left\langle 
\sigma_0, \ldots, \sigma_{n-2} \;\middle|\;
\begin{array}{ll}
\sigma_i \sigma_j = \sigma_j \sigma_i & \text{for } |i-j| \ge 2,\\[4pt]
\sigma_i \sigma_j \sigma_i = \sigma_j \sigma_i \sigma_j & \text{for } |i-j| = 1
\end{array}
\right\rangle^+
\]

Here, the relation is denoted $R_{B_n^+}$ and the equivalence, $\equiv_n^+$.
\end{definition}

\begin{definition}
    We define the set of generators in $B_n$ (and $B^+_n$) as $S_n=\{\sigma_i | 0 \leq i \leq n-2\}$. Then the set of finite-braid group words is $BW_n := (S_n\cup\overline{S_n})^*$ and the set of finite-braid monoid words is $BW_n^+:=S_n^*$.
\end{definition}

Let us now move to the Lean implementation. 
For infinite braid groups and monoids, the generators are kept track of via their index, of type $\mathbb{N}$.
For $B_n$ and $B_n^+$, the indices of the generators are natural numbers, but bounded above by $n-1$ (in Lean, this is the predecessor of $n$, \lstinline{n.pred}). We thus use the \lstinline{Fin} type. It is a bundled type: an object of type \lstinline{Fin n.pred} contains a natural number $k$ and a proof that $k<n-1$. Aside from the change in type of the indices, the implementations of the finite braid groups and braid monoids is almost exactly akin to that of their infinite counterparts. Thus, there are clear similarities between the finite and infinite braid groups (and respective monoids). When we prove properties about one, it would be nice to have some kind of ``converter'' to port them to the other group/monoid. 

This is especially true as working with the finite braid groups/monoids is difficult in Lean. First, we are always juggling a hypothesis of the form $k<n-1$. Moreover, natural number subtraction is truncated, which makes many arithmetic proofs long. In the case of the braid relations, where one often must prove something of the form $i+2\leq j$, this complicates things enormously. Thus, we work within the finite braid monoid as sparingly as possible. Luckily, the relation on finite braids is a restriction on that for infinite braids. 

\begin{theorem} \leanlink{https://github.com/hannahfechtner/braids_better/blob/098309403cda787eaf1fa954e07ff0ca3b2927ee/BraidProject/BraidGroup.lean#L317-L323}
    For braid words $a,b$, if $a\equiv_n b$, then $a\equiv_m b$ for all $m>n$, and also $a\equiv b$. If $a,b$ are braid \textit{monoid} words, and if $a\equiv^+_n b$, then $a\equiv^+_m b$ for all $m>n$, and also $a\equiv^+ b$. \leanlink{https://github.com/hannahfechtner/braids_better/blob/098309403cda787eaf1fa954e07ff0ca3b2927ee/BraidProject/BraidMonoidFin.lean#L272-L279}
\end{theorem}

There is a small hiccup here : when we say $a\equiv_n b$ our formalization requires $a, b\in B_n$, but when we say $a\equiv_k b$ the formalization requires $a, b\in B_k$. It is quick to see that we have $BW_0\subset BW_1 \subset BW_2 ...\subset BW$, and so this poses little trouble to the paper-and-pencil mathematician. In Lean, we must cast from \lstinline{Fin n.pred} to \lstinline{Fin m.pred} or $\mathbb{N}$. 

When can we transfer an equivalence from the infinite braid group or monoid to the corresponding finite structure? Let us begin with the monoid case. 

\begin{theorem}\label{thm:braid_monoid_inf_to_fin} \leanlink{https://github.com/hannahfechtner/braids_better/blob/098309403cda787eaf1fa954e07ff0ca3b2927ee/BraidProject/ConvertMonoidInfToFin.lean#L18-L46}
    Consider $a,b\in BW^+$ such that $a\equiv^+ b$. Let $\mathsf{idx}(a)=i_a$ (that is, the largest generator index occurring in the word $a$ is $i_a$) and $\mathsf{idx}(b)=i_b$. Then for any $n\geq \max(i_a,i_b)$, we have $a\equiv^+_n b$.
\end{theorem}

The proof proceeds by induction on the derivation $a\equiv^+ b$. In our formalization, we must cast the infinite braid monoid element ``down" to a finite braid monoid element. We call this procedure \lstinline{FreeMonoid.mapNatToFin} :

\leanlink{https://github.com/hannahfechtner/braids_better/blob/098309403cda787eaf1fa954e07ff0ca3b2927ee/BraidProject/Additions/FreeMonoid.lean#L170-L171}
\begin{lstlisting}
def FreeMonoid.mapNatToFin (n : ℕ) (a : FreeMonoid ℕ) (bound : ∀ x ∈ a, x < n) : FreeMonoid (Fin n) 
\end{lstlisting}

The above theorems will eventually be used to apply our solver to the finite case and show that it is correct for finite braid group inputs. We leave the finite braid groups/monoids aside until the final chapter -- the algorithm is built for the infinite braid group, and then proved to work for the finite group as well. We are not yet able to give a theorem analogous to Theorem \ref{thm:braid_monoid_inf_to_fin} for the braid groups. In fact, as the very last step of our project, we will use the correctness specification of our solver to prove this.

In sum, we have defined the infinite and finite braid groups, as well as a number of structures needed to implement and reason about subword reversing (presented monoids and Semi-Thue systems). 

\section{Localization}
We next elucidate the connection between the braid group and the ``braid'' monoid by constructing an embedding of $B^+$ into $B$. In this chapter, we discuss Ore localization, the procedure by which this is done. We specialize it to our particular use case, first to localizing by the entire monoid, then to localizing by an entire \textit{presented} monoid. We merely state the required conditions on the (presented) monoid, and leave the proofs to the following two chapters. (Note that there is no discussion of the finite braid groups/monoids).

\subsection{Ore Localization}
Given a monoid $M$ and a submonoid $M_0$ of $M$, a localization procedure produces an associated group of fractions, essentially ``adding inverses'' for every element of $M_0$. The usual procedure  requires a commutative monoid. However, braid concatenation is not commutative.

For example, one may see that $\sigma_{1}\sigma_{2}$ does not equal $\sigma_{2}\sigma_{1}$ through properties of braid relations (no relations match $\sigma_{1}\sigma_{2}$). Or visually, neither of these two braids may be deformed into the other :

\begin{center}
    \begin{tikzpicture}
\pic[
braid/.cd,
every strand/.style={ultra thick},
gap=0.2,
strand 1/.style={Rhodamine},
strand 2/.style={YellowGreen},
] {braid={s_1 s_2}};
\end{tikzpicture}
\hspace{1cm}  and   \hspace{1cm}
\begin{tikzpicture}
\pic[
braid/.cd,
every strand/.style={ultra thick},
gap=0.2,
strand 1/.style={Rhodamine},
strand 2/.style={YellowGreen},
] {braid={s_2 s_1}};
\end{tikzpicture}
\end{center}

Instead, we use the Ore localization, which requires different monoid properties. We denote the (left)\footnote{There also exists a \textit{right} Ore localization $M[M_0^{-1}]$, used by Dehornoy. Here we choose to use the left localization, as that is what has already been implemented in Mathlib. Either can be used for the braid group, as it is cancellative and has common multiples on both sides.} Ore localization of $M$ by $M_0$ as $[M_0^{-1}]M$.

The following two conditions suffice for left Ore localization; we will refer to them as Ore conditions throughout this project:
\begin{definition}
    We say that a monoid $M$ and a submonoid $M_0$ of $M$ satisfy the Ore conditions if the following two criteria hold : 

\begin{enumerate}
    \item Left Common Multiples: For any $a \in M$, $b\in M_0$, there exist $c\in M_0$ and $d\in M$ such that $c a = d b$.
    \item Right Cancellativity : For $a,b \in M$ and $c\in M_0$, if $ac=bc$, then $a=b$.
\end{enumerate}
\end{definition}

The second condition can be slightly weakened.\footnote{We could merely require that common factors on the left be turned in to common factors on the right. That is, for $a, b \in  M$, $c\in M_0$, if $ac=bc$ then there exists $d\in M_0$ such that $da=db$.} Ore localization with a weakened cancellativity condition has already been implemented in Lean by Jakob von Raumer, with assistance from Kevin Klinge and Andrew Yang.\footnote{Their work can be found \hyperlink{https://github.com/leanprover-community/mathlib4/blob/f57201f9e0f3097e72753cdcc99f9917ea055339/Mathlib/RingTheory/OreLocalization/Basic.lean}{here}.} It is found under ring theory; the file first gives the case of a monoid. We will specialize this code to the case of localizing by the entire monoid, and then to presented monoids. We briefly summarize the mathematics from the formalization of von Raumer.

We begin with a monoid $M$ and a submonoid $M_0\subseteq M$ satisfying the Ore conditions.  We consider all pairs of the form $(x, y)$ with $x\in M$ and $y\in M_0$. We give the following equivalence relation on pairs:

\begin{definition} [$\sim$]
    We say $(x_1, y_1)\sim(x_2, y_2)$ iff there exist $a\in M$, $b\in M_0$ such that $ax_1=bx_2$ and $ay_1=by_2$. 
\end{definition}

\begin{theorem}
    $\sim$ is an equivalence relations on pairs $(x, y)$.
\end{theorem}

The proofs of the theorem follows directly from the definition of $\sim$. Next, we consider equivalence classes of pairs under $\sim$. These equivalence classes will be the elements of the localization. We use the abbreviation $\frac{x}{y}$ the equivalence class of $(x,y)$ under $\sim$.
\newcommand\sbullet[1][.75]{\mathbin{\vcenter{\hbox{\scalebox{#1}{$\bullet$}}}}}

\begin{definition}\label{def:ore_localization}
    Given a monoid $M$ and a submonoid $M_0$of $M$ satisfying the Ore conditions, we define the monoid $[M_0^{-1}]M$ as the set $\{\frac{x}{y}|\, x\in M,y\in M_0 \}$, with $\frac{1_M}{1_{M_0}}$ as the identity (note that $1_M$ is in fact equal to $1_{M_0}$) and a binary operation $\sbullet$ defined as follows: $\frac{x}{y} \sbullet \frac{z}{w} = \frac{az}{by}$ where $a\in M$ and $b\in M_0$ are such that $aw=bx$.
\end{definition}

The proof that the definition is independent of the choice of common multiples $a$ and $b$ is straightforward, as is that for associativity of the operator. We next consider a mapping from $M$ into $[M_0^{-1}]M$. 

\begin{theorem} \label{thm:ore_numeratorHom}
    The function $h$ from $M$ into $[M_0^{-1}]M$, defined by $h(x):=\frac{x}{1}$, is a homomorphism.
\end{theorem}

\begin{proof} We have $h(x_1x_2)= \frac{x_1x_2}{1}$. To calculate $h(x_1)\sbullet h(x_2)=\frac{x_1}{1}\sbullet\frac{x_2}{1}$, following definition \ref{def:ore_localization} we look for a left common multiple of $x_1$ and $1$. As $1\cdot x_1=x_1\cdot 1$, we may let $a=x_1$ and $b=1$. Hence $\frac{x_1}{1}\sbullet\frac{x_2}{1}=\frac{x_1\cdot x_2}{1\cdot 1}$. Thus, $h(x_1x_2)=h(x_1)\sbullet h(x_2)$.
\end{proof}

We now add one more condition. Together with our original Ore conditions, these form the cancellative Ore conditions:
\begin{enumerate}
    \setcounter{enumi}{2}
    \item Left Cancellativity: For $a,b\in M$ and $c\in M_0$, if $ca=cb$, then $a=b$.
\end{enumerate}

\noindent In the context of the cancellative Ore conditions, $h$ is an embedding of $M$ into $[M_0^{-1}]M$:
\begin{theorem} \label{thm:into_ore_hom_injective} \leanlink{https://github.com/hannahfechtner/braids_better/blob/098309403cda787eaf1fa954e07ff0ca3b2927ee/BraidProject/Additions/OreLocalization.lean#L3-L8}
Given a monoid $M$ with a submonoid $M_0$ satisfying the \textit{cancellative} Ore conditions, the function $h(x):=\frac{x}{1}$ is injective. 
\end{theorem}
\begin{proof}
    Assume  $h(x_1)=h(x_2)$; that is,$\frac{x_1}{1}=\frac{x_2}{1}$. Then $(x_1, 1)\sim(x_2, 1)$. There exist $a\in M, b\in M_0$ such that $ax_1=bx_2$ and $a\cdot1=b\cdot1$. Hence $a=b$ and so $bx_1=bx_2$. By left-cancellativity, $x_1=x_2$ as desired. 
\end{proof}

We now lay the full cancellative Ore conditions aside for the moment. For this project, we first specialize to the case of localizing by the entire monoid, viewed as a submonoid of itself. Mathlib already contains the typeclass \lstinline{OreSet}, which represents a pair of a monoid and a submonoid satisfying the (weakened) Ore conditions. Thus, we want to show that if a monoid has left common multiples and right cancellativity, then it paired with the full submonoid of itself is an instance of \lstinline{OreSet}.  To take in the hypotheses that the monoid M satisfies the left common multiple and right cancellativity criteria, we use mixins. These are typeclasses that extend a structure to add in an extra property. For example, we define

\leanlink{https://github.com/hannahfechtner/braids_better/blob/098309403cda787eaf1fa954e07ff0ca3b2927ee/BraidProject/Additions/Mixins.lean#L4-L5}
\begin{lstlisting}
class IsCommonLeftMultipleMul (M : Type*) [Mul M] where
  common_left_multiple : ∀ a b : M, ∃ c d : M, c * a = d * b
\end{lstlisting}

The above code is nonconstructive, as the existential quantifier merely asserts the existence of the common multiples. It does not provide them as data. This is in line with general Mathlib conventions. However, the existing Ore localization code is set up to be used constructively, because \lstinline{OreSet} has fields for the common-multiple data. We take advantage of this design, and give both constructive and non-constructive versions for braids. The constructive version will be useful for a concurrent project, developing a widget to visualize braids in the Lean InfoView. The widget will use the data from constructive proofs to display common multiples, equivalences, canonical forms, etc. Here is the analogous constructive version of \lstinline{IsCommonLeftMultipleMul} :

\leanlink{https://github.com/hannahfechtner/braids_better/blob/098309403cda787eaf1fa954e07ff0ca3b2927ee/BraidProject/DataCarrying/Mixins.lean#L3-L6}
\begin{lstlisting}
class HasCommonLeftMultipleData (M : Type*) [Mul M] where
  cl₁ : M → M → M
  cl₂ : M → M → M
  cl_spec : ∀ a b : M, cl₂ a b * a = cl₁ a b * b
\end{lstlisting}

We continue our discussion in the non-data-carrying setting only. The nonconstructive version is used in the implementation of the braid word problem solver, for Ore localization is used only in the proof of correctness of the algorithm, not in its definition. As such, it need not be constructive. We thus use the following mixin hypotheses added to a monoid $M$:

\leanlink{https://github.com/hannahfechtner/braids_better/blob/098309403cda787eaf1fa954e07ff0ca3b2927ee/BraidProject/OreLocalizationCombined.lean#L13}
\begin{lstlisting} 
variable {M : Type*} [Monoid M] [IsRightCancelMul M] [IsCommonLeftMultipleMul M]
\end{lstlisting}

\noindent With the above assumed, we are able to show that the full submonoid, \lstinline{(⊤ : Submonoid M)}, is an instance of \lstinline{OreSet} : 

\leanlink{https://github.com/hannahfechtner/braids_better/blob/098309403cda787eaf1fa954e07ff0ca3b2927ee/BraidProject/OreLocalizationCombined.lean#L16-L23}
\begin{lstlisting}
noncomputable instance : OreLocalization.OreSet (⊤ : Submonoid M) where ...
\end{lstlisting}
 The proof uses \lstinline{Classical.choose} to noncomputably provide the witnesses that the constructive \lstinline{OreSet} requires, and \lstinline{Classical.choose_spec} to assert that those witnesses have the desired property of being common multiples. 
 
 We now may define the Ore localization of $M$ by itself, as Lean can infer that the required Ore conditions are met from the above two lines of code:

\leanlink{https://github.com/hannahfechtner/braids_better/blob/098309403cda787eaf1fa954e07ff0ca3b2927ee/BraidProject/OreLocalizationCombined.lean#L44}\begin{lstlisting}
abbrev "OreLocalizationSelf" => OreLocalization (⊤ : Submonoid M) M
\end{lstlisting}

\begin{theorem} \label{thm:ore_group} \leanlink{https://github.com/hannahfechtner/braids_better/blob/098309403cda787eaf1fa954e07ff0ca3b2927ee/BraidProject/OreLocalizationCombined.lean#L46-L53}
    For any monoid $M$ satisfying the Ore conditions, $[M^{-1}]M$ has a group structure.
\end{theorem}

\begin{proof}
    We already have shown that $[M^{-1}]M$ has a monoid structure. All that remains is to define inverses. For any $\frac{x}{y}\in [M^{-1}]M$, we define $(\frac{x}{y})^{-1}:=\frac{y}{x}$. We omit the straightforward proof that $\frac{y}{x}$ satisfies the properties of an inverse.
\end{proof}

Next, we consider \textit{presented} monoids. 
We aim to plug in a presented monoid satisfying the Ore conditions to the existing Ore localization API. It is now clear why mixins are required. If we instead took in the hypotheses by defining an extension of the \lstinline{Monoid} class satisfying the Ore hypotheses called \lstinline{OreMonoid}, much like \lstinline{CancelMonoid}, we would run into an issue.  Both the \lstinline{PresentedMonoid} and \lstinline{OreMonoid} instances would carry with them a \lstinline{.toMonoid} function. If we try to give Lean a type \lstinline{M} with instances \lstinline{[PresentedMonoid M]} \lstinline{[OreMonoid M]}, Lean can infer the \lstinline{Monoid} instance on \lstinline{M} in two different ways. This is referred to as a ``bad diamond''. 

\begin{center}
\begin{tikzcd}
                           & \tt{Monoid} \arrow[ld] \arrow[rd] &                      \\
\tt{PresentedMonoid} \arrow[rd] &                              & \tt{OreMonoid} \arrow[ld] \\
                           & M                            &                     
\end{tikzcd}
\end{center}

Thus, we choose to use mixins to take in the hypotheses. 
\begin{theorem}  \label{thm:ore_presented_monoid_to_presented_group}
\leanlink{https://github.com/hannahfechtner/braids_better/blob/098309403cda787eaf1fa954e07ff0ca3b2927ee/BraidProject/OreLocalizationCombined.lean#L200-L210}
Let $M = \langle S \,| R\rangle ^+$ satisfy Ore conditions. Then $[M^{-1}]M$ is isomorphic to the presented group on those same generators and relations, $\langle S \,| R\rangle $. 
\end{theorem}

The proof will proceed by a careful inspection of the following universal properties of the Ore localization and of presented groups; we begin by defining those universal properties.  For simplicity, we let $\PML(S,R)$ denote $\left[(\langle S\, |  R\rangle ^+)^{-1}\right] \langle S\, |  R\rangle^+$.

\begin{definition} [Lift]
    Given a map $f$ from a set $S$ to a monoid $M$, we may lift $f$ to a function from the free group on $S$ to $M$, which we will denote $\overline{f}$. We denote the lift of $f$ to a function from the free monoid on $S$ as $\overline{f}^+$
\end{definition}

\begin{theorem} \label{thm:ore_localization_universal_property} [Universal property of the Ore localization of a presented monoid] \leanlink{https://github.com/hannahfechtner/braids_better/blob/098309403cda787eaf1fa954e07ff0ca3b2927ee/BraidProject/OreLocalizationCombined.lean#L81-L104}
    Consider a presented monoid $M=\langle S\mid R\rangle^+$ satisfying the Ore conditions.
For any group $G$ and any map $f:S\to G$ respecting the relations $R$
(for $(r_1,r_2)\in R$, we have
$\overline{f}^+(r_1)=\overline{f}^+(r_2)$),
there exists a unique homomorphism
\[
f^*:[M^{-1}]M\to G
\]
such that, for every $s\in S$,
\[
f^*\left(\frac{\llbracket s\rrbracket_R^+}{1}\right)
=f(s).
\]
\end{theorem}
\begin{proof}
    This is merely a specialization of the universal property for Ore localizations, which is again a specialization of that for quotients.
\end{proof}
\begin{remark}
    Since the identity map on $\PML(S,R)$ is a group
    homomorphism, the uniqueness clause of the universal
    property implies that any group homomorphism
    $g : \PML(S,R) \to \PML(S,R)$ satisfying
    $g\left(\frac{\llbracket s\rrbracket_R^+}{1}\right) = \frac{\llbracket s\rrbracket_R^+}{1} $ for all $ s \in S$ must be the identity map.
\end{remark}

\noindent Our next task will be to define two homomorphisms between the presented group $\langle S\, | R\rangle$ and the Ore localization of the associated presented monoid $\langle S\, | R\rangle ^+$.

\begin{definition} \label{def:ore_to_presented_group_maps} \leanlink{https://github.com/hannahfechtner/braids_better/blob/098309403cda787eaf1fa954e07ff0ca3b2927ee/BraidProject/OreLocalizationCombined.lean#L146-L148}
We define a homomorphism $h_1$ from $\langle S\, |  R\rangle $ to $\PML(S, R)$ using the universal property for presented groups to lift the function $f(s)=\frac{\llbracket s\rrbracket_{R}^+ :}{1_{{\langle S, R \rangle}^{+}}} : S \to \PML(S, R)$. The function $f$ is the composition of $h$ from Theorem \ref{thm:ore_numeratorHom} (recall $h(s)=\frac{s}{1}$) with $\llbracket \cdot \rrbracket _R^+$ (in Lean, this is \lstinline{PresentedMonoid.of} when applied to a single generator).

\leanlink{https://github.com/hannahfechtner/braids_better/blob/098309403cda787eaf1fa954e07ff0ca3b2927ee/BraidProject/OreLocalizationCombined.lean#L108-L112} We define a homomorphism $h_2$ from $\PML(S, R)$ to $\langle S\, |  R\rangle $ by using the universal property for Ore localizations of presented monoids to lift the function $m \mapsto \llbracket m\rrbracket_R$. 
\begin{remark}\label{rem:maps_to_itself} \leanlink{https://github.com/hannahfechtner/braids_better/blob/098309403cda787eaf1fa954e07ff0ca3b2927ee/BraidProject/OreLocalizationCombined.lean#L120-L132}
    Note that for any $s\in S^*$, $h_2\left(\frac{\llbracket s\rrbracket^+_R}{1}\right)=\llbracket s\rrbracket_R$. The proof proceeds by induction on $s$; the case where $s$ is a single generator is a consequence of theorem \ref{thm:ore_localization_universal_property}.
\end{remark}

\begin{center}
\begin{tikzcd}
                  & |[draw, circle,minimum size=5.8em]| {\PML(S, R)} \arrow[rr, "h_2", bend left=49] 
                  &  & |[draw, circle,minimum size=5.8em]| {\langle S\, | R\rangle } \arrow[ll, "h_1", bend left=49] \\
|[draw, circle]| {\langle S\, | R\rangle ^+} \arrow[ru, "h"] & & &
\end{tikzcd}
\end{center}
\end{definition}
Finally, we put the pieces together to complete the proof of Theorem \ref{thm:ore_presented_monoid_to_presented_group}. Composing $h_1$ with $h_2$ gives a homomorphism from $\langle S\, |  R\rangle $ to $\langle S\, |  R\rangle $ -- but by uniqueness portion of the universal property on presented group (Theorem \ref{thm:presented_group_universal_property}), that must be the identity! Similarly, composing $h_2$ with $h_1$ gives the identity by Theorem \ref{thm:ore_localization_universal_property}. Thus, $h_1$ is a bijection (and $h_2$ is as well). Hence we have an isomorphism between $\PML(S, R)$ and $\langle S\, |  R\rangle $. Thus, the localization of a presented monoid is isomorphic to the presented group on those same generators/relations, as desired. $\Box$

As to the Lean implementation of the above proof, we must note that the API for presented monoids takes a relation as a function of two \lstinline{FreeMonoid} arguments, while the API for presented groups requires the relations to be given as a set of \lstinline{FreeGroup} elements. To translate from the former to the latter, we use the function

\leanlink{https://github.com/hannahfechtner/braids_better/blob/098309403cda787eaf1fa954e07ff0ca3b2927ee/BraidProject/Additions/PresentedGroup.lean#L13-L15}
\begin{lstlisting}
def free_group_set_of_function (rels : FreeMonoid α → FreeMonoid α → Prop) : Set (FreeGroup α) 
\end{lstlisting}
We further provide a number of lemmas about the interaction of \lstinline{free_group_set_of_function} with lifts and the relations. For example, consider a function $f$ from $S$ to a group $G_1$. $\overline{f}$ applied to any element of \lstinline{free_group_set_of_function rels} equals one if and only if $\overline{f}^+ (r_1)=\overline{f}^+(r_2)$ for all $(r_1, r_2)$ for which \lstinline{rels} holds: 

\leanlink{https://github.com/hannahfechtner/braids_better/blob/098309403cda787eaf1fa954e07ff0ca3b2927ee/BraidProject/Additions/PresentedGroup.lean#L17-L37}
\begin{lstlisting}
theorem free_group_set_of_function_lift_eq_one_iff {G₁ : Type} [Group G₁]
        {rels : FreeMonoid α → FreeMonoid α → Prop} (f : α → G₁) :
    (∀ r ∈ free_group_set_of_function rels, (FreeGroup.lift f) r = 1) ↔
    (∀ r₁ r₂, rels r₁ r₂ → FreeMonoid.lift f r₁ = FreeMonoid.lift f r₂)
\end{lstlisting}

Such lemmas help us prove the required hypotheses for the universal properties in the definition of $h_1$ and $h_2$. The verification that $h_1$ and $h_2$ are inverses of one another is straightforward, showing that $\PML(S, R)$ and $\langle S\, |  R\rangle $ are isomorphic.

Returning to the mathematical development, we aim to show that for a cancellative presented monoid, equivalence in the associated presented group of positive words implies equivalence in the presented monoid of those same words.  For the rest of this section, we discuss only presented monoids which satisfy the \textit{cancellative} Ore conditions. (We say that the monoid $M$ satisfies the (cancellative) Ore conditions if $M$ considered with the full submonoid of itself satisfies the (cancellative) Ore conditions.)

\begin{definition} \label{def:iota}\leanlink{https://github.com/hannahfechtner/braids_better/blob/098309403cda787eaf1fa954e07ff0ca3b2927ee/BraidProject/OreLocalizationCombined.lean#L222-L225}
    We define the function $\iota : \langle S\,| R\rangle^+\to  \langle S\,| R\rangle$ by $\iota :=h_2\circ h$, where $h$ is the function from Theorem \ref{thm:ore_numeratorHom}.
\end{definition}
\begin{theorem}\label{thm:iota_embedding} \leanlink{https://github.com/hannahfechtner/braids_better/blob/098309403cda787eaf1fa954e07ff0ca3b2927ee/BraidProject/OreLocalizationCombined.lean#L230-L235}
    For a monoid $\langle S\,| R\rangle^+$ satisfying the cancellative Ore conditions, $\iota$ is an embedding of $\langle S\,| R\rangle^+$ into $\langle S\,| R\rangle$.
\end{theorem}
\begin{proof}
    By Theorem \ref{thm:into_ore_hom_injective}, $h$ is an injective homomorphism. $h_2$ is as well -- it is a homomorphism by its definition, and we exhibited above that $h_1$ is its two-sided inverse. Composition preserves both injectivity and the property of being a homomorphism, so $\iota$ is an embedding of $\langle S \,| R\rangle^+$ into $\langle S \,| R\rangle$. 
\end{proof}

Finally, we may show that equivalence in the group implies equivalence in the monoid for positive words :
\begin{theorem}\label{thm:group_equiv_monoid_equiv_ore} \leanlink{https://github.com/hannahfechtner/braids_better/blob/098309403cda787eaf1fa954e07ff0ca3b2927ee/BraidProject/OreLocalizationCombined.lean#L237-L248}
Let $\langle S\, | R\rangle$ be a left-cancellative monoid satisfying the Ore conditions. If $a\equiv_R b$ for \textbf{positive} words $a, b \in S^*$ (those containing no inverses)\footnote{Note that in our implementation we would need to use \lstinline{free_group_set_of_function}, and thus we more precisely assume $a\equiv_{\mathtt{free\_group\_set\_of\_function}(R)} b$ instead of the looser $a\equiv_R b$.}, then $a\equiv_R^+b$. 
\end{theorem}
\begin{proof}
    We know that $\llbracket a\rrbracket = \llbracket b \rrbracket$. By Remark \ref{def:ore_to_presented_group_maps} about $h_2$, we know that $\llbracket a\rrbracket = h_2(\frac{\llbracket a\rrbracket^+}{1})$, and $\llbracket b\rrbracket = h_2(\frac{\llbracket b\rrbracket^+}{1})$. Thus $\frac{\llbracket a\rrbracket^+}{1}=\frac{\llbracket b\rrbracket^+}{1}$. By definition of $h$, we have $h(\llbracket a\rrbracket^+)=h(\llbracket b\rrbracket^+)$, and so by $h$'s injectivity, $\llbracket a\rrbracket^+=\llbracket b\rrbracket^+$. Thus $a\equiv_R^+b$ as desired. 
\end{proof}

The next two chapters will be dedicated to showing the braid monoid satisfies the cancellative Ore conditions. Once we have done so, we will be able to embed the braid monoid (no scare quotes now!) into the braid group. Then Theorem \ref{thm:group_equiv_monoid_equiv_ore} can be used to state that if two positive braid words are equivalent, they are also equivalent in the braid monoid. 

\begin{theorem}\label{thm:braid_group_to_monoid_equiv_of_pos} \leanlink{https://github.com/hannahfechtner/braids_better/blob/098309403cda787eaf1fa954e07ff0ca3b2927ee/BraidProject/BraidLocalization.lean#L155-L164}
    For positive $a, b \in BW$ such that $a\equiv_{R_B}b$, we have $a\equiv^+_{R_{B^+}}b$.
\end{theorem}

\noindent That is, there exists a proof of their equivalence in which undercrossings never appear.\cprotect\footnote{This is not a simple application of Theorem \ref{thm:group_equiv_monoid_equiv_ore}. Our formalizations of the relations $R_B$ and $R_{B^+}$ do not match correctly -- we do not know that \lstinline{free_group_set_of_function} applied to our implementation of $R_{B^+}$ (recall, this is the simplified $R^S_{B^+}$) equals our implementation of $R_B$. That is, we do not have \lstinline{free_group_set_of_function (braid_monoid_rels_inf) = braidRelationInf}. As an instance of an Artin-Tits relation, $R_B$ (\lstinline{braidRelationInf}) is symmetric and contains the element $1$. Our direct implementation of $R^s_{B^+}$ (\lstinline{braid_monoid_rels_inf}) is not symmetric and contains no occurrence of reflexivity (which would be used to generate the identity element for $R_B$). Thus, we must implement a non-simplified monoid relation $R_M$ as follows.

\begin{definition}
    We implement the (infinite) braid monoid relation $R_{M}$ by the following :
    \begin{lstlisting}
inductive braid_monoid_rels_inf_one_refl_symm : FreeMonoid ℕ → FreeMonoid ℕ → Prop
  | adjacent (i j : ℕ) (h : i.dist j = 1) :
      braid_monoid_rels_inf_one_refl_symm (of i * of j * of i) (of j * of i * of j)
  | separated (i j : ℕ) (h : i.dist j ≥ 2) :
      braid_monoid_rels_inf_one_refl_symm (of i * of j) (of j * of i)
  | basic : braid_monoid_rels_inf_one_refl_symm 1 1
    \end{lstlisting}
Note that the \lstinline{basic} case could be replaced by any instance of reflexivity. 
\end{definition}
\begin{lemma}
    We have our desired condition that 
    
    \lstinline{free_group_set_of_function (braid_monoid_rels_inf_one_refl_symm) = braidRelationInf}.
\end{lemma}
\begin{lemma}
    The presented monoids with generating relations $R_M$ and $R^S_{B^+}$ are isomorphic.
\end{lemma}

\noindent Thus, we may return to Theorem \ref{thm:group_equiv_monoid_equiv_ore} and safely apply it in this case.}

The following (and final) theorem is not necessary for our project, but is a nice consequence of the above discussion :
 
\begin{theorem}\label{thm:fraction_form} \leanlink{google.com}
    Let $\langle S\, | R\rangle$ be a left-cancellative monoid satisfying the Ore conditions. Let $c$ be an element of $\langle S\, | R\rangle$.  Then for some $a, b\in \langle S \, | R \rangle ^+$, $c=(\iota(b))^{-1}\iota(a)$. 
\end{theorem}
\begin{proof}
    By surjectivity of $h_2$, we know that there is some $x\in\PML(S, R)$ such that $h_2(x)=c$. By definition of $\PML$, $x=\frac{r}{s}$ for $r, s\in \langle S\,|R\rangle^+$. We claim that $c=(\iota(s))^{-1}\iota(r)$. 
    
\begin{align*}
c  &= h_2\left(\frac{r}{s}\right) \\
&= h_2\left(\frac{1\cdot r}{1 \cdot s}\right) \\
&= h_2\left(\frac{1}{s} \sbullet \frac{r}{1}\right)
  && \text{letting $a=1$ and $b=1$ in definition \ref{def:ore_localization} of $\sbullet$} \\
&= h_2\left(\left(\frac{s}{1}\right)^{-1} \sbullet \frac{r}{1}\right)
  && \text{definition of inverse in \ref{thm:ore_group}} \\
&= \left(h_2\left(\frac{s}{1}\right)\right)^{-1} \cdot h_2\left(\frac{r}{1}\right)
  && \text{$h_2$ is a homomorphism} \\
&= \left(h_2\left(h(s)\right)\right)^{-1} \cdot h_2\left(h(r)\right)
  && \text{definition of $h$} \\
&= \left(\iota(s)\right)^{-1} \iota(r)
  && \text{definition of $\iota$}
\end{align*}
\end{proof}

Thus, specialized to braids, we see that any braid is equivalent to one in positive-negative form (that is, all overcrossings appearing before any undercrossings).

\section{Common Multiples}
In this chapter, we show the first of the two required conditions of Ore localization : common multiples. The proof for the common-multiples property follows Dehornoy exactly.\footnote{Dehornoy distinguishes between the braid monoid \textit{words} $\underline{\sigma}_{i,j}$ and $\underline{\Delta}_n$ and the braid monoid \textit{elements} $\sigma_{i,j}$ and $\Delta_n$ respectively. Here, for simplicity's sake, we speak only of the braid monoid words and omit their underlines.} Here we give a brief sketch -- if no proof is alluded to, it is by an inductive argument. 

\begin{definition} \leanlink{https://github.com/hannahfechtner/braids_better/blob/098309403cda787eaf1fa954e07ff0ca3b2927ee/BraidProject/AcrossStrands.lean#L121-L123}
    The $\sigma_{i,j}$ braid monoid word represents the $i^{\text{th}}$ strand moving to position $j$.  Given $i,j\in \mathbb{N}$, we define
\[
{\sigma}_{i,j} =
\begin{cases}
\varepsilon & \text{if } i = j,\\[6pt]
\sigma_i \cdot \sigma_{i+1,j} & \text{if } i < j,\\[6pt]
\sigma_{i-1} \cdot \sigma_{i-1,j} & \text{if } i > j.
\end{cases}
\]
\end{definition}
\noindent For example, $\sigma_{6,2}=\sigma_5\sigma_4\sigma_3\sigma_2$. As a physical braid,\footnote{We can only make this geometric connection once the localization procedure has been completed and we have an embedding of the braid monoid into the braid group} we see the strand in position 6 (remember, we are 0-indexing) moving to position 2 :
\begin{center}
\begin{tikzpicture}
\pic[
braid/.cd,
every strand/.style={ultra thick},
gap=0.2,
strand 3/.style={Mulberry},
strand 4/.style={Periwinkle},
strand 5/.style={TealBlue},
strand 6/.style={YellowGreen},
strand 7/.style={Rhodamine},
] {braid={s_6 s_5 s_4 s_3}};
\end{tikzpicture}
\end{center}

\begin{definition} \leanlink{https://github.com/hannahfechtner/braids_better/blob/098309403cda787eaf1fa954e07ff0ca3b2927ee/BraidProject/FlipBraid.lean#L66-L69}
    The $\Delta_{n}$ braid monoid word represents a half twist of the first $n+1$ strands in the braid.  We define
\begin{align*}
{\Delta}_0 &= \varepsilon \\
\Delta_{n+1} &= \sigma_{0,n+1}\;\Delta_n
\end{align*}

\end{definition}
\noindent For example, $\Delta_4=\sigma_0\sigma_1\sigma_2\sigma_3\sigma_0\sigma_1\sigma_2\sigma_0\sigma_1\sigma_0$. 

\begin{center}
\begin{tikzpicture}
\pic[
braid/.cd,
every strand/.style={ultra thick},
gap=0.2,
strand 1/.style={Mulberry},
strand 2/.style={Periwinkle},
strand 3/.style={TealBlue},
strand 4/.style={YellowGreen},
strand 5/.style={Rhodamine},
] {braid={s_1 s_2 s_3 s_4 s_1 s_2 s_3 s_1 s_2 s_1}};
\end{tikzpicture}
\end{center}

We next note that any generator can be completed into a delta braid of sufficiently large index :
\begin{theorem} \label{thm:partial_generator} \leanlink{https://github.com/hannahfechtner/braids_better/blob/098309403cda787eaf1fa954e07ff0ca3b2927ee/BraidProject/FlipBraid.lean#L182-L190}
    For $n,i\in \mathbb{N}$ such that $0<i<n$, there exists a value in $BW^+$, which we will denote $\partial_n(i)$, such that $\sigma_i\partial_n(i)\equiv^+{\Delta}_n$.
\end{theorem}

\noindent Thus, for any $k\in \mathbb{N}$ with $k>i$ and $k>j$, $\Delta_{k}$ a common multiple for $\sigma_i$ and $\sigma_j$.

\begin{definition}
    The flip function $f_n$ is defined on generators of $BW^+$ by $f_n(\sigma_i)=\sigma_{n-1-i}$ (where $i<n$). We may lift $f_n$ to a word in $BW^+$ (bounded by $n$) by applying $f_n$ to each of the word's generators. By abuse of notation, we denote this also by $f_n$.
\end{definition}
\noindent The flip function helps us move generators, and hence entire words in $BW^+$, past a half-twist.

\begin{lemma} \label{lem:move_flip} \leanlink{https://github.com/hannahfechtner/braids_better/blob/098309403cda787eaf1fa954e07ff0ca3b2927ee/BraidProject/FlipBraid.lean#L169-L180}
 For any $a\in BW^+$ strictly bounded above by $n$, $\;a{\Delta}_{n}\equiv^+{\Delta}_{n} f_n(a)$.
\end{lemma}

\noindent We will use the flip function in the witness to the following  theorem:
\begin{theorem} \label{thm:partial} \leanlink{https://github.com/hannahfechtner/braids_better/blob/098309403cda787eaf1fa954e07ff0ca3b2927ee/BraidProject/FlipBraid.lean#L267-L301}
    Let $a\in BW^+$ have length at most $\ell$. Let the set of  indices of the generators in $a$ be strictly bounded above by $n$. Then there exists a word $p_{\ell, n}(a)$ such that $a \; p_{\ell, n}(a)\equiv^+({\Delta}_n)^l$.
\end{theorem}

\begin{proof}(Sketch)
    Let $a=\sigma_{a_1}\sigma_{a_2}...\sigma_{a_{\ell-1}}\sigma_{a_\ell}$. Define $$p_{\ell, n}(a) = p_{\ell, n}(\sigma_{a_1}\sigma_{a_2}...\sigma_{a_{\ell-1}}\sigma_{a_\ell}) =\partial_n(a_\ell)\; f_n(\partial_n(a_{\ell-1})) \; \partial_n(a_{\ell-2}) \; f_n(\partial_n(a_{\ell-3}))...\footnote{The ellipses do not continue on indefinitely. Writing it in this manner avoids casework on the parity of $l$, which determines whether the flip function is applied to $\partial_i(a_0))$}$$

    For example, when $\ell =3$, 
    $$p_{3,n}(a) = p_{3,n}(\sigma_{a_1}\sigma_{a_2}\sigma_{a_3}) =\partial_n(a_3)\; f_n(\partial_n(a_2)) \; \partial_n(a_1)$$. 

    To show correctness of this particular witness, we calculate
    \begin{align*}
        a\; p_{3,n}(a) & = \sigma_{a_1}\sigma_{a_2}\sigma_{a_3} \partial_n(a_3)\; f_n(\partial_n(a_2)) \; \partial_n(a_1) \\
        &\equiv^+ \sigma_{a_1}\sigma_{a_2}\Delta_n\; f_n(\partial_n(a_2)) \; \partial_n(a_1) &\text{by Lemma \ref{thm:partial_generator}} \\
        &\equiv^+ \sigma_{a_1}\sigma_{a_2}\;\partial_n(a_2)\Delta_n \; \partial_n(a_1) &\text{by Lemma \ref{lem:move_flip}} \\
        &\equiv^+ \sigma_{a_1}\Delta_n\Delta_n \, \partial_n(a_1) &\text{by Lemma \ref{thm:partial_generator}} \\
        & = \;\ \sigma_{a_1}\Delta_n\Delta_n \, f_n(f_n(\partial_n(a_1))) &f_n\text{ is an involution} \\
        &\equiv^+ \sigma_{a_1}\, \Delta_i \,f_n(\partial_i(a_1)) \Delta_n &\text{by Lemma \ref{lem:move_flip}} \\
        &\equiv^+ \sigma_{a_1}\, \partial_n(a_1) \, \Delta_n \Delta_n &\text{by Lemma \ref{lem:move_flip}} \\
        &\equiv^+ \Delta_n\Delta_n \Delta_n &\text{by Lemma \ref{thm:partial_generator}} \\
        &=(\Delta_n)^3
    \end{align*}

    The full proof of correctness proceeds by induction on $\ell$. 
\end{proof}

\begin{theorem} \label{thm:right_common_multiples}(Right Common Multiples) \leanlink{https://github.com/hannahfechtner/braids_better/blob/098309403cda787eaf1fa954e07ff0ca3b2927ee/BraidProject/CommonMultiples.lean#L10-L19}
    For all $a$ and $b$ in $BW^+$, there exist $c$ and $d$ in $BW^+$ such that $ac\equiv^+bd$. 
\end{theorem}
\begin{proof}
Set $\ell$ as the maximum of the lengths of $a$ and $b$. Fix some $n\in\mathbb{N}$ strictly greater than the index of the largest generator occurring in either $a$ or $b$. Then by Theorem \ref{thm:partial} we have 
$$a\; p_{\ell,n}(a) \equiv^+ (\Delta_n)^\ell\equiv^+b\; p_{\ell,n}(b)$$ 
so $p_{\ell,n}(a)$ and $p_{\ell,n}(b)$ are the desired witnesses.
\end{proof}

While not needed for Ore localization, showing left common multiples is simple, thanks to the reversing function $r$ introduced in Lemma \ref{lem:reverse_equiv}.

\begin{theorem} \label{thm:left_common_multiples} (Left Common Multiples) \leanlink{https://github.com/hannahfechtner/braids_better/blob/098309403cda787eaf1fa954e07ff0ca3b2927ee/BraidProject/CommonMultiples.lean#L28-L34}
    For all $a$ and $b$ in $BW^+$, there exist $c$ and $d$ in $BW^+$ such that $ca\equiv^+db$. 
\end{theorem}
\begin{proof}
Apply Theorem \ref{thm:right_common_multiples} to $r(a)$ and $r(b)$ to obtain $c'$ and $d'$ with $r(a)\cdot c'\equiv^+r(b)\cdot d'$. Then use $r(c')$ and $r(d')$ as our witnesses. We aim to show that $r(c') \cdot a \equiv ^+ r(d')\cdot b $. By Lemma \ref{lem:reverse_equiv}, it suffices to show that $r(r(c')\cdot a) \equiv ^+ r(r(d')\cdot b )$. But by definition of $r$, this simplifies to $r(a)\;c'\equiv^+r(b)\;d'$, which we have. 

\end{proof}

We have thus shown the first of two Ore conditions for the braid monoid. Cancellativity is a bit trickier; it is proved using an alternative, diagrammatic representation of equivalence in the braid monoid, called a grid. This structure will be discussed in the following chapter.

\section{Grids} 

We next prove the remaining condition for Ore localization -- cancellativity of the braid monoid. To do so, following Dehornoy, we model braid monoid equivalence using rectangular grid diagrams\footnote{Grids are a type of van Kampen diagram, planar diagrams which exhibit the fact that a given word in a presented group is equivalent to the identity element. Van Kampen diagrams were first introduced by their namesake in 1933 \cite{van_kampen_lemmas_1933}. A full exposition can be found in the classic text by Lyndon and Schupp \cite{lyndon_combinatorial_2001}, which Nicholas Touikan summarizes in a visual, accessible manner in his online textbook \cite{touikan_van_2024}.}. Later, these grids will be given further interpretation and reused in the proofs of both correctness and termination of the subword reversing algorithm.


\subsection{Definition and Basic Properties}
Here is an example of a grid, with source $(\sigma_2\sigma_1, \sigma_0\sigma_1)$ (seen down the left side and across the top side) and target $(\sigma_0\sigma_1\sigma_2, \sigma_2\sigma_1\sigma_0)$ (seen across the bottom side and down the right side). We will see that the labels on all paths from top-left to bottom-right are equivalent in the braid monoid.

\begin{center}
\begin{tikzcd}
{} \arrow[dd, "\sigma_2"'] \arrow[rr, "\sigma_0"]  &                           & {} \arrow[rr, "\sigma_1"] \arrow[dd, "\sigma_2"]   &                                                    & {} \arrow[d, "\sigma_2"]    \\
                                                   &                           &                                                    &                                                    & {} \arrow[d, "\sigma_1"]    \\
{} \arrow[rr, "\sigma_0"'] \arrow[dd, "\sigma_1"'] &                           & {} \arrow[r, "\sigma_1"'] \arrow[d, "\sigma_1"]    & {} \arrow[r, "\sigma_2"'] \arrow[d, "\varepsilon"] & {} \arrow[d, "\varepsilon"] \\
                                                   &                           & {} \arrow[d, "\sigma_0"] \arrow[r, "\varepsilon"'] & {} \arrow[r, "\sigma_2"'] \arrow[d, "\sigma_0"]    & {} \arrow[d, "\sigma_0"]    \\
{} \arrow[r, "\sigma_0"']                          & {} \arrow[r, "\sigma_1"'] & {} \arrow[r, "\varepsilon"']                       & {} \arrow[r, "\sigma_2"']                          & {}                         
\end{tikzcd}
\end{center}

Intuitively, a grid is defined inductively as either a cell, or two grids abutted horizontally or vertically so that the relevant edges match. We give the data for a grid as a 4-tuple of braid monoid words, listing out the left, top, bottom, and right sides. We will later see that this uniquely determines the inner edges of the grid. We begin with the definition of cells:

\begin{definition} \leanlink{https://github.com/hannahfechtner/braids_better/blob/098309403cda787eaf1fa954e07ff0ca3b2927ee/BraidProject/Grid/Basic.lean#L12-L17}
The set of \textbf{cells} $C$ is the smallest subset of $(BW^+)^4$ such that:
\begin{itemize}
    \item $\langle \varepsilon, \varepsilon,\varepsilon,\varepsilon\rangle \in C$. Visually, \begin{tikzcd}
{} \arrow[r, "\varepsilon"] \arrow[d, "\varepsilon"'] & {} \arrow[d, "\varepsilon"] \\
{} \arrow[r, "\varepsilon"']                   & {}               
\end{tikzcd} 
    \item $\langle \varepsilon, \sigma _i,\sigma_i,\varepsilon\rangle \in C$ for any $i\in\mathbb{N}$. Visually,
\begin{tikzcd}
{} \arrow[r, "\sigma_{i}"] \arrow[d, "\varepsilon"'] & {} \arrow[d, "\varepsilon"] \\
{} \arrow[r, "\sigma_{i}"']                   & {}               
\end{tikzcd}

\item $\langle \sigma_i,\varepsilon,\varepsilon,\sigma_i\rangle\in C$ for any $i\in \mathbb{N}$. Visually, 
\begin{tikzcd}
{} \arrow[r, "\varepsilon"] \arrow[d, "\sigma_{i}"'] & {} \arrow[d, "\sigma_{i}"] \\
{} \arrow[r, "\varepsilon"']                   & {}               
\end{tikzcd}

\item $\langle\sigma_i,\sigma_i,\varepsilon,\varepsilon\rangle\in C$ for any $i\in \mathbb{N}$. Visually,
\begin{tikzcd}
{} \arrow[r, "\sigma_{i}"] \arrow[d, "\sigma_{i}"'] & {} \arrow[d, "\varepsilon"] \\
{} \arrow[r, "\varepsilon"']                   & {}               
\end{tikzcd}
\item For $i,j\in \mathbb{N}$ with $|i-j|=1$, $\langle\sigma_j,\sigma_i,\sigma_i\sigma_j,\sigma_j\sigma_i\rangle\in C$. Visually,
\begin{tikzcd}
{} \arrow[rr, "\sigma_{i}"] \arrow[dd, "\sigma_{j}"'] &              & {} \arrow[d, "\sigma_{j}"] \\
                         &              & {} \arrow[d, "\sigma_{i}"] \\
{} \arrow[r, "\sigma_{i}"']             & {} \arrow[r, "\sigma_{j}"'] & {}          
\end{tikzcd}
\item For $i,j\in\mathbb{N}$ with $|i-j|>1$, $\langle\sigma_j,\sigma_i,\sigma_i,\sigma_j\rangle\in C$. Visually,
\begin{tikzcd}
{} \arrow[r, "\sigma_{i}"] \arrow[d, "\sigma_j"'] & {} \arrow[d, "\sigma_j"] \\
{} \arrow[r, "\sigma_{i}"']                   & {}               
\end{tikzcd}
\end{itemize}
\end{definition}






Often we simplify our diagrams combining adjacent arrows and labels. This does not change the cell or structure; it is merely notation. Thus, for the second-to-last case of the above definition, we could draw:
\begin{center}
\begin{tikzcd}
{} \arrow[r, "\sigma_i"] \arrow[d, "\sigma_j"'] & {} \arrow[d, "\sigma_j\sigma_i"] \\
{} \arrow[r, "\sigma_i\sigma_j"']               & {}                              
\end{tikzcd}
\end{center}

\begin{definition} \leanlink{https://github.com/hannahfechtner/braids_better/blob/098309403cda787eaf1fa954e07ff0ca3b2927ee/BraidProject/Grid/Basic.lean#L9-L19}
The set of \textbf{grids} $G$ is defined inductively as follows
\begin{itemize}
\item Any cell is a grid; $C\subseteq G$

\item Horizontal abutting: If $\langle a, b, c,d\rangle, \langle d,e,f,g\rangle\in G$, then $\langle a,be,cf,g\rangle\in G$.

Visually, if \begin{tikzcd}
{} \arrow[r, "b"] \arrow[d, "a"'] & {} \arrow[d, "d"] \\
{} \arrow[r, "c"']                & {}               
\end{tikzcd} and 
\begin{tikzcd}
{} \arrow[r, "e"] \arrow[d, "d"'] & {} \arrow[d, "g"] \\
{} \arrow[r, "f"']                & {}               
\end{tikzcd}
are grids, then so too is
\begin{tikzcd}
{} \arrow[r, "b"] \arrow[d, "a"'] & {} \arrow[d, "d"] \arrow[r, "e"] & {} \arrow[d, "g"] \\
{} \arrow[r, "c"']                & {} \arrow[r, "f"']               & {}               
\end{tikzcd}\item
Vertical abutting: If $\langle a, b, c,d\rangle, \langle e,c,f,g\rangle\in G$, then $\langle ae,b,f,dg\rangle\in G$.

Visually, if \begin{tikzcd}
{} \arrow[r, "b"] \arrow[d, "a"'] & {} \arrow[d, "d"] \\
{} \arrow[r, "c"']                & {}               
\end{tikzcd} and 
\begin{tikzcd}
{} \arrow[r, "c"] \arrow[d, "e"'] & {} \arrow[d, "g"] \\
{} \arrow[r, "f"']                & {}               
\end{tikzcd}
are grids, then so too is
\begin{tikzcd}
{} \arrow[r, "b"] \arrow[d, "a"']  & {} \arrow[d, "d"] \\
{} \arrow[r, "c"'] \arrow[d, "e"'] & {} \arrow[d, "g"] \\
{} \arrow[r, "f"']                 & {}               
\end{tikzcd}

\end{itemize}
\end{definition}
We may also draw a grid omitting its inner arrows, simply showing the four sides and the words labelling them. Later we will prove a uniqueness property for grids (Theorem \ref{thm:grid_uniqueness}), so this representation does not omit unrecoverable information. The last-drawn grid would thus be

\begin{center}
\begin{tikzcd}
{} \arrow[r, "b"] \arrow[d, "ae"'] & {} \arrow[d, "dg"] \\
{} \arrow[r, "f"']               & {}               
\end{tikzcd}
\end{center}

\begin{definition}
    Given a grid $\langle a,b,c,d\rangle\in G$, its \textbf{source} (later also called spine) is $(a,b)$ and its \textbf{target} is $(c,d)$
\end{definition}
\noindent For example, the above grid has source $(ae, b)$ and target $(f, dg)$. 

\begin{theorem}\label{thm:grid_equiv} \leanlink{https://github.com/hannahfechtner/braids_better/blob/098309403cda787eaf1fa954e07ff0ca3b2927ee/BraidProject/Grid/Basic.lean#L58-L77}
If $\langle a,b,c,d\rangle\in G$, $ac\equiv^+bd$.
\end{theorem}

\begin{proof}
    By structural induction on the grid.
\end{proof}

\begin{definition} \label{def:diag_length}
    Given a grid $\langle a,b,c,d\rangle\in G$, the \textbf{length of its diagonal} is the length of $ac$. 
\end{definition}
\begin{remark} \leanlink{https://github.com/hannahfechtner/braids_better/blob/098309403cda787eaf1fa954e07ff0ca3b2927ee/BraidProject/Grid/Basic.lean#L85-L89}
    The length of a grid's diagonal is also equal to $|bd|$, as $ac\equiv^+bd$ by Theorem \ref{thm:grid_equiv} and equivalent words have the same length by Definition \ref{def:braid_length}. 
\end{remark}

We next amass simple properties of grids that will be helpful in proving more complex theorems. All proofs are via structural induction on the given grid.
\begin{lemma} \label{lem:grid_swap}\leanlink{https://github.com/hannahfechtner/braids_better/blob/098309403cda787eaf1fa954e07ff0ca3b2927ee/BraidProject/Grid/Basic.lean#L23-L34}
If $\langle a, b, c, d\rangle \in G$, then $\langle b, a, d, c\rangle \in G$. 

Visually, given \begin {tikzcd}{} \arrow[r, "b"] \arrow[d, "a"'] & {} \arrow[d, "d"] \\
{} \arrow[r, "c"']                & {}               
\end{tikzcd}, we also have 
\begin {tikzcd}{} \arrow[r, "a"] \arrow[d, "b"'] & {} \arrow[d, "c"] \\
{} \arrow[r, "d"']                & {}               
\end{tikzcd}
\end{lemma}
\begin{lemma}\label{lem:grid_splittability}(Splitting)\mbox{} \leanlink{https://github.com/hannahfechtner/braids_better/blob/098309403cda787eaf1fa954e07ff0ca3b2927ee/BraidProject/Grid/Basic.lean#L91-L165}
Given $\langle a_1a_2,b,c,d\rangle \in G$, there exist grids $\langle a_1,b,m,d_1 \rangle$ and $\langle a_2,m,c,d_2\rangle$ with $d=d_1d_2$. 

That is,
\begin{tikzcd}
{} \arrow[d, "a_1"'] \arrow[r, "b"] & {} \arrow[dd, "d"] \\
{} \arrow[d, "a_2"']                &                    \\
{} \arrow[r, "c"']                  & {}                
\end{tikzcd}
can be split horizontally into two subgrids as so:
\begin{tikzcd}
{} \arrow[d, "a_1"'] \arrow[r, "b"]  & {} \arrow[d, "d_1"] \\
{} \arrow[d, "a_2"'] \arrow[r, "m"'] & {} \arrow[d, "d_2"] \\
{} \arrow[r, "c"']                   & {}                 
\end{tikzcd}

\noindent Similarly, given $\langle a,b_1b_2,c,d\rangle \in G$, there exist grids $\langle a,b_1,c_1,m \rangle$ and $\langle m,b_2,c_2,d\rangle$ with $c=c_1c_2$. 

That is, \begin{tikzcd}
{} \arrow[r, "b_1"] \arrow[d, "a"'] & {} \arrow[r, "b_2"] & {} \arrow[d, "d"] \\
{} \arrow[rr, "c"']                 &                     & {}               
\end{tikzcd}
can be split vertically into two subgrids so:
\begin{tikzcd}
{} \arrow[r, "b_1"] \arrow[d, "a"'] & {} \arrow[r, "b_2"] \arrow[d, "m"] & {} \arrow[d, "d"] \\
{} \arrow[r, "c_1"']                & {} \arrow[r, "c_2"']               & {}               
\end{tikzcd} Note that Lemma \ref{lem:grid_swap} allows vertical splitting to be derived from horizontal splitting.
\end{lemma}
\begin{lemma}[Determinative Spines] \leanlink{https://github.com/hannahfechtner/braids_better/blob/098309403cda787eaf1fa954e07ff0ca3b2927ee/BraidProject/Grid/DeterminativeSpine.lean#L11-L207} \mbox{}
For braid monoid words $a,b,c,d$ and natural numbers $i,j$,
\begin{itemize}
    \item If $\langle \varepsilon, b, c, d\rangle\in G$, then $c=b$ and $d=\varepsilon$.
    \item If $\langle a, \varepsilon, c, d\rangle\in G$, then $c=\varepsilon$ and $d=a$.
    \item If $\langle a, a, c, d\rangle\in G$, then $c=\varepsilon$ and $d=\varepsilon$.
    \item If $\langle \sigma_i, \sigma_j, c, d\rangle\in G$ with $|i-j|=1$, then $c=\sigma_j\sigma_i$ and $d=\sigma_i\sigma_j$.
    \item If $\langle \sigma_i, \sigma_j, c, d\rangle\in G$ with $|i-j|>1$, then $c=\sigma_j$ and $d=\sigma_i$.
\end{itemize}
\end{lemma}

\begin{theorem} (Uniqueness of sinks) \label{thm:grid_uniqueness}
\leanlink{https://github.com/hannahfechtner/braids_better/blob/098309403cda787eaf1fa954e07ff0ca3b2927ee/BraidProject/Grid/DeterminativeSpine.lean#L52-L67}
    Any two grids with the same source must have the same sink.
\end{theorem}

\begin{proof} (Sketch)
    Consider two grids $g_1=\langle a, b, c, d\rangle$ and $g_2=\langle a, b, c', d'\rangle$. We show $c=c'$ and $d=d'$ by induction on the $g_1$. The base cases are shown by determinative spines, and the inductive cases split $g_2$ and apply the inductive hypothesis to the resulting pieces.
\end{proof}
\subsection{Stability}

We define a property of grids called stability which will be essential to build grids from monoid equivalence, and thus to use grid properties to show cancellativity in the braid monoid.

\begin{definition}[Stability] \leanlink{https://github.com/hannahfechtner/braids_better/blob/098309403cda787eaf1fa954e07ff0ca3b2927ee/BraidProject/Grid/Stability.lean#L8-L10}
A grid $\langle a, b, c, d\rangle$ is called \textbf{stable} if for all $a', b'\in BW^+$  such that $a\equiv^+a'$ and $b\equiv^+b'$, there exist $c', d'\in BW^+$ such that we have a grid $\langle a', b', c', d'\rangle$ with $c\equiv^+c'$ and $d\equiv^+d'$
\end{definition}

The Lean proof that all grids are stable is done by a triple induction, first on the length of the grid's diagonal. We then induct structurally on the equivalence between $a$ and $a'$ and then that of $b$ and $b'$. Dehornoy, instead of these two structural inductions, defines a length metric on derivations witnessing equivalence and inducts on the sum of the lengths of the respective derivations associated with $a\equiv a'$ and $b \equiv b'$. At this stage of our Lean proof, we have kept a very simple definition of a derivation that is \lstinline{Prop}-valued in Lean (it carries only the fact that a derivation exists, not the data of the derivation itself) and thus cannot recurse on the derivation's construction to calculate its length. Instead, we induct structurally, one derivation at a time. This does not increase the total number of cases significantly. We separate out a large portion of the proof into the following lemma:

\begin{lemma}\label{stable_braid_rel} \leanlink{https://github.com/hannahfechtner/braids_better/blob/098309403cda787eaf1fa954e07ff0ca3b2927ee/BraidProject/Grid/Stability.lean#L188-L193}
    Let $(w, x)\in R_{B^+}$ and $k\in \mathbb{N}$. All grids having source $(\sigma_k, w)$ or $(\sigma_k, x)$ are stable. 
\end{lemma}
\begin{proof} (sketch)
    We discuss the proof for $(\sigma_k, w)$; the proof for $(\sigma_k, x)$ proceeds nearly identically. Consider an arbitrary grid $g = \langle \sigma_k,w,c,d\rangle$. We proceed by casework on $w$, and then on the distances between $k$ and each of the indices appearing in $w$. For example, let $w$ be a commutative relation of the form $\sigma_i\sigma_j$ with $|i-j|\geq 2$, with $k$ such that $|i-k|=1$ and $|j-k|=1$. For clarity in the following diagram, we instantiate this with $i=2, \, j=4,\, k=3$.  We know there exists a grid $g'=\langle\sigma_3,\sigma_2\sigma_3,\sigma_2\sigma_3\sigma_4\sigma_3\sigma_2, \sigma_3\sigma_4\sigma_2\sigma_3\rangle$, as evidenced by the diagram below which constructs it from single cells : 
    
\begin{center}
\begin{tikzcd}
{} \arrow[rr, "\sigma_2"] \arrow[dddd, "\sigma_3"'] &                           & {} \arrow[dd, "\sigma_3"] \arrow[rrr, "\sigma_4"] &                                                  &                           & {} \arrow[d, "\sigma_3"] \\
                                                    &                           &                                                   &                                                  &                           & {} \arrow[d, "\sigma_4"] \\
                                                    &                           & {} \arrow[dd, "\sigma_2"] \arrow[r, "\sigma_4"]   & {} \arrow[rr, "\sigma_3"] \arrow[dd, "\sigma_2"] &                           & {} \arrow[d, "\sigma_2"] \\
                                                    &                           &                                                   &                                                  &                           & {} \arrow[d, "\sigma_3"] \\
{} \arrow[r, "\sigma_2"']                           & {} \arrow[r, "\sigma_3"'] & {} \arrow[r, "\sigma_4"']                         & {} \arrow[r, "\sigma_3"']                        & {} \arrow[r, "\sigma_2"'] & {}                      
\end{tikzcd}
\end{center}

Thus, by uniqueness of sinks (Theorem \ref{thm:grid_uniqueness}) applied to $g$ and $g'$, we must have $c = \sigma_2\sigma_3\sigma_4\sigma_3\sigma_2$ and $d = \sigma_3\sigma_4\sigma_2\sigma_3$.
We consider $a'$ and $b'$ such that $\sigma_3\equiv^+a'$ and $\sigma_2\sigma_4\equiv ^+b'$. By Lemma \ref{lem:generator_equiv}, $a'=\sigma_3$. By Lemma \ref{lem:equiv_relation}, either $b'=\sigma_2\sigma_4$ or $b'=\sigma_4\sigma_2$. In the first case, we take $c'=c$ and $d'=d$. We use the original grid $\langle \sigma_3,\sigma_2\sigma_4,c,d\rangle$ and quickly finish the proof as the required equivalences are true by reflexivity. In the second case, we take $c'=\sigma_4\sigma_3\sigma_2\sigma_3\sigma_4$ and $d'=\sigma_3\sigma_2\sigma_4\sigma_3$ and the following grid:
\begin{center}
\begin{tikzcd}
{} \arrow[rr, "\sigma_4"] \arrow[dddd, "\sigma_3"'] &                           & {} \arrow[dd, "\sigma_3"] \arrow[rrr, "\sigma_2"] &                                                  &                           & {} \arrow[d, "\sigma_3"] \\
                                                    &                           &                                                   &                                                  &                           & {} \arrow[d, "\sigma_2"] \\
                                                    &                           & {} \arrow[dd, "\sigma_4"] \arrow[r, "\sigma_2"]   & {} \arrow[rr, "\sigma_3"] \arrow[dd, "\sigma_4"] &                           & {} \arrow[d, "\sigma_4"] \\
                                                    &                           &                                                   &                                                  &                           & {} \arrow[d, "\sigma_3"] \\
{} \arrow[r, "\sigma_4"']                           & {} \arrow[r, "\sigma_3"'] & {} \arrow[r, "\sigma_2"']                         & {} \arrow[r, "\sigma_3"']                        & {} \arrow[r, "\sigma_4"'] & {}                      
\end{tikzcd}
\end{center}
A quick calculation verifies
\begin{align*}
    c&=\sigma_2\sigma_3\sigma_4\sigma_3\sigma_2 \\
    &\equiv^+ \sigma_2\sigma_4\sigma_3\sigma_4\sigma_2 \\
    &\equiv^+ \sigma_4\sigma_2\sigma_3\sigma_4\sigma_2 \\
    &\equiv^+ \sigma_4\sigma_2\sigma_3\sigma_2\sigma_4 \\
    &\equiv^+ \sigma_4\sigma_3\sigma_2\sigma_3\sigma_4\\
    & = c'
\end{align*}

\noindent A similar calculation verifies $d\equiv^+d'$ and completes the proof.
\end{proof}
\begin{theorem} \label{thm:stability} \leanlink{https://github.com/hannahfechtner/braids_better/blob/098309403cda787eaf1fa954e07ff0ca3b2927ee/BraidProject/Grid/Stability.lean#L293}
All grids are stable.
\end{theorem}
For the proof of stability, we present one illustrative case, drawn from the inductive case of the first induction (on diagonal length). Our inductive hypothesis says that for some $l$, we know that all grids with length strictly less than $l$ are stable. We consider an arbitrary grid $\langle a, b, c, d\rangle$ with diagonal length $l$.

We select the case where $a\equiv^+a$ is true by reflexivity (so $a=a'$), and where we are in the one-step reduction case of the induction on $b\equiv^+b'$. Thus, we assume that $b=b_1wb_2$ and $b'=b_1xb_2$, where $(w, x)\in R_{B^+}$. Splitting the grid vertically, we obtain
\[
\begin{tikzcd}
{} \arrow[r, "b_1"] \arrow[d, "a"']
& {} \arrow[r, "w"] \arrow[d, "d_1"]
& {} \arrow[r, "b_2"] \arrow[d, "f"]
& {} \arrow[d, "d"] \\
{} \arrow[r, "c_1"']
& {} \arrow[r, "e"']
& {} \arrow[r, "c_2"']
& {}
\end{tikzcd}
\]

 \noindent where $c=c_1ec_2$. We then consider the case where $d_1$ has length at least one. We write $d_1=\sigma_kd_2$, and obtain the following grid and label the subgrids as follows (using capital letters for readability): 
 \begin{center}
\begin{tikzcd}
{} \arrow[r, "b_1"] \arrow[dd, "a"'] 
& {} \arrow[d, "\sigma_k"'] \arrow[r, "w"]        
   \arrow[dr, phantom, "G_2"] 
& {} \arrow[d, "f_1"] \arrow[r, "b_2"] 
& {} \arrow[dd, "d"] \\
{} \arrow[r, phantom, "G_1"] 
& {} \arrow[d, "d_2"'] \arrow[r, "g"'] 
   \arrow[dr, phantom, "G_3"] 
& {} \arrow[d, "f_2"] \arrow[r, phantom, "G_4"] 
& {}                  \\
{} \arrow[r, "c_1"']                 
& {} \arrow[r, "e"']                             
& {} \arrow[r, "c_2"']                 
& {}
\end{tikzcd}
\end{center}

\noindent Note that $f=f_1f_2$. The length of $G_3$ is $|d_2e|$. Calculating, we see that we may apply the induction hypothesis to $G_3$ :
$$|d_2e|<|\sigma_kd_2e|\leq|b_1\sigma_kd_2ec_2|=|ac_1ec_2|=|ac|=l$$

\noindent $G_3$ is hence stable, and a similar argument shows that $G_4$ is as well.\footnote{Note that it is entirely possible that the diagonal lengths of $G_1$ and $G_2$ may equal $l$} By Lemma \ref{stable_braid_rel}, we know that $G_2$ is also stable.

Unfortunately, lemmas showing preservation of stability under appending are as difficult to prove as Theorem \ref{thm:stability} itself, as they require inductions on equivalences. Instead, we proceed manually, calculating the required grids and equivalences for this specific case. 

First, by stability of $G_2$, we have ${G_2}'$ with source $(\sigma_k, x)$ and target $g', {f_1}'$ where $g\equiv^+ g'$ and $f_1\equiv^+{f_1}'$. By stability of  $G_3$, we have ${G_3}'$ with source $(d_2, g')$ and target $(e', {f_2}')$ where $e\equiv^+e'$ and $f_2\equiv^+{f_2}'$. We then consider $G_4$. Note that $f_1f_2\equiv^+{f_1}'{f_2}'$. By stability of $G_4$, we have a grid ${G_4}'$ with source $({f_1}'{f_2}', b_2)$ and target $({c_2}', d')$ with ${c_2}'\equiv^+ c_2$ and $d'\equiv^+ d$. We build the following grid:

 \begin{center}
\begin{tikzcd}
{} \arrow[r, "b_1"] \arrow[dd, "a"'] 
& {} \arrow[d, "\sigma_k"'] \arrow[r, "x"]        
   \arrow[dr, phantom, "{G_2}'"] 
& {} \arrow[d, "{f_1}'"] \arrow[r, "b_2"] 
& {} \arrow[dd, "d'"] \\
{} \arrow[r, phantom, "G_1"] 
& {} \arrow[d, "d_2"'] \arrow[r, "g'"'] 
   \arrow[dr, phantom, "{G_3}'"] 
& {} \arrow[d, "{f_2}'"] \arrow[r, phantom, "{G_4}'"] 
& {}                  \\
{} \arrow[r, "c_1"']                 
& {} \arrow[r, "e'"']                             
& {} \arrow[r, "{{c_2}'}"']                 
& {}
\end{tikzcd}
\end{center}

We have just seen from the stability of $G_4$ that $d\equiv^+d'$. We must show that ${c_1}'e'{c_2}'\equiv^+c$. This is immediate from the facts given in the previous paragraph.

As we have seen in the above example, the proof can branch at the induction on the length and the inductions on the equivalences $a\equiv^+a'$ and $b\equiv^+b'$, and the lemma can branch at all the combinatorial options for the distance a single generator's index can lie in relation to the indices of elements of the braid monoid relations.  Luckily, we do not have worst-case branching behavior, which would lead to a bit more than 100 cases. Our lemma has fourteen cases, and the main proof has twelve, for a total of 26. This is approximately the same as Dehornoy's original proof. 

\begin{theorem} \label{thm:grid_one_one_iff_equiv} \leanlink{https://github.com/hannahfechtner/braids_better/blob/098309403cda787eaf1fa954e07ff0ca3b2927ee/BraidProject/Grid/Properties.lean#L8-L12}
    There exists a grid with source $(a,b)$ and target $(\varepsilon, \varepsilon)$ if and only if $a\equiv^+b$.
\end{theorem}
\begin{proof}

The forward direction is immediate by Theorem \ref{thm:grid_equiv}. For the backwards direction, assume that $a\equiv^+b$. It is quick to show by induction on $a$ that there is a grid $\langle a,a, \varepsilon, \varepsilon\rangle$. By Theorem \ref{thm:stability}, $\langle a,a, \varepsilon, \varepsilon\rangle$ is stable, so we know that there is a grid $\langle a,b,c',d'\rangle $ where $c'\equiv^+\varepsilon$ and $d'\equiv^+\varepsilon$. Hence $c'=\varepsilon$ and $d'=\varepsilon$ (length is preserved under braid monoid equivalence), so we have a grid $\langle a, b, \varepsilon, \varepsilon\rangle$ as desired. 
   
\end{proof}

\noindent We are now ready to prove cancellativity. 

\subsection{Cancellativity}
\begin{theorem} \leanlink{https://github.com/hannahfechtner/braids_better/blob/098309403cda787eaf1fa954e07ff0ca3b2927ee/BraidProject/Grid/Properties.lean#L14-L37}
    The braid monoid is left-cancellative; that is, if $ca\equiv^+cb$, then $a\equiv^+b$.
\end{theorem}
\begin{proof}
    Assume $ca=cb$. Then by Theorem \ref{thm:grid_one_one_iff_equiv}, there is a grid $\langle ca, cb, \varepsilon, \varepsilon\rangle$. By repeated applications of splitting and determinative spines, the grid can be written as such:
\begin{center}
    \begin{tikzcd}
{} \arrow[r, "c"] \arrow[d, "c"']            & {} \arrow[r, "b"] \arrow[d, "\varepsilon"] & {} \arrow[d, "\varepsilon"] \\
{} \arrow[r, "\varepsilon"'] \arrow[d, "a"'] & {} \arrow[r, "b"] \arrow[d, "a"']          & {} \arrow[d, "\varepsilon"] \\
{} \arrow[r, "\varepsilon"']                 & {} \arrow[r, "\varepsilon"']               & {}                         
\end{tikzcd}
\end{center}
The lower right subgrid  gives us $\langle a, b, \varepsilon, \varepsilon \rangle\in G$. Again by Theorem \ref{thm:grid_one_one_iff_equiv}, we know $a\equiv^+b$.
\end{proof}
\begin{theorem} \label{thm:right_cancellativity} \leanlink{https://github.com/hannahfechtner/braids_better/blob/098309403cda787eaf1fa954e07ff0ca3b2927ee/BraidProject/Grid/Properties.lean#L39-L42}
The braid monoid is right-cancellative; that is, if $ac\equiv^+bc$, then $a\equiv^+b$. This is proved from left-cancellativity using a trick with the reversing function similar to that of Theorem \ref{thm:left_common_multiples}.
\end{theorem}

Thus, we have shown the necessary criteria to prove Theorem \ref{thm:braid_group_to_monoid_equiv_of_pos}, and may conclude that for positive $a, b \in BW$ such that $a\equiv_{R_B}b$, we have $a\equiv^+_{R_{B^+}}b$.

\subsection{Grid Existence}
Beyond cancellativity, stability (via Theorem \ref{thm:grid_one_one_iff_equiv}) is also used to show the existence of grids.
\begin{theorem} \label{thm:grid_existence} \leanlink{https://github.com/hannahfechtner/braids_better/blob/098309403cda787eaf1fa954e07ff0ca3b2927ee/BraidProject/Grid/Properties.lean#L52-L61}
    For any $a,b\in BW^+$ there are $c,d\in BW^+$ such that there exists a grid with source $(a,b)$ and target $(c,d)$.
\end{theorem}

\begin{proof}
    Consider some $a,b\in BW^+$. We know there exist $c_1,d_1$ such that $ac_1\equiv^+bd_1$ by Theorem \ref{thm:right_common_multiples}. Then by Theorem \ref{thm:grid_one_one_iff_equiv} there is a grid $\langle ac_1, bd_1,\varepsilon, \varepsilon\rangle$. Then by repeatedly splitting that grid, there must exist $c,d,e,f\in BW^+$ and the following grid:
    \begin{center}
    \begin{tikzcd}
{} \arrow[r, "b"] \arrow[d, "a"']            & {} \arrow[r, "d_1"] \arrow[d, "d"] & {} \arrow[d, "\varepsilon"] \\
{} \arrow[r, "c"'] \arrow[d, "c_1"'] & {} \arrow[r, "f"] \arrow[d, "e"']          & {} \arrow[d, "\varepsilon"] \\
{} \arrow[r, "\varepsilon"']                 & {} \arrow[r, "\varepsilon"']               & {}                         
\end{tikzcd}
\end{center}
The top-left grid gives us $\langle a, b, c, d\rangle \in G$ as desired. 
\end{proof}

\section{Grids and Rewriting}\label{chap:grids_rw}

We have used grids fruitfully in the proof of the cancellativity in the braid monoid, as part of the localization procedure connecting the braid monoid to the braid group. This portion of the project is now complete.

We find continued use for grids in the next stage of the project, in which we discuss the reversing procedure, as informally discussed in the introduction.  Recall the reversing relation : 

\begin{definition}\label{def:reversing_rel} \leanlink{https://github.com/hannahfechtner/braids_better/blob/098309403cda787eaf1fa954e07ff0ca3b2927ee/BraidProject/Relations.lean#L8-L12}
    \noindent The reversing relation $R_{\mathsf{rev}}$ is the smallest relation such that
\begin{itemize}
    \item For all $i\in \mathbb{N}$, $(\overline{\sigma_i}\,\sigma_i,\epsilon)\in R_{\mathsf{rev}}$
    \item For all $i,j\in \mathbb{N}$ such that $|i-j|=1$,  $(\overline{\sigma_i}\,\sigma_j,\sigma_j\sigma_i\overline{\sigma_j}\,\overline{\sigma_i})\in R_{\mathsf{rev}}$
    \item For all $i,j\in \mathbb{N}$ such that $|i-j|>1$, $(\overline{\sigma_i}\,\sigma_j,\sigma_j\overline{\sigma_i})\in R_{\mathsf{rev}}$
\end{itemize}
\end{definition}

Note that the reversing relation preserves equivalence under the braid relations -- a quick induction shows the following : 

\begin{theorem} \label{thm:reversing_to_braid_group_equivalence} \leanlink{https://github.com/hannahfechtner/braids_better/blob/098309403cda787eaf1fa954e07ff0ca3b2927ee/BraidProject/SemiThueReversing.lean#L149-L184}
    For any $x,y\in BW$, if $x\xrightarrow[R_{\mathrm{rev}}]{*}y$, then $x\equiv y$ in the braid group.
\end{theorem}

\noindent We aim to connect reversing to grids in the following manner:

\begin{theorem}\label{thm:big_one}
    There exists a grid $\langle a, b, c,d\rangle$ if and only if $\overline{a}b\xrightarrow[R_{\mathrm{rev}}]{*}c\overline{d}$.
\end{theorem}

We will thus be able to use properties of grids (such as existence and uniqueness) to prove properties of our reversing algorithm (such as, respectively, termination and correctness). In fact it turns out that the above theorem is not quite as precise as we will need for termination. We will need to define size metrics on Semi-Thue derivations and on grids to eventually provide a natural-number (and hence well-founded) measure that strictly decreases as we progress through our algorithm. We first present our proof of the theorem in its original form, and then augment it with the size metric in the following chapter. The astute reader will note that the linked code for this section occasionally contains material regarding size (in Lean, a suffix of the form \lstinline{.length} for derivations and grids), and certain code structure (such as \lstinline{Type} in place of \lstinline{Prop}). This may be safely ignored for now, and will be discussed in detail in the following chapter. 

Returning to Theorem \ref{thm:big_one} as we have presented it here, the forward direction is nearly immediate by induction on grids:

\begin{theorem}\label{thm:big_one_easy}\leanlink{https://github.com/hannahfechtner/braids_better/blob/098309403cda787eaf1fa954e07ff0ca3b2927ee/BraidProject/SemiThueReversing.lean#L186-L204}
    If there exists a grid $\langle a, b, c,d\rangle$ , then $\overline{a}b\xrightarrow[R_{\mathrm{rev}}]{*}c\overline{d}$.
\end{theorem}

The reverse poses a much thornier problem. Patrick Dehornoy returned to the issue time and again. In one of his earliest papers, he took a different approach to showing termination for the subword reversing algorithm.

In his article ``\textit{Deux propriétés des groupes de tresses}'' from 1992 \cite{dehornoy_deux_1992}, Dehornoy created an extended alphabet which grouped together letters from the original alphabet into single letters in the new alphabet. For example, he added to our $S^*$ the letters $A=\sigma_1\sigma_2$ and $B=\sigma_2\sigma_1$. This was designed so that any application of reversing did not change the length of the word in the new alphabet. Thus, we see that $\overline{\sigma_1}\sigma_2\xrightarrow[R_{\mathrm{rev}}]{*} \sigma_2\sigma_1\overline{\sigma_2}\,\overline{\sigma_1}=B\overline{A}$, and so the length in the new alphabet remains 2. He carefully defines the new alphabet and shows that it is closed under applying reversing rewrites to negative-positive pairs. Thus, we have termination for reversing over the new alphabet. He then shows that this implies termination for reversing over the original alphabet. 

I have not attempted to implement this proof, and so I cannot speak to its rigor. It appears to be a solid, reliable proof. Nevertheless, there is still interest in the grid-based proof, due to its clarity of exposition and connection with other mathematical structures such as van Kampen diagrams. Clearly, Dehornoy preferred the grid representation, as nearly all of his following papers on the matter (approximately 10) discuss grids instead of the extended alphabet. He had a certain freedom in experimenting with different, intuitive options for grids, knowing that he had shown termination by a different method earlier. 

In his 2019 textbook, he gives the following sketch [variable names changed to match ours] :

\begin{displayquote}
    ... when $\overline a b \;\xrightarrow[R_{\mathrm{rev}}]{*} c\overline{d}$, consider a finite sequence of reversing steps leading from $\overline a b$ to $c\overline{d}$. Let $\Gamma_0$ be the diagram obtained by writing vertically and from top to bottom the letters of $a$ and horizontally from left to right those of $b$ starting from a common summit. Then following the sequence of elementary reversings starting with $\overline a b$ corresponds to constructing step by step from $\Gamma_0$ a grid $\Gamma$ with source $(a, b)$ and target $(c, d)$. CITE
\end{displayquote}

The sketch is compelling -- one can easily envision the grid being built step-by-step. Earlier papers of Dehornoy's which discuss reversing grids give additional ideas. In a 2011 paper \cite{dehornoy_subword_2011}, he begins with the case of general van Kampen diagrams \footnote{There is another, well-established method of constructing general van Kampen diagrams (using a balloon diagram; see \url{https://ntouikan.ext.unb.ca/MATH6022/IntroCGGT/html_output/section-14.html}), which differs from Dehornoy's description}:

\begin{displayquote}
Lemma 1.1 (folklore). If $(S, R)$ is a semigroup presentation, then two words $w$ and $w'$ of $S^{*}$
are $R$-equivalent if and only if there exists an $(S, R)$-van Kampen diagram for $(w, w')$.

Proof (sketch). If $(w_0, ..., w_p)$ is an $R$-derivation from $w$ to $w'$, then drawing paths labeled
with the successive words $w_i$ one below the other and identifying the unchanged letters yields
a van Kampen diagram for $(w, w')$. Conversely, if $K$ is a van Kampen diagram for $(w, w')$,
one obtains an $R$-derivation from $w$ to $w'$ by enumerating the labels in a sequence of paths
from the source of $K$ to its target that diﬀer by one face at a time.
\end{displayquote}

Again, this feels intuitively justifiable, but does not suffice for a rigorous proof. In that same article (written before the textbook came out), he also gives a description of how one might extend this method to rectangular grids, using dotted lines to identify edges that would have touched in a regular van Kampen diagram (essentially, $\varepsilon$-labelled edges in a grid, but missing their orientation).\footnote{Dehornoy had also used the curved lines for $\varepsilon$-labelled edges (no arrows, unsigned) in a paper from 2000. \cite{dehornoy_completeness_2000}.} The below figure is from page 7 of Dehornoy's 2011 ``\textit{The Subword Reversing Method}'' : 

\begin{center}
    \includegraphics[width=0.7\linewidth]{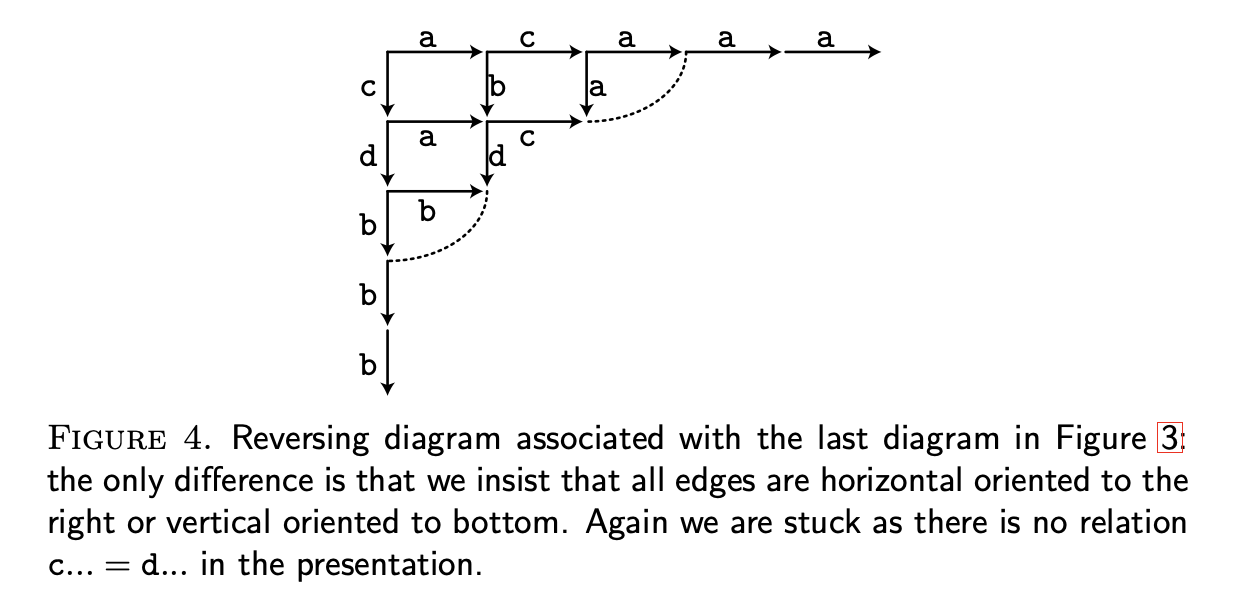}
\end{center}

In 2006, Dehornoy returned to the idea of the extended alphabet, and showed grids where each cell edge corresponded to a letter of the extended alphabet \cite{dehornoy_word_2006}. The below figure is from page 10 of Dehornoy and Wiest's ``\textit{On Word Reversing in Braid Groups}'' : 

\begin{center}
    \includegraphics[width=0.7\linewidth]{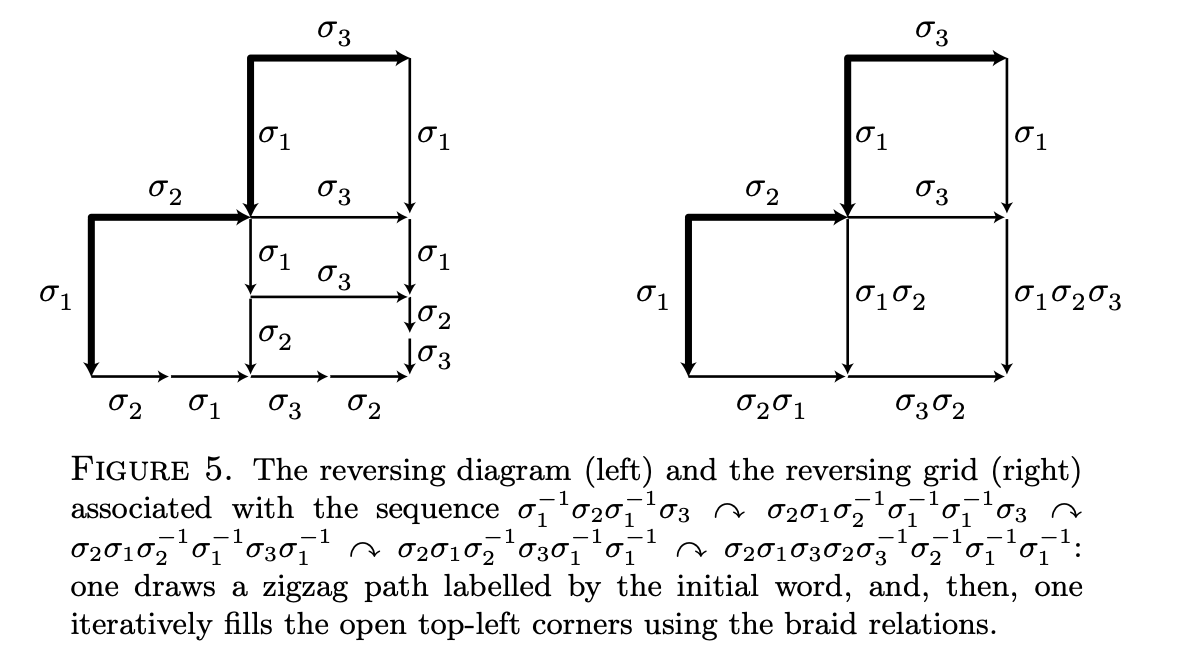}
\end{center}

This approach is tempting. However, it was not immediately applicable to our algorithm's termination, as we are working over the original alphabet. The proof translating termination from one alphabet to the other is non-trivial, and a metric defined on the extended-alphabet grid would not increase with each one of our individual rewrites. In any case, Dehornoy's sketch of extended-alphabet grids would have needed a more detailed treatment for formalization, and would have still involved the two largest parts of our argument in this chapter : translating one rewriting sequence for another and a definition of a grid-under-construction. In some sense, such extended-alphabet grids are less satisfying, as they hide many steps of the algorithm, all combined into a single cell. 

Returning to our definition of grids, my early naive attempts based on keeping track of the cells filled in unfortunately fail to capture the complexity of the construction required. It certainly feels so simple, so intuitive, and just out of reach.  Below, we give a novel (to the best of my knowledge) proof of this claim. The proof defines two intermediate structures: grid-style rewriting and partial grids. First, we convert the normal rewriting to grid-style rewriting. Then we show grid-style rewriting corresponds to a partial grid. Finally, we connect the idea of partial grids to regular grids. 

\subsection{Reading Grids Backwards}

We begin by enlarging our definition of grids.
At the moment, we have only discussed reading words off grids following the direction of the arrows. We can expand our interpretation of a grid structure by defining a procedure for reading \textit{against} the direction of an arrow. 

Instead of the labels on arrows representing elements of $BW^+$ (the set of braid monoid words), we now allow them to represent elements of $BW$ (the set of braid group words). 
We note that we are discussing words, not elements of the braid group. To make this clear, following Dehornoy, we continue to write $\overline{\sigma_i}$ in place of $\sigma_i^{-1}$. Notably, $\overline{\sigma_i}\sigma_i\neq\varepsilon$, for they are clearly different strings. However, $\overline{\sigma_i}\sigma_i\equiv\varepsilon$. 

We obtain such $\overline{\sigma_i}$ from the grid by reading against the direction of the arrows. Thus, we can read two paths from the bottom-left to the top-right in the following grid. The upper path is denoted $\overline{\sigma_i}\, \sigma_j$ and the lower is $\sigma_l\overline{\sigma_k}\,.$
\begin{center}
    \begin{tikzcd}
{} \arrow[d, "\sigma_i"'] \arrow[r, "\sigma_j"] & {} \arrow[d, "\sigma_k"] \\
{} \arrow[r, "\sigma_l"']                        & {}                      
\end{tikzcd}
\end{center}

For a side labelled with a word of length $>1$, we read the letters of the word in opposite order, and add an overbar to each letter. Consider the following grid:

\begin{center}
    \begin{tikzcd}
{} \arrow[d, "\sigma_i"'] \arrow[r, "\sigma_l"] & {} \arrow[d, "\sigma_i"] \\
{} \arrow[d, "\sigma_j"']                       & {} \arrow[d, "\sigma_j"] \\
{} \arrow[d, "\sigma_k "']                      & {} \arrow[d, "\sigma_k"] \\
{} \arrow[r, "\sigma_l"']                       & {}                      
\end{tikzcd}
\end{center}

The upper bottom-left to top-right path would be read as  $\overline{\sigma_k}\, \overline{\sigma_j} \, \overline{\sigma_i}\sigma_l$, and the lower $\sigma_l\overline{\sigma_k}\, \overline{\sigma_j} \, \overline{\sigma_i}$. 

\noindent Given a word $w=\sigma_{a_0}\sigma_{a_1}...\sigma_{a_n}$, let $\overline{w}$ denote $\overline{\sigma_{a_n}}\, ...\overline{\sigma_{a_1}}\, \overline{\sigma_{a_0}}$

Notably, we have only given examples of reading a downwards-pointing arrow upwards. We could read left-to-right arrows backwards as well, using the same procedure. However, we will not find a use for this in the following algorithm.

\leanlink{https://github.com/hannahfechtner/braids_better/blob/098309403cda787eaf1fa954e07ff0ca3b2927ee/BraidProject/Additions/SignedList.lean#L1}In Lean, we implement such words as elements of \lstinline{List (Nat x Bool)}. For example,  \lstinline{[(4, false), (3, true)]} represents $\overline{\sigma_4}\,\sigma_3$. Note that we choose \lstinline{List} over \lstinline{FreeMonoid} -- this is a nod to Theorem \ref{thm:big_one}, which discusses Semi-Thue rewriting, not presented monoids.

\subsection{Empty Arrows}

Recall that we omitted $\varepsilon$-labelled arrows in our definition of grids. We did not distinguish between the two following grids:
\begin{tikzcd}
{} \arrow[r, "\sigma_2"] \arrow[d, "\sigma_0"'] & {} \arrow[d, "\sigma_0"] \\
{} \arrow[r, "\sigma_2"']                       & {}                      
\end{tikzcd}
and
\begin{tikzcd}
{} \arrow[r, "\sigma_2"] \arrow[d, "\sigma_0"'] & {} \arrow[d, "\sigma_0"] \arrow[r, "\varepsilon"] & {} \arrow[d, "\sigma_0"] \\
{} \arrow[r, "\sigma_2"']                       & {} \arrow[r, "\varepsilon"']                      & {}                      
\end{tikzcd}. For our proof of theorem \ref{thm:big_one}, the more granular detail proof provided by tracking empty arrows will prove crucial. So, we again add more detail to the definition of our grid. 

We add a new symbol $\mathds{1}$ to our alphabet. $\mathds{1}$ represents an $\varepsilon$-arrow. 
\begin{definition} \leanlink{https://github.com/hannahfechtner/braids_better/blob/098309403cda787eaf1fa954e07ff0ca3b2927ee/BraidProject/Additions/SignedOptionList.lean#L4}
    The \textbf{augmented braid generators} are given by the set $S^{\mathds{1}} := \{\mathds{1}\}
   \cup
   \bigcup_{i\in\mathbb{N}} \{\sigma_i\}$

    The \textbf{augmented braid alphabet} is the set 
     \[    \mathscr{A}^{\mathds{1}}:= S^{\mathds{1}}\cup\overline{S^{\mathds{1}}} =
   \{\mathds{1}, \overline{\mathds{1}}\}
   \cup
   \bigcup_{i\in\mathbb{N}} \{\sigma_i, \overline{\sigma_i}\}
 \]
 
    The set of \textbf{augmented braid words} $BW^{\mathds{1}}$ is given by
    \[
BW^{\mathds{1}}
:= \left(\mathscr{A}^{\mathds{1}}\right)^*=
\left(
   \{\mathds{1}, \overline{\mathds{1}}\}
   \cup
   \bigcup_{i\in\mathbb{N}} \{\sigma_i, \overline{\sigma_i}\}
\right)^*
    \]

        \leanlink{https://github.com/hannahfechtner/braids_better/blob/098309403cda787eaf1fa954e07ff0ca3b2927ee/BraidProject/Additions/SignedOptionList.lean#L151-L165} A subset of $BW^{\mathds{1}}$, the set of \textbf{positive augmented braid words} $BW^{\mathds{1}, +}$ is given by
    \[
BW^{\mathds{1}, +}
:= \left(S^{\mathds{1}}\right)^*
    \]
    
\end{definition}

\begin{remark}
    \leanlink{https://github.com/hannahfechtner/braids_better/blob/098309403cda787eaf1fa954e07ff0ca3b2927ee/BraidProject/Additions/SignedOptionList.lean#L167-L181}We occasionally find use for the terminology ``negative augmented braid word'', which is an element of $\left(\overline{S^{\mathds{1}}}\right)^*$.
\end{remark}

In Lean, $BW^{\mathds{1}}$ is represented by the type \lstinline{List (Option ℕ × Bool)}. The \lstinline{Bool} keeps track of whether there is an overline on the generator; a value of \lstinline{none} for \lstinline{Option ℕ} indicates the generator is a $\mathds{1}$, and a value of \lstinline{some n} for \lstinline{Option ℕ} means the generator is $\sigma_n$.

\begin{definition} \leanlink{https://github.com/hannahfechtner/braids_better/blob/098309403cda787eaf1fa954e07ff0ca3b2927ee/BraidProject/Additions/SignedList.lean#L158-L159}
    We define the following functions to map between $BW$ and $BW^{\mathds{1}}$:

    $i^{\mathds{1}} : BW \to BW^{\mathds{1}}$ is simply the inclusion map. 
    
    \leanlink{https://github.com/hannahfechtner/braids_better/blob/098309403cda787eaf1fa954e07ff0ca3b2927ee/BraidProject/Additions/SignedOptionList.lean#L123-L127} $s^{\mathds{1}} :  BW^{\mathds{1}} \to BW$ is the lift of \[
f(x) =
\begin{cases}
\sigma_i & \text{if } x = \sigma_i, \\
\varepsilon & \text{if } x = \mathds{1}.
\end{cases}
\]
\end{definition}

\begin{lemma}\label{lem:s_after_i} \leanlink{https://github.com/hannahfechtner/braids_better/blob/098309403cda787eaf1fa954e07ff0ca3b2927ee/BraidProject/Additions/SignedOptionList.lean#L246-L252}
    The function $s^{\mathds{1}}$ is a left-inverse of $i^{\mathds{1}}$. 
\end{lemma}
\subsection{An Example}

Recall the part of Theorem \ref{thm:big_one} to be proved :  if $\overline{a}b\xrightarrow[R_{\mathrm{rev}}]{*}c\overline{d}$, then there exists a grid $\langle a, b, c,d\rangle$. To motivate the definitions and other decisions in the upcoming proof, we present an illustrative example.  Let $a = \sigma_2\sigma_1$, $b=\sigma_0\sigma_1$, $c= \sigma_0\sigma_1\sigma_2$, and $d=\sigma_2 \sigma_1 \sigma_0$. Consider the following derivation, one of many which witnesses $\overline{a}b\xrightarrow[R_{\mathrm{rev}}]{*}c\overline{d}$ :
\begin{align*}
\overline{\sigma_2\sigma_1} \sigma_0\sigma_1 
=&\, \overline{\sigma_1}\,\overline{\sigma_2} \sigma_0\sigma_1 \\
 \to_{R_{\mathsf{rev}}} &\,\overline{\sigma_1}\,\sigma_0 \overline {\sigma_2}\sigma_1 & (1) \\
\to_{R_{\mathsf{rev}}} & \,\sigma_0\sigma_1 \overline{\sigma_0}\,\overline{\sigma_1}\, \overline {\sigma_2} \sigma_1 & (2) \\
\to_{R_{\mathsf{rev}}}& \,\sigma_0\sigma_1 \overline{\sigma_0}\,\overline{\sigma_1} \sigma_1 \sigma_2\overline {\sigma_1}\, \overline{\sigma_2} & (3)\\
\to_{R_{\mathsf{rev}}}& \,\sigma_0\sigma_1 \overline{\sigma_0}\sigma_2\overline {\sigma_1} \,\overline{\sigma_2} & (4)\\
\to_{R_{\mathsf{rev}}}& \, \sigma_0\sigma_1 \sigma_2\overline{\sigma_0}\, \overline {\sigma_1} \,\overline{\sigma_2} &(5)\\
= &\, \sigma_0\sigma_1 \sigma_2\overline{\sigma_2 \sigma_1 \sigma_0}
\end{align*}

We aim to construct a grid $\langle \sigma_2\sigma_1, \sigma_0\sigma_1, \sigma_0\sigma_1\sigma_2, \sigma_2\sigma_1\sigma_0\rangle$ : 
\begin{tikzcd}
{} \arrow[d, "\sigma_2\sigma_1"'] \arrow[r, "\sigma_0\sigma_1"] & {} \arrow[d, "\sigma_2\sigma_1\sigma_0"] \\
{} \arrow[r, "\sigma_0\sigma_1\sigma_2"']                       & {}                                      
\end{tikzcd}. We begin by writing out the spine of the grid : 

\begin{tikzcd}
{} \arrow[dd, "\sigma_2"'] \arrow[rr, "\sigma_0"] &  & {} \arrow[rr, "\sigma_1"] &  & {} \\
                                                  &  &                           &  &    \\
{} \arrow[dd, "\sigma_1"']                        &  &                           &  &    \\
                                                  &  &                           &  &    \\
{}                                                &  &                           &  &   
\end{tikzcd}

\noindent The first rewrite of $\overline{\sigma_2}\sigma_0$ to $\sigma_0\overline{\sigma_2}$ corresponds to inserting the cell 
\begin{tikzcd}
{} \arrow[r, "\sigma_0"] \arrow[d, "\sigma_2"'] & {} \arrow[d, "\sigma_2"] \\
{} \arrow[r, "\sigma_0"']                       & {}                      
\end{tikzcd}:

\begin{tikzcd}
{} \arrow[dd, "\sigma_2"'] \arrow[rr, "\sigma_0"]  &  & {} \arrow[rr, "\sigma_1"] \arrow[dd, "\sigma_2"] &  & {} \\
                                                   &  &                                                  &  &    \\
{} \arrow[rr, "\sigma_0"'] \arrow[dd, "\sigma_1"'] &  & {}                                               &  &    \\
                                                   &  &                                                  &  &    \\
{}                                                 &  &                                                  &  &   
\end{tikzcd}

We define the frontier of a partially-completed grid as its southeasternmost bottom-left-to-top-right path. The frontier of the above diagram is $\overline{\sigma_1}\sigma_0\overline{\sigma_2}\sigma_1$, exactly the same as step (1) of the reversing derivation. 

This procedure continues easily for the next three steps of the derivation, lines $2$ through $4$, with corresponding cells fitting into the frontier of diagram neatly:

\noindent 
\begin{tikzcd}
{} \arrow[dd, "\sigma_2"'] \arrow[rr, "\sigma_0"]  &                           & {} \arrow[rr, "\sigma_1"] \arrow[dd, "\sigma_2"] &  & {} \\
                                                   &                           &                                                  &  &    \\
{} \arrow[rr, "\sigma_0"'] \arrow[dd, "\sigma_1"'] &                           & {} \arrow[d, "\sigma_1"]                         &  &    \\
                                                   &                           & {} \arrow[d, "\sigma_0"]                         &  &    \\
{} \arrow[r, "\sigma_0"']                          & {} \arrow[r, "\sigma_1"'] & {}                                               &  &   
\end{tikzcd}, 
\begin{tikzcd}
{} \arrow[dd, "\sigma_2"'] \arrow[rr, "\sigma_0"]  &                           & {} \arrow[rr, "\sigma_1"] \arrow[dd, "\sigma_2"] &                           & {} \arrow[d, "\sigma_2"] \\
                                                   &                           &                                                  &                           & {} \arrow[d, "\sigma_1"] \\
{} \arrow[rr, "\sigma_0"'] \arrow[dd, "\sigma_1"'] &                           & {} \arrow[r, "\sigma_1"'] \arrow[d, "\sigma_1"]  & {} \arrow[r, "\sigma_2"'] & {}                       \\
                                                   &                           & {} \arrow[d, "\sigma_0"]                         &                           &                          \\
{} \arrow[r, "\sigma_0"']                          & {} \arrow[r, "\sigma_1"'] & {}                                               &                           &                         
\end{tikzcd}, and
\begin{tikzcd}
{} \arrow[dd, "\sigma_2"'] \arrow[rr, "\sigma_0"]  &                           & {} \arrow[rr, "\sigma_1"] \arrow[dd, "\sigma_2"]   &                                                    & {} \arrow[d, "\sigma_2"] \\
                                                   &                           &                                                    &                                                    & {} \arrow[d, "\sigma_1"] \\
{} \arrow[rr, "\sigma_0"'] \arrow[dd, "\sigma_1"'] &                           & {} \arrow[r, "\sigma_1"'] \arrow[d, "\sigma_1"]    & {} \arrow[r, "\sigma_2"'] \arrow[d, "\varepsilon"] & {}                       \\
                                                   &                           & {} \arrow[d, "\sigma_0"] \arrow[r, "\varepsilon"'] & {}                                                 &                          \\
{} \arrow[r, "\sigma_0"']                          & {} \arrow[r, "\sigma_1"'] & {}                                                 &                                                    &                         
\end{tikzcd}
We encounter an issue when moving to line 5. There, we replace $\overline{\sigma_0}\sigma_2$ with $\sigma_2\overline{\sigma_0}$. The corresponding cell is 
\begin{tikzcd}
{} \arrow[r, "\sigma_2"] \arrow[d, "\sigma_0"'] & {} \arrow[d, "\sigma_0"] \\
{} \arrow[r, "\sigma_2"']                       & {}                      
\end{tikzcd}. But this does not fit into the frontier! The frontier edges \begin{tikzcd}
{} \arrow[d, "\sigma_0"']   \\
{}                                          
\end{tikzcd} and \begin{tikzcd}
{} \arrow[r, "\sigma_2"]  & {}  \\
{}                    & {}                      
\end{tikzcd} are separated by $\varepsilon$-labelled edges.

To insert this cell, we must first add in the following two cells (which we will call \textit{trivial} cells), which do not correspond to a step of the derivation :
\begin{center}
\begin{tikzcd}
{} \arrow[dd, "\sigma_2"'] \arrow[rr, "\sigma_0"]  &                           & {} \arrow[rr, "\sigma_1"] \arrow[dd, "\sigma_2"]   &                                                    & {} \arrow[d, "\sigma_2"]    \\
                                                   &                           &                                                    &                                                    & {} \arrow[d, "\sigma_1"]    \\
{} \arrow[rr, "\sigma_0"'] \arrow[dd, "\sigma_1"'] &                           & {} \arrow[r, "\sigma_1"'] \arrow[d, "\sigma_1"]    & {} \arrow[r, "\sigma_2"'] \arrow[d, "\varepsilon"] & {} \arrow[d, "\varepsilon"] \\
                                                   &                           & {} \arrow[d, "\sigma_0"] \arrow[r, "\varepsilon"'] & {} \arrow[r, "\sigma_2"'] \arrow[d, "\sigma_0"]    & {}                          \\
{} \arrow[r, "\sigma_0"']                          & {} \arrow[r, "\sigma_1"'] & {} \arrow[r, "\varepsilon"']                       & {}                                                 &                            
\end{tikzcd}
\end{center}
Further, note that adding in trivial cells must be done recursively. To see why, consider the following snippet of a grid : 
\begin{center}
\begin{tikzcd}
{} \arrow[d, "\sigma_1"'] \arrow[r, "\sigma_1"]    & {} \arrow[d, "\varepsilon"] \arrow[r, "\sigma_2"] & {} \arrow[r, "\sigma_1"] & {} \\
{} \arrow[r, "\varepsilon"']  & {}                                                &                          &    
\end{tikzcd}
\end{center}

\noindent When we add in the appropriate trivial cell \begin{tikzcd}
{} \arrow[r, "\sigma_2"] \arrow[d, "\varepsilon"'] & {} \arrow[d, "\varepsilon"] \\
{} \arrow[r, "\sigma_2"']                       & {}                      
\end{tikzcd} we obtain the following grid snippet:

\begin{center}
\begin{tikzcd}
{} \arrow[d, "\sigma_1"'] \arrow[r, "\sigma_1"] & {} \arrow[d, "\varepsilon"] \arrow[r, "\sigma_2"] & {} \arrow[r, "\sigma_1"] \arrow[d, "\varepsilon"] & {} \\
{} \arrow[r, "\varepsilon"']                    & {} \arrow[r, "\sigma_2"']                         & {}                                                &   
\end{tikzcd}
\end{center}
This snippet now requires another trivial cell to be inserted.

Returning to our main example, we may now complete the grid following the rewrite in step (5), obtaining
\begin{center}
\begin{tikzcd}
{} \arrow[dd, "\sigma_2"'] \arrow[rr, "\sigma_0"]  &                           & {} \arrow[rr, "\sigma_1"] \arrow[dd, "\sigma_2"]   &                                                    & {} \arrow[d, "\sigma_2"]    \\
                                                   &                           &                                                    &                                                    & {} \arrow[d, "\sigma_1"]    \\
{} \arrow[rr, "\sigma_0"'] \arrow[dd, "\sigma_1"'] &                           & {} \arrow[r, "\sigma_1"'] \arrow[d, "\sigma_1"]    & {} \arrow[r, "\sigma_2"'] \arrow[d, "\varepsilon"] & {} \arrow[d, "\varepsilon"] \\
                                                   &                           & {} \arrow[d, "\sigma_0"] \arrow[r, "\varepsilon"'] & {} \arrow[r, "\sigma_2"'] \arrow[d, "\sigma_0"]    & {} \arrow[d, "\sigma_0"]    \\
{} \arrow[r, "\sigma_0"']                          & {} \arrow[r, "\sigma_1"'] & {} \arrow[r, "\varepsilon"']                       & {} \arrow[r, "\sigma_2"']                          & {}                         
\end{tikzcd}
\end{center}

This has frontier $\sigma_0\sigma_1\mathds{1}\sigma_2\overline{\sigma_0}\,\overline{\mathds{1}}\,\overline{\sigma_1}\,\overline{\sigma_2}$. Note that $s^{\mathds{1}}(\sigma_0\sigma_1\mathds{1}\sigma_2\overline{\sigma_0}\,\overline{\mathds{1}}\,\overline{\sigma_1}\,\overline{\sigma_2})=\sigma_0\sigma_1\sigma_2\overline{\sigma_0}\,\overline{\sigma_1}\,\overline{\sigma_2}$, which equals step $(5)$ of the derivation.

To give a rigorous proof, we must first give a definition of partially-constructed grids. Then, we must develop a procedure for inserting trivial cells. We \textit{could} enact this procedure on the partial grids. It turns out to be simpler to translate the given derivation to one tracking $\varepsilon$'s explicitly, and then have a clean, straightforward translation from the new derivation to partial grids. This is largely due to strings being simpler to work with than two-dimensional diagrams. We define this new derivation in the following section.

\subsection {Grid-Style Rewriting}

We define a relation $R_{\mathsf{grid}}$ which tracks $\varepsilon$-labelled arrows, which are represented by $\mathds{1}$ in our alphabet.
This relation is defined to correspond neatly with grid cells; if we have $R_{\mathsf{grid}}(x, y)$, we know there is a cell with $x$ as the left-top border and $y$ as the bottom-right border (reading, as always in the section, north-east, and using an overbar to denote travelling opposite an arrow's direction).

\begin{definition} \leanlink{https://github.com/hannahfechtner/braids_better/blob/098309403cda787eaf1fa954e07ff0ca3b2927ee/BraidProject/Relations.lean#L80-L84}
The relation $R_{\mathsf{nontriv}}$ is very similar to $R_{\mathsf{rev}}$. It differs only in the first case (see below); rewriting under $R_{\mathsf{rev}}$ simply deletes an occurrence of $\overline{\sigma_i}\sigma_i$, whereas under $R_{\mathsf{nontriv}}$, $\overline{\sigma_i}\sigma_i$ would be replaced by $\mathds{1}\overline{\mathds{1}}$. $R_{\mathsf{nontriv}}$ is the smallest relation satisfying

\hspace{1cm} $R_{\mathsf{nontriv}} ( \overline{\sigma}_i \sigma_i,\mathds{1}\overline{\mathds{1}})$ for all $i \in \mathbb{N}$. This corresponds to the cell \begin{tikzcd}
{} \arrow[r, "\sigma_i"] \arrow[d, "\sigma_i"'] & {} \arrow[d, "\varepsilon"] \\
{} \arrow[r, "\varepsilon"']                   & {}               
\end{tikzcd}

\hspace{1cm} $R_{\mathsf{nontriv}} ( \overline{\sigma_i}\sigma_j,\sigma_j\sigma_i\overline{\sigma_j}\:\overline{\sigma_i})$ when $|i-j|=1$. This corresponds to the cell \begin{tikzcd}
{} \arrow[rr, "\sigma_j"] \arrow[dd, "\sigma_i"'] &                           & {} \arrow[d, "\sigma_i"] \\
                                                  &                           & {} \arrow[d, "\sigma_j"] \\
{} \arrow[r, "\sigma_j"']                         & {} \arrow[r, "\sigma_i"'] & {}                      
\end{tikzcd}

\hspace{1cm} $R_{\mathsf{nontriv}} ( \overline{\sigma_i}\sigma_j,\sigma_j\overline{\sigma_i})$ when $|i-j|\geq 2$. This corresponds to the cell \begin{tikzcd}
{} \arrow[r, "\sigma_j"] \arrow[d, "\sigma_i"'] & {} \arrow[d, "\sigma_i"] \\
{} \arrow[r, "\sigma_j"']                   & {}               
\end{tikzcd}

    \noindent \leanlink{https://github.com/hannahfechtner/braids_better/blob/098309403cda787eaf1fa954e07ff0ca3b2927ee/BraidProject/Relations.lean#L88-L91} The relation $R_{\mathsf{triv}}$ tracks ``trivial'' rewrites which do not occur in an $R_{\mathsf{rev}}$-- derivation. $R_{\mathsf{triv}}$ is the smallest relation satisfying
    
    \hspace{1cm} $R_{\mathsf{triv}} ( \overline{\mathds{1}} \sigma_i,\sigma_i \overline{\mathds{1}})$ for all $i\in \mathbb{N}$. This corresponds to the cell \begin{tikzcd}
{} \arrow[r, "\sigma_{i}"] \arrow[d, "\varepsilon"'] & {} \arrow[d, "\varepsilon"] \\
{} \arrow[r, "\sigma_{i}"']                   & {}               
\end{tikzcd}

\hspace{1cm} $R_{\mathsf{triv}} ( \overline{\sigma}_i \;\mathds{1},\mathds{1}\:\overline{\sigma}_i)$ for all $i\in \mathbb{N}$. This corresponds to the cell \begin{tikzcd}
{} \arrow[r, "\varepsilon"] \arrow[d, "\sigma_i"'] & {} \arrow[d, "\sigma_i"] \\
{} \arrow[r, "\varepsilon"']                   & {}               
\end{tikzcd}

\hspace{1cm} $R_{\mathsf{triv}} ( \overline{\mathds{1}} \;\mathds{1}, \;\mathds{1} \; \overline{\mathds{1}})$. This corresponds to the cell \begin{tikzcd}
{} \arrow[r, "\varepsilon"] \arrow[d, "\varepsilon"'] & {} \arrow[d, "\varepsilon"] \\
{} \arrow[r, "\varepsilon"']                   & {}               
\end{tikzcd}

\noindent  \leanlink{https://github.com/hannahfechtner/braids_better/blob/098309403cda787eaf1fa954e07ff0ca3b2927ee/BraidProject/Relations.lean#L45-L52} Finally, $R_{\mathsf{grid}}$ is the union of $R_{\mathsf{triv}}$ and $R_{\mathsf{nontriv}}$.
\end{definition}

We next make precise the close relationship between $R_{\mathsf{rev}}$ and $R_{\mathsf{nontriv}}$ :
\begin{theorem}\label{thm:rev_to_nontriv_relation} \leanlink{https://github.com/hannahfechtner/braids_better/blob/098309403cda787eaf1fa954e07ff0ca3b2927ee/BraidProject/Relations.lean#L181-L195}
    If $R_{\mathsf{rev}}(x,y)$, there exists $y'\in BW^{\mathds{1}}$ such that $R_{\mathsf{nontriv}}(i^{\mathds{1}}(x),y')$. This $y'$ happens to be unique, and have the property that $s^{\mathds{1}}(y')=y$.
\end{theorem}
We next give terminology for augmented braid words which can or can not be reduced by the above relations.
\begin{definition}

    An element $b\in BW^{\mathds{1}}$ is \textbf{trivially-reducible} if it is of the form $b_1xb_2$, where $(x,y)\in R_{\mathsf{triv}}$.

    Otherwise, it is called \textbf{trivially-irreducible}. \leanlink{https://github.com/hannahfechtner/braids_better/blob/098309403cda787eaf1fa954e07ff0ca3b2927ee/BraidProject/Irreducibility.lean#L50-L52} This is equivalent to saying $b$ has no subwords of the form $\overline{\sigma_i}\mathds{1}$, $\overline{\mathds{1}}\sigma_i$, or $\overline{\mathds{1}}\mathds{1}$.

    We have analogous definitions for nontrivially-reducible and nontrivially-irreducible, and $R_{\mathsf{grid}}$ - reducible and $R_{\mathsf{grid}}$ - irreducible.
\end{definition}

\begin{lemma}\label{lem:i_irr} \leanlink{https://github.com/hannahfechtner/braids_better/blob/098309403cda787eaf1fa954e07ff0ca3b2927ee/BraidProject/Irreducibility.lean#L82-L97}
    For any $b\in BW$, $i^{\mathds{1}}(b)$ is trivially-irreducible. 
\end{lemma}

\noindent We give a specific fact about irreducibility to be used in a later proof :
\begin{lemma}\label{lem:pos_neg_irr} \leanlink{https://github.com/hannahfechtner/braids_better/blob/098309403cda787eaf1fa954e07ff0ca3b2927ee/BraidProject/Irreducibility.lean#L181-L285}
    If $g$ is trivially-irreducible and $s^{\mathds{1}}(g)=c\overline{d}$ for some $c,d\in BW^{+}$, then there exist $e,f\in BW^{\mathds{1}, +}$ such that $g=e\overline{f}$.
\end{lemma}

\noindent We name the important words of length two :

\begin{definition}
    Given letters $\alpha,\beta\in S^{\mathds{1}}$, a word of the form $\overline{\alpha}\beta$ is called an \textbf{open pair}. A word in which at least one of $\alpha$ or $\beta$ equals $\mathds{1}$ is called a \textbf{trivial open pair}. A word of the form $\overline{\sigma_i}\sigma_j$ is called a \textbf{sigma-sigma open pair}. We use the adjective ``open'' to signify there is an opening for a grid cell in the corresponding diagram.

    We have the analogous vocabulary for pairs of the form $\alpha\overline{\beta}$, with ``closed'' replacing ``open''.
\end{definition}

\subsection{Building Grid-Style Rewriting}

This section's main result is the generation of an $R_{\mathsf{grid}}$ derivation from an $R_{\mathsf{rev}}$ derivation : 
\begin{theorem} \label{thm:new_rw} \leanlink{https://github.com/hannahfechtner/braids_better/blob/098309403cda787eaf1fa954e07ff0ca3b2927ee/BraidProject/ReversingToGridStyle.lean#L137-L152}
    Let $a, b, c, d \in BW^+$. If $\overline{a}b \xrightarrow[R_{\mathrm{rev}}]{*}c\overline{d}$, then there exist positive augmented words $e, f$ such that $i^{\mathds{1}}(\overline{a})i^{\mathds{1}}(b)\xrightarrow[R_{\mathrm{grid}}]{*}e\overline{f}$ . Moreover, deleting all instances of $\mathds{1}$ and $\overline{\mathds{1}}$ from $e$ leaves us with $c$; that is, $s^{\mathds{1}}(e)=c$. Similarly, $s^{\mathds{1}}(f)=d$.
\end{theorem}

Before we begin the proof, let us return to our example from earlier in the chapter. Recall the derivation:

\begin{align*}
\overline{\sigma_2\sigma_1} \sigma_0\sigma_1 
=&\, \overline{\sigma_1}\,\overline{\sigma_2} \sigma_0\sigma_1 \\
 \to_{R_{\mathsf{rev}}} &\,\overline{\sigma_1}\,\sigma_0 \overline {\sigma_2}\sigma_1 & (1) \\
\to_{R_{\mathsf{rev}}} & \,\sigma_0\sigma_1 \overline{\sigma_0}\,\overline{\sigma_1}\, \overline {\sigma_2} \sigma_1 & (2) \\
\to_{R_{\mathsf{rev}}}& \,\sigma_0\sigma_1 \overline{\sigma_0}\,\overline{\sigma_1} \sigma_1 \sigma_2\overline {\sigma_1}\, \overline{\sigma_2} & (3)\\
\to_{R_{\mathsf{rev}}}& \,\sigma_0\sigma_1 \overline{\sigma_0}\sigma_2\overline {\sigma_1} \,\overline{\sigma_2} & (4)\\
\to_{R_{\mathsf{rev}}}& \, \sigma_0\sigma_1 \sigma_2\overline{\sigma_0}\, \overline {\sigma_1} \,\overline{\sigma_2} &(5)\\
= &\, \sigma_0\sigma_1 \sigma_2\overline{\sigma_2 \sigma_1 \sigma_0}
\end{align*}

The following derivation will suffice for Theorem \ref{thm:new_rw}, with $e = \sigma_0\sigma_1 \mathds{1} \sigma_2$ and $f = \sigma_2 \sigma_1 \mathds{1} \sigma_0$

\begin{align*}
\overline{\sigma_2\sigma_1} \sigma_0\sigma_1 
=&\, \overline{\sigma_1}\,\overline{\sigma_2} \sigma_0\sigma_1 \\
 \to_{R_{\mathsf{rev}}} &\,\overline{\sigma_1}\,\sigma_0 \overline {\sigma_2}\sigma_1 & (1) \\
\to_{R_{\mathsf{rev}}} & \,\sigma_0\sigma_1 \overline{\sigma_0}\,\overline{\sigma_1}\, \overline {\sigma_2} \sigma_1 & (2) \\
\to_{R_{\mathsf{rev}}}& \,\sigma_0\sigma_1 \overline{\sigma_0}\,\overline{\sigma_1} \sigma_1 \sigma_2\overline {\sigma_1}\, \overline{\sigma_2} & (3)\\
\to_{R_{\mathsf{rev}}}& \,\sigma_0\sigma_1 \overline{\sigma_0} \,\mathds{1}\overline{\mathds{1}}\sigma_2\overline {\sigma_1} \,\overline{\sigma_2} & (4)\\
\to_{R_{\mathsf{rev}}}& \,\sigma_0\sigma_1 \mathds{1} \overline{\sigma_0} \,\overline{\mathds{1}}\sigma_2\overline {\sigma_1} \,\overline{\sigma_2} & (4.1)\\
\to_{R_{\mathsf{rev}}}& \,\sigma_0\sigma_1 \mathds{1} \overline{\sigma_0} \sigma_2\,\overline{\mathds{1}}\,\overline {\sigma_1} \,\overline{\sigma_2} & (4.2)\\
\to_{R_{\mathsf{rev}}}& \,\sigma_0\sigma_1 \mathds{1} \sigma_2 \overline{\sigma_0} \,\overline{\mathds{1}}\,\overline {\sigma_1} \,\overline{\sigma_2} &(5)\\
= &\, \sigma_0\sigma_1 \mathds{1} \sigma_2\overline{\sigma_2 \sigma_1 \mathds{1} \sigma_0}
\end{align*}

Here, the work lies in moving the $\mathds{1}$s in lines $(4)$ and $(4.1)$. Then, we are relieved to see that in line $(4.2)$, we are poised to make the rewrite in line $(5)$.

Moving on from this example, we begin by defining a general procedure for performing the trivial reductions (i.e. moving the $\mathds{1}$ and $\overline{\mathds{1}}s$). The simplest way to define the procedure for trivial reduction is to repeatedly search for and rewrite the first occurrence of a trivial pair. We can give a lexicographic ordering on our alphabet, and in that way show this procedure does terminate.\footnote{For example, \leanlinkfootnote{https://github.com/hannahfechtner/braids_better/blob/098309403cda787eaf1fa954e07ff0ca3b2927ee/BraidProject/AlternativeChoices/AlphabetRel.lean\#L93} we might use the following well-founded relation : 
\begin{itemize}
\item $\sigma_i< \sigma_j$ iff $i<j$

\item $\overline {\sigma_i}< \overline{\sigma_j}$ iff $i<j$ 

\item $\sigma_i< \overline{\sigma_j}$  for all natural numbers $i, j$ 

\item $\mathds{1}< \overline{\mathds{1}}$

\item $\sigma_i< \mathds{1}$ for all natural numbers $i$

\item $\overline{\mathds{1}}< \overline{\sigma_i}$ for all natural numbers $i$
\end{itemize}

\leanlinkfootnote{https://github.com/leanprover-community/mathlib4/blob/a2f40e9ce2b3d1ebc9a5e057463f56fb87f119a8/Mathlib/Data/List/Shortlex.lean\#L37-L41} We extend this to the associated shortlex ordering over strings. Although this method was not ultimately used in the project, it did lead to the implementation of the shortlex order, now in Mathlib.} However, we will need to prove a number of properties about this procedure's result (such as trivial irreducibility), so we need a definition more amenable to inductive proofs. 

We first define a helper function \lstinline{concatenate_reduction}, which appends a single letter of our alphabet $\mathscr{A}^{\mathds{1}}$ to a word $b\in BW^{\mathds{1}}$ and applies any  trivial reductions involving said letter. We will in the future only apply this function to a trivially-irreducible $b$. 

\begin{definition} \leanlink{https://github.com/hannahfechtner/braids_better/blob/098309403cda787eaf1fa954e07ff0ca3b2927ee/BraidProject/MoveOnes.lean#L9-L20}
    Given a single letter $\alpha\in \mathscr{A}^{\mathds{1}}$ and a word $b\in BW^{\mathds{1}}$, we define 
    
    \lstinline{concatenate_reduction a b} as
    \begin{itemize}
        \item if $\alpha$ is $\overline{\mathds{1}}$ and $b=\sigma_ib_1$, return $\sigma_i\cdot($\lstinline{concatenate_reduction a b₁}$)$. (This case rewrites $\overline{\mathds{1}}\sigma_i$ to $\sigma_i\overline{\mathds{1}}$, and then continues on with any further trivial rewrites involving $\alpha$)
        \item if $\alpha$ is negative (an element of $\overline{S^{\mathds{1}}}$) and $b=\mathds{1}b_1$,  return $\mathds{1} \cdot ($\lstinline{concatenate_reduction a b₁}$)$. (In this case, $\alpha=\overline{\sigma_i}$ or $\alpha=\overline{\mathds{1}}$. So, this case either rewrites $\overline{\sigma_i}\mathds{1}$ to $\mathds{1}\overline{\sigma_i}$ or $\overline{\mathds{1}}\mathds{1}$ to $\mathds{1}\overline{\mathds{1}}$, and then continues on with any further trivial rewrites involving $\alpha$)
        
        \item otherwise, return $\alpha b$
    \end{itemize}
\end{definition}
    
\begin{lemma}\label{lem:concat_reduce_irr} \leanlink{https://github.com/hannahfechtner/braids_better/blob/098309403cda787eaf1fa954e07ff0ca3b2927ee/BraidProject/MoveOnes.lean#L282-L328}
    If $b$ is trivially-irreducible, so too is \lstinline{concatenate_reduction a b}. The proof proceeds by induction on $b$, and then casework on $\alpha$ and the first element of $b$.
\end{lemma}

\begin{theorem}\label{lem:concat_reduce_st_triv} \leanlink{https://github.com/hannahfechtner/braids_better/blob/098309403cda787eaf1fa954e07ff0ca3b2927ee/BraidProject/MoveOnes.lean#L158-L161}
    $\alpha b\xrightarrow[R_{\mathrm{triv}}]{*}$\lstinline{ concatenate_reduction a b}. The proof proceeds in the same style as that of the above theorem.
\end{theorem}

We are now prepared to define the algorithm \lstinline{move_ones} which trivially-reduces a given augmented braid word.

\begin{definition} \leanlink{https://github.com/hannahfechtner/braids_better/blob/098309403cda787eaf1fa954e07ff0ca3b2927ee/BraidProject/MoveOnes.lean#L209}
    Given $b\in BW^{\mathds{1}}$, we define \lstinline{move_ones b} as follows :
    \begin{itemize}
        \item If $b$ has length zero, then \lstinline{move_ones b} returns the empty word $\varepsilon$. 
        \item If $b$ has length at least one, it has the form $\alpha b_1$ for some $\alpha \in\mathscr{A}^{\mathds{1}}$ and $b_1\in BW^{\mathds{1}}$. Let $c$ be the result of the recursive call \lstinline{move_ones b₁}. Then we return \lstinline{concatenate_reduction a c}. 
    \end{itemize}
    In Lean, this reads as
    \leanlink{\https://github.com/hannahfechtner/braids_better/blob/098309403cda787eaf1fa954e07ff0ca3b2927ee/BraidProject/MoveOnes.lean#L209}\begin{lstlisting}
def move_ones (L : List (Option ℕ × Bool)) :=
  match L with
  | [] => []
  | head :: tail => concatenate_reduction head (move_ones tail)
\end{lstlisting}
\end{definition}

\begin{theorem}\label{thm:move_ones_equiv} \leanlink{https://github.com/hannahfechtner/braids_better/blob/098309403cda787eaf1fa954e07ff0ca3b2927ee/BraidProject/MoveOnes.lean#L265-L269}
    $b\xrightarrow[R_{\mathrm{triv}}]{*}$\lstinline{ move_ones b}. The proof proceeds by induction, using Lemma \ref{lem:concat_reduce_st_triv}. Notably, as $R_{\mathsf{triv}}\subseteq R_{\mathsf{grid}}$, $b\xrightarrow[R_{\mathrm{grid}}]{*}$\lstinline{ move_ones b}
\end{theorem}

\begin{theorem} \label{thm:move_ones_delete} \leanlink{https://github.com/hannahfechtner/braids_better/blob/098309403cda787eaf1fa954e07ff0ca3b2927ee/BraidProject/MoveOnes.lean#L271-L280}
\lstinline{move_ones} does not affect the relative positions of any $\sigma_i$. $s^{\mathds{1}}($\lstinline{move_ones b}$)=s^{\mathds{1}}(b)$
\end{theorem}

\begin{theorem}\label{thm:move_ones_irr} \leanlink{https://github.com/hannahfechtner/braids_better/blob/098309403cda787eaf1fa954e07ff0ca3b2927ee/BraidProject/MoveOnes.lean#L330-L334}
    \lstinline{move_ones} produces a trivially-irreducible word; for any $b\in BW^{\mathds{1}},$ \lstinline{move_ones b} is trivially-irreducible.
\end{theorem}

Let us take stock of where we are in terms of our example. The \lstinline{move_ones} algorithm effects the work done in steps $(4)-(4.2)$. We next noted that we were poised for the next step from the $R_{\mathsf{rev}}$ - derivation, able to make the rewrite from step $(5)$. The property of trivial-irreducibility just mentioned will in fact suffice, but not straightforwardly so. Instead, we introduce another notion which more immediately connects to preparedness for a nontrivial rewrite. 

We begin by defining a property \lstinline{pairAppears}.

\begin{definition} \leanlink{https://github.com/hannahfechtner/braids_better/blob/098309403cda787eaf1fa954e07ff0ca3b2927ee/BraidProject/Irreducibility.lean#L12-L14}
    An element $b\in BW^{\mathds{1}}$ has the \lstinline{pairAppears} property if for every substring $\overline{\sigma_i}\sigma_j$ appearing in $s^{\mathds{1}}(b)$, $\overline{\sigma_i}\sigma_j$ is a substring of $b$. 
\end{definition}

This is not quite enough to define the desired property. \leanlink{https://github.com/hannahfechtner/braids_better/blob/098309403cda787eaf1fa954e07ff0ca3b2927ee/BraidProject/Irreducibility.lean#L335} For example, the string $\overline{\sigma_1}\,\overline{\mathds{1}}\sigma_2\mathds{1}\,\overline{\sigma_1}\sigma_2$ satisfies \lstinline{pairAppears}, but the first occurrence of $\overline{\sigma_1}\sigma_2$ remains separated by a ${\mathds{1}}$. We define a stronger condition based on \lstinline{pairAppears} as follows:
\begin{definition} \leanlink{https://github.com/hannahfechtner/braids_better/blob/098309403cda787eaf1fa954e07ff0ca3b2927ee/BraidProject/Irreducibility.lean#L38}
    An element $b\in BW^{\mathds{1}}$ has the \lstinline{pairsTogether} property if every substring of $b$ satisfies \lstinline{pairAppears}.
\end{definition}

Showing that an augmented braid word is trivially-irreducible is an easy way to show that it satisfies \lstinline{pairsTogether}, thanks to the following theorem :
\begin{theorem}\label{thm:pts_of_irr} \leanlink{https://github.com/hannahfechtner/braids_better/blob/098309403cda787eaf1fa954e07ff0ca3b2927ee/BraidProject/Irreducibility.lean#L353-L476}
    Any trivially-irreducible $b\in BW^{\mathds{1}}$ satisfies \lstinline{pairsTogether}
\end{theorem}

\begin{proof} \textit{(Sketch)}
    We show that any irreducible $b$ satisfies \lstinline{pairAppears} by induction on the length of $b$, and then casework on its first and second elements. Then we quickly show that any subword of a trivially-irreducible word is itself trivially-irreducible, and combine this fact with the definition of \lstinline{pairsTogether} to obtain our conclusion.
\end{proof}

Note that \lstinline{pairsTogether} is strictly weaker than trivial irreducibility. For example, $\sigma_2\overline{\mathds{1}}\mathds{1}$ satisfies \lstinline{pairsTogether} vacuously, but is not trivially-irreducible.

We give two specific lemmas about open-pair infixes of trivially-irreducible strings, which will be used for the proof of Theorem \ref{thm:new_rw} : 
\begin{lemma}\label{lem:infix_of_irr} \leanlink{https://github.com/hannahfechtner/braids_better/blob/098309403cda787eaf1fa954e07ff0ca3b2927ee/BraidProject/PairInfixIrreducibleCases.lean#L117-L139}
    Let $w\in BW^{\mathds{1}}$ be trivially-irreducible, and assume $s^{\mathds{1}}(w)$ has an infix of the form $\overline{\sigma_i}\sigma_j$ -- more specifically, assume $s^{\mathds{1}}(w)=c\overline{\sigma_i}\sigma_jd$ for $c,d\in BW^{\mathds{1}}$. Then there exist $w_1, w_2\in BW^{\mathds{1}}$ such that $w = w_1 \overline{\sigma_i}\sigma_j w_2$, $c=s^{\mathds{1}}(w_1)$, and $d=s^{\mathds{1}}(w_2)$.
\end{lemma}
\begin{proof} (sketch)
    The proof proceeds by stacking induction on $w$ with lemmas also proved by induction on w (for example, a lemma where ``infix'' is replaced by ``prefix''), and critically relies on the trivial-irreducibility of $w$. Otherwise, we would only know that there exist $w_1, w_2, w_m\in BW^{\mathds{1}}$ such that $w = w_1 w_m w_2$ with $s^{\mathds{1}}(w_m)=\overline{\sigma_i}\sigma_j$.
\end{proof}
\begin{lemma} \label{lem:list_split}\leanlink{https://github.com/hannahfechtner/braids_better/blob/098309403cda787eaf1fa954e07ff0ca3b2927ee/BraidProject/PairInfixIrreducibleCases.lean#L141-L163}
    Given trivially-irreducible $a, b\in BW^{\mathds{1}}$ and any $c,d\in BW$, if $s^{\mathds{1}}(a) \overline{\sigma_i} \sigma_j s^{\mathds{1}}(b)=c \,\overline{\sigma_i} \sigma_j \,d$, then we have one of the following three cases, depending on if $c$ is equal to, longer than, or shorter than $s^{\mathds{1}}(a)$ :
    \begin{itemize}
    \item Case 1 : $s^{\mathds{1}}(a)= c$ and $s^{\mathds{1}}(b)=d$.
    
    \begin{tikzpicture}[
  box/.style={draw, minimum height=.5cm, align=center},
  every node/.style={inner sep=0pt}
]

\def\a{4.5}   
\def\b{2.0}   
\def\c{4.5}   

\def\d{4.5}   
\def\e{2.0}   
\def\f{4.5}   

\def\topY{0}
\def\botY{-.7}

\node[box, minimum width=\a cm] at ({\a/2}, \topY)
  {$s^{\mathds{1}}(a)$};

\node[box, minimum width=\b cm] at ({\a+\b/2}, \topY)
  {$\overline{\sigma_i}\,\sigma_j$};

\node[box, minimum width=\c cm] at ({\a+\b+\c/2}, \topY)
  {$s^{\mathds{1}}(b)$};

\node[box, minimum width=\d cm] at ({\d/2}, \botY)
  {$c$};

\node[box, minimum width=\e cm] at ({\d+\e/2}, \botY)
  {$\overline{\sigma_i}\,\sigma_j$};

\node[box, minimum width=\f cm] at ({\d+\e+\f/2}, \botY)
  {$d$};

\end{tikzpicture}
        \item Case 2 : We have $b_1, b_2\in BW^{\mathds{1}}$ such that $b=b_1\overline{\sigma_i}\sigma_jb_2$, $c=s^{\mathds{1}}(a)\overline{\sigma_i}\sigma_js^{\mathds{1}}(b_1)$, and $d=s^{\mathds{1}}(b_2)$

\begin{tikzpicture}[
  box/.style={draw, minimum height=.5cm, align=center},
  every node/.style={inner sep=0pt}
]

\def\a{4.5}   
\def\b{2.0}   
\def\c{1.3}   
\def\cc{2.0}
\def\ccc{1.2}

\def\d{7.8}   
\def\e{2.0}   
\def\f{1.2}   

\def\topY{0}
\def\botY{-.7}

\node[box, minimum width=\a cm] at ({\a/2}, \topY)
  {$s^{\mathds{1}}(a)$};

\node[box, minimum width=\b cm] at ({\a+\b/2}, \topY)
  {$\overline{\sigma_i}\,\sigma_j$};

\node[box, minimum width=\c cm] (C) at ({\a+\b+\c/2}, \topY)
  {$s^{\mathds{1}}(b_1)$};

\node[box, minimum width=\cc cm] at ({\a+\b+\c+\cc/2}, \topY)
  {$\overline{\sigma_i}\,\sigma_j$};

\node[box, minimum width=\ccc cm] (E) at ({\a+\b+\c+\cc+\ccc/2}, \topY)
  {$s^{\mathds{1}}(b_2)$};

\node[box, minimum width=\d cm] at ({\d/2}, \botY)
  {$c$};

\node[box, minimum width=\e cm] at ({\d+\e/2}, \botY)
  {$\overline{\sigma_i}\,\sigma_j$};

\node[box, minimum width=\f cm] at ({\d+\e+\f/2}, \botY)
  {$d$};

\draw[decorate,decoration={brace,amplitude=6pt}]
  ([yshift=4pt] C.north west) -- ([yshift=4pt] E.north east)
  node[midway, above=8pt] {$s^{\mathds{1}}(b)$};
  
\end{tikzpicture}

\item Case 3 : We have $a_1, a_2\in BW^{\mathds{1}}$ such that $a=a_1\overline{\sigma_i}\sigma_ja_2$, $c=s^{\mathds{1}}(a_1)$, and $d=s^{\mathds{1}}(a_2)\overline{\sigma_i\sigma_j}s^{\mathds{1}}(b)$

\begin{tikzpicture}[
  box/.style={draw, minimum height=.5cm, align=center},
  every node/.style={inner sep=0pt}
]

\def\a{1.2}   
\def\b{2.0}   
\def\c{1.3}   
\def\cc{2.0}
\def\ccc{4.5}

\def\d{1.2}   
\def\e{2.0}   
\def\f{7.8}   

\def\topY{0}
\def\botY{-.7}

\node[box, minimum width=\a cm] (A) at ({\a/2}, \topY)
  {$s^{\mathds{1}}(a_1)$};

\node[box, minimum width=\b cm] at ({\a+\b/2}, \topY)
  {$\overline{\sigma_i}\,\sigma_j$};

\node[box, minimum width=\c cm] (C) at ({\a+\b+\c/2}, \topY)
  {$s^{\mathds{1}}(a_2)$};

\node[box, minimum width=\cc cm] at ({\a+\b+\c+\cc/2}, \topY)
  {$\overline{\sigma_i}\,\sigma_j$};

\node[box, minimum width=\ccc cm] at ({\a+\b+\c+\cc+\ccc/2}, \topY)
  {$s^{\mathds{1}}(b)$};

\node[box, minimum width=\d cm] at ({\d/2}, \botY)
  {$c$};

\node[box, minimum width=\e cm] at ({\d+\e/2}, \botY)
  {$\overline{\sigma_i}\,\sigma_j$};

\node[box, minimum width=\f cm] at ({\d+\e+\f/2}, \botY)
  {$d$};

\draw[decorate,decoration={brace,amplitude=6pt}]
  ([yshift=4pt] A.north west) -- ([yshift=4pt] C.north east)
  node[midway, above=8pt] {$s^{\mathds{1}}(a)$};
  
\end{tikzpicture}
    \end{itemize}
\end{lemma}
\begin{proof} (sketch)
    Because $\overline{\sigma_i}\neq\sigma_j$, we know $c$ cannot have length $s^{\mathds{1}}(a)\pm 1$. This reduces our possible cases from five to three. Within each case, Lemma \ref{lem:infix_of_irr} does the majority of the work. 
\end{proof}
\noindent With the necessary machinery set up, we prove a generalization of Theorem \ref{thm:new_rw}.

\begin{theorem}\label{thm:general_new_rw} \leanlink{https://github.com/hannahfechtner/braids_better/blob/098309403cda787eaf1fa954e07ff0ca3b2927ee/BraidProject/ReversingToGridStyle.lean#L115-L121}
    Given $a,b\in BW$, if $a \xrightarrow[R_{\mathrm{rev}}]{*}b$ then there exists a trivially-irreducible $g\in BW^{\mathds{1}}$ such that $i^{\mathds{1}}(a)\xrightarrow[R_{\mathrm{rev}}]{*}g$ and $s^{\mathds{1}}(g)=b$.
\end{theorem}
\begin{proof}
Assume $ a \xrightarrow[R_{\mathrm{rev}}]{*}b$. Thus, $a \underset{R_{\mathsf{rev}}}{\rightarrowtriangle} b$. We proceed by induction on that derivation.
\begin{itemize}
    \item Reflexivity : In this case, $a = b$. Let $g=i^{\mathds{1}}(b)$. Then we have $ i^{\mathds{1}}(a)\xrightarrow[R_{\mathrm{grid}}]{*}i^{\mathds{1}}(b)$ by reflexivity. By Lemma \ref{lem:i_irr}, $i^{\mathds{1}}(b)$ is trivially-irreducible, and by Lemma \ref{lem:s_after_i} $s^{\mathds{1}}\left( i^{\mathds{1}}(b) \right)= b$. 
    \item Transitive one-step reduction : In this case, we have $c,d, x, y\in BW$ such that $a \underset{R_{\mathsf{rev}}}{\rightarrowtriangle} c x d$, $R_{\mathsf{rev}}(x, y)$, and $b = c y d$. Our inductive hypothesis states that there is some trivially-irreducible $e\in BW^{\mathds{1}}$ such that $i^{\mathds{1}}(a) \underset{R_{\mathsf{grid}}}{\rightarrowtriangle} e$ and $s^{\mathds{1}}(e)=c x d$.

    As $R_{\mathsf{rev}}(x, y)$, $x=\overline{\sigma_i}\sigma_j$ for some $i$ and $j$; moreover, there is some $y'\in BW^{\mathds{1}}$ such that $R_{\mathsf{grid}}(x,y')$ holds by Theorem \ref{thm:rev_to_nontriv_relation}.

    Thus, $s^{\mathds{1}}(e)=c \overline{\sigma_i} \sigma_j d$. As $e$ is trivially-irreducible, it satisfies \lstinline{pairsTogether} by Theorem \ref{thm:pts_of_irr}. Hence $\overline{\sigma_i} \sigma_j$ is a substring of $e$. Let $e = w \overline{\sigma_i} \sigma_jt$. Substituting this in to the equation at the start of the paragraph,
    $s^{\mathds{1}}(w \overline{\sigma_i} \sigma_jt)=c \overline{\sigma_i} \sigma_j d$. Simplifying, we obtain
    $s^{\mathds{1}}(w) \overline{\sigma_i} \sigma_js^{\mathds{1}}(t)=c \overline{\sigma_i} \sigma_j d$.

    As $w$ and $t$ are both trivially-irreducible (being substrings of $e$), by Lemma \ref{lem:list_split} we have the following three cases:
    \begin{itemize}
        \item $s^{\mathds{1}}(w)= c$ and $s^{\mathds{1}}(t)=d$.  Let $g =$ \lstinline{move_ones (w ++ y' ++ t)}. We know 
$i^{\mathds{1}}(a) \underset{R_{\mathsf{grid}}}{\rightarrowtriangle} w x t$, so since $R_{\mathsf{grid}}(x,y')$, $i^{\mathds{1}}(a) \underset{R_{\mathsf{grid}}}{\rightarrowtriangle} w y' t$. Thus $i^{\mathds{1}}(a) \xrightarrow[R_{\mathrm{grid}}]{*} w y' t$. By Theorem \ref{thm:move_ones_equiv}, $w y' t\xrightarrow[R_{\mathrm{grid}}]{*}$\lstinline{ move_ones w ++ y' ++ t}. Then by transitivity, $i^{\mathds{1}}(a)\xrightarrow[R_{\mathrm{grid}}]{*}$\lstinline{ move_ones w ++ y' ++ t}. We thus have $i^{\mathds{1}}(a)\xrightarrow[R_{\mathrm{grid}}]{*} g$ as desired. We know that $g=$\lstinline{ move_ones w ++ y' ++ t} is trivially-irreducible by Theorem \ref{thm:move_ones_irr}. Lastly, 
\begin{align*}
    s^{\mathds{1}}(g) & = s^{\mathds{1}}(\text{\lstinline{move_ones w ++ y' ++ t}}) \\
    & = s^{\mathds{1}}(wy't) &\text{(by Theorem \ref{thm:move_ones_delete})} \\
    & = s^{\mathds{1}}(w)s^{\mathds{1}}(y')s^{\mathds{1}}(t) \\
    & = c \cdot s^{\mathds{1}}(y') \cdot d \\
    & = c y d & \text{(by Theorem \ref{thm:rev_to_nontriv_relation})}\\
    & = b
\end{align*} as desired

\item In this case, $t=t_1\overline{\sigma_i}\sigma_jt_2$, $c=s^{\mathds{1}}(w)\overline{\sigma_i}\sigma_js^{\mathds{1}}$, and $d=s^{\mathds{1}}(t_2)$ : 

We use \lstinline{move_ones (w ++ x ++ t1 ++ y' ++ t2)}; the proof proceeds from here analogously to that of the previous case.

\item The third and final case is entirely symmetrical with the previous

    \end{itemize}
    
\end{itemize}
\end{proof}

\noindent Finally, we give the proof of Theorem \ref{thm:new_rw}. Recall we must show that if $\overline{a}b \xrightarrow[R_{\mathrm{rev}}]{*}c\overline{d}$, then there exist $e, f \in BW^{\mathds{1}, +}$ such that $i^{\mathds{1}}(\overline{a})i^{\mathds{1}}(b)\xrightarrow[R_{\mathrm{grid}}]{*}e\overline{f}$ with $s^{\mathds{1}}(e)=c$ and $s^{\mathds{1}}(f)=d$.
\begin{proof} \leanlink{https://github.com/hannahfechtner/braids_better/blob/098309403cda787eaf1fa954e07ff0ca3b2927ee/BraidProject/ReversingToGridStyle.lean#L142-L152}
    Let $a, b, c, d \in BW^{+}$. Assume $\overline{a}b \xrightarrow[R_{\mathrm{rev}}]{*}c\overline{d}$. By Theorem \ref{thm:general_new_rw}, there exists some trivially-irreducible $g\in BW^{\mathds{1}}$ such that $i^{\mathds{1}}(\overline{a}b)\xrightarrow[R_{\mathrm{grid}}]{*}g$ and $s^{\mathds{1}}(g)=c\overline{d}$.
    By Lemma \ref{lem:pos_neg_irr}, there exist positive $e,f\in BW^{\mathds{1}}$ such that $g=e\overline{f}$. Thus, $i^{\mathds{1}}(\overline{a})i^{\mathds{1}}(b)\xrightarrow[R_{\mathrm{grid}}]{*}e\overline{f}$. 

    Moreover, note that $s^{\mathds{1}}(g)=s^{\mathds{1}}(e\overline{f})=s^{\mathds{1}}(e)s^{\mathds{1}}(\overline{f})$, and $s^{\mathds{1}}(e)$ is positive and $s^{\mathds{1}}(\overline{f})$ negative. Thus, as $s^{\mathds{1}}(e)s^{\mathds{1}}(\overline{f})=c\overline{d}$, we must have 
    $s^{\mathds{1}}(e)=c$ and $s^{\mathds{1}}(f)=d$ as desired.
\end{proof}

 We now know that we may produce an $R_{\mathsf{grid}}$--derivation from an $R_{\mathsf{rev}}$--derivation. We next connect this $R_{\mathsf{grid}}$--derivation to a grid under construction, towards our goal of building a grid $\langle a, b, c, d\rangle$.

\subsection{Partial Grids}
A partial grid generalizes the notion of a grid to include ``unfinished'' grids, which are missing cells in the bottom/right area. Let's use the following as a motivating example:

\begin{center}
\begin{tikzcd}
{} \arrow[rrrr, "\sigma_2"] \arrow[dddd, "\sigma_1"'] &                    &                                    &  & {} \arrow[dd, "\sigma_1"] \arrow[rr, "\sigma_2"] &                                             & {} \arrow[d, "\sigma_1"] \\
                                        &                    &                                    &  &                                    &                                             & {} \arrow[d, "\sigma_2"] \\
                                        &                    &                                    &  & {} \arrow[dd, "\sigma_2"] \arrow[r, "\sigma_2"]  & {} \arrow[r, "\sigma_1"] \arrow[dd, "\varepsilon"] & {}                \\
                                        &                    &                                    &  &                                    &                                             &                   \\
{} \arrow[rr, "\sigma_2"'] \arrow[dd, "\sigma_1"']    &                    & {} \arrow[rr, "\sigma_1"'] \arrow[d, "\sigma_1"] &  & {} \arrow[r, "\varepsilon"']        & {}                                          &                   \\
                                        &                    & {} \arrow[d, "\sigma_2"]                  &  &                                    &                                             &                   \\
{} \arrow[r, "\sigma_2"']                      & {} \arrow[r, "\sigma_1"'] & {}                                 &  &                                    &                                             &                  
\end{tikzcd}
\end{center}

\noindent Partial grids will be defined as an inductive subset of $(BW^\mathds{1})^5$. Before we give the formal definition, here is an intuitive definition to keep in mind. We define a partial grid only by its ``frame'' : the left side (read upwards against its direction), top side, the bottom frontier (the section that runs straight along the bottom), the middle frontier (the zigzag-like portion in the middle, read from southwest to northeast with downwards-pointing arrows read upwards against their direction), and the right frontier (the section which runs down the far-righthand side, read upwards against its direction). 

The frame of the above grid is given by the tuple $\langle \overline{\sigma_1}\,\overline{\sigma_1}, \sigma_2\sigma_2, \sigma_2\sigma_1, \overline{\sigma_2} \,\overline{\sigma_1} \sigma_1, \mathds{1}, \overline{\mathds{1}} \sigma_1, \overline{\sigma_2}\,\overline{\sigma_1}\rangle $, seen visually as
\begin{center}
\begin{tikzcd}
{} \arrow[dd, "\sigma _1\sigma _1"'] \arrow[rr, "\sigma_2\sigma_2"] &                     & {} \arrow[d, "\sigma_1\sigma_2"] \\
                                    &                     & {}                \\
{} \arrow[r, "\sigma_1\sigma_2"']                  & {} \arrow[ru, "\overline{\sigma_2}\,\overline{\sigma_1} \sigma_1 \mathds{1} \overline{\mathds{1}}\sigma_1 \overline{\sigma_2}\,\overline{\sigma_1}"', dash, squiggly] &                  
\end{tikzcd}
\end{center}
\noindent In general, we read a partial grid 
\begin{tikzcd}
{} \arrow[dd, "a"'] \arrow[rr, "b"] &                     & {} \arrow[d, "d"] \\
                                    &                     & {}                \\
{} \arrow[r, "c"']                  & {} \arrow[ru, "m"', dash, squiggly] &                  
\end{tikzcd}
as $\langle \overline{a}, b, c, m, \overline{d}\rangle$.

\noindent We draw a partial grid $\langle a, b, c, m, d\rangle$ as
\begin{tikzcd}
{} \arrow[dd, "\overline{a}"'] \arrow[rr, "b"] &                     & {} \arrow[d, "\overline{d}"] \\
                                    &                     & {}                \\
{} \arrow[r, "c"']                  & {} \arrow[ru, "m"', dash, squiggly] &   
\end{tikzcd}. We will soon show that both $a$ and $d$ must be negative augmented braid words, so $\overline{a}$ and $\overline{d}$ are positive augmented braid words (no generators with overlines appear in them). 


Often, we speak of the (total) frontier of a partial grid; this is the concatenation of the bottom frontier, middle frontier, and right frontier. In our example, it would be $\sigma_1\sigma_2\cdot \overline{\sigma_2}\,\overline{\sigma_1} \sigma_1 \mathds{1} \overline{\mathds{1}}\sigma_1 \overline{\sigma_2}\,\overline{\sigma_1}\cdot \overline{\sigma_2}\, \overline{\sigma_1} $.  

We now give a formal, inductive definition of the set of partial grids $G_P \subset (BW^{\mathds{1}})^5 $. We provide diagrams to guide the reader's intuition; they do not form part of the definition.
\begin{definition} \label{def:partial_grid} \leanlink{https://github.com/hannahfechtner/braids_better/blob/098309403cda787eaf1fa954e07ff0ca3b2927ee/BraidProject/PartialGrid/Basic.lean#L9-L30}
The set of partial grids $G_P$ is the smallest subset of $ (BW^{\mathds{1}})^5 $ satisfying
\begin {itemize}
\item \textbf{(a cell)} : For any $a,b,c,d\in BW^{\mathds{1}, +}$, if $(\overline{a}\,b, \;c\,\overline{d})\in R_{\mathsf{grid}}$, we have $\langle \overline{a},\; b,\; c, \; \varepsilon,\; \overline{d}\rangle \in G_P$

Visually, \begin{tikzcd}
{} \arrow[r, "b"] \arrow[d, "a"'] & {} \arrow[d, "d"] \\
{} \arrow[r, "c"']                & {}               
\end{tikzcd}

\item \textbf{(a spine)} : For any $a,b\in BW^{\mathds{1},+}\,\backslash \, \{\varepsilon\}$, $\langle \overline{a},\; b,\; \varepsilon, \; \overline{a}\, b,\; \varepsilon\rangle \in G_P$ 

Visually,
\begin{tikzcd}
{} \arrow[d, "a"'] \arrow[r, "b"] & {} \\
{}                                &   
\end{tikzcd} for non-empty $a$ and $b$

\item \textbf{(horizontal concatenation to a partial grid with empty middle frontier)} :

\begin{tikzcd}
{} \arrow[rr, "b_1"] \arrow[dd, "\overline{a_1}"'] &  & {} \arrow[dd, "\overline{d_1}"] \arrow[rr, "b_2"] &                     & {} \arrow[d, "\overline{d_2}"] \\
                                    &  &                                    &                     & {}                \\
{} \arrow[rr, "c_1"']                 &  & {} \arrow[r, "c_2"']                 & {} \arrow[ru, "m_2"', dash, squiggly] &                  
\end{tikzcd}
becomes
\begin{tikzcd}
{} \arrow[dd, "\overline{a_1}"'] \arrow[rrrr, "b_1\cdot b_2"] &  &  &                       & {} \arrow[d, "\overline{d_2}"] \\
                                                   &  &  &                       & {}                  \\
{} \arrow[rrr, "c_1\cdot c_2"']                    &  &  & {} \arrow[ru, "m_2"', dash, squiggly] &                    
\end{tikzcd}

Given $\langle a_1, b_1, c_1, \varepsilon, d_1\rangle \in G_P$ and $\langle d_1, b_2, c_2, m_2, d_2\rangle \in G_P$, we have $\langle a_1, b_1\cdot b_2, c_1\cdot c_2, m_2, d_2\rangle \in G_P$.

\item \textbf{(horizontal concatenation to a partial grid with non-empty middle frontier)} ($m_1$ is non-empty) : 

\begin{tikzcd}
{} \arrow[ddd, "\overline{a_1}"'] \arrow[rr, "b_1"] &                     & {} \arrow[rr, "b_2"] \arrow[dd, "\overline{d_1}"] &                     & {} \arrow[d, "\overline{d_2}"] \\
                                     &                     &                                    &                     & {}                \\
                                     &                     & {} \arrow[r, "c_2"']                 & {} \arrow[ru, "m_2"', dash, squiggly] &                   \\
{} \arrow[r, "c_1"']                   & {} \arrow[ru, "m_1"', dash, squiggly] &                                    &                     &                  
\end{tikzcd}
 becomes
\begin{tikzcd}
{} \arrow[ddd, "\overline{a_1}"'] \arrow[rrrr, "b_1 \cdot b_2"] &                                              &  &  & {} \arrow[d, "\overline{d_2}"] \\
                                                   &                                              &  &  & {}                  \\
                                                   &                                              &  &  &                     \\
{} \arrow[r, "c_1"']                               & {} \arrow[rrruu, "m_1 \cdot c_2 \cdot m_2"', dash, squiggly] &  &  &                    
\end{tikzcd}

Given $\langle a_1, b_1, c_1, m_1, d_1\rangle \in G_P$ and $\langle d_1, b_2, c_2, m_2, d_2\rangle \in G_P$, we have $\langle a_1, b_1\cdot b_2, c_1, m_1 \cdot c_2 \cdot m_2, d_2\rangle \in G_P$.

\item \textbf{(vertical concatenation to a partial grid with empty middle frontier)}  :

\begin{tikzcd}
{} \arrow[rr, "b_1"] \arrow[dd, "\overline{a_1}"'] &                     & {} \arrow[dd, "\overline{d_1}"] \\
                                    &                     &                    \\
{} \arrow[rr, "c_1"] \arrow[dd, "\overline{a_2}"'] &                     & {} \arrow[d, "\overline{d_2}"]  \\
                                    &                     & {}                 \\
{} \arrow[r, "c_2"']                  & {} \arrow[ru, "m_2"', dash, squiggly] &                   
\end{tikzcd}
becomes
\begin{tikzcd}
{} \arrow[rr, "b_1"] \arrow[dddd, "\overline{a_1}\cdot\overline{a_2} = \overline{a_2\cdot a_1}"'] &                       & {} \arrow[ddd, "\overline{d_2 \cdot d_1}"] \\
                                                   &                       &                                 \\
                                                   &                       &                                 \\
                                                   &                       & {}                              \\
{} \arrow[r, "c_2"']                               & {} \arrow[ru, "m_2"'] &                                
\end{tikzcd}

Given $\langle a_1, b_1, c_1, \varepsilon, d_1\rangle \in G_P$ and $\langle a_2, c_1, c_2, m_2, d_2\rangle \in G_P$, we have $\langle a_2\cdot a_1, b_1, c_2, m_2,d_2\cdot d_1\rangle \in G_P$. Note the order of the first element and the final element : $a_2\cdot a_1$ and $d_2\cdot d_1$. We are reading \textit{against} the direction of the arrow on the vertical side of the partial grid, so the order is flipped.

\item \textbf{(vertical concatenation to a partial grid with non-empty middle frontier)} 

($m_1$ is non-empty) :

\begin{tikzcd}
{} \arrow[dd, "\overline{a_1}"'] \arrow[rrr, "b_1"] &                     &                                    & {} \arrow[d, "\overline{d_1}"] \\
                                     &                     &                                    & {}                \\
{} \arrow[rr, "c_1"] \arrow[dd, "\overline{a_2}"']  &                     & {} \arrow[d, "\overline{d_2}"] \arrow[ru, "m_1"', dash, squiggly] &                   \\
                                     &                     & {}                                 &                   \\
{} \arrow[r, "c_2"']                   & {} \arrow[ru, "m_2"', dash, squiggly] &                                    &                  
\end{tikzcd}
becomes
\begin{tikzcd}
{} \arrow[dddd, "\overline{a_2\cdot a_1}"'] \arrow[rrr, "b_1"] &                                            &  & {} \arrow[d, "\overline{d_1}"] \\
                                         &                                            &  & {}                  \\
{}                    &                                            &  &                     \\
                                         &                                            &  &                     \\
{} \arrow[r, "c_2"']                     & {} \arrow[rruuu, "m_2\cdot d_2\cdot m_1"', dash, squiggly] &  &                    
\end{tikzcd}

Given $\langle a_1, b_1, c_1, m_1, d_1\rangle \in G_P$ and $\langle a_2, c_1, c_2, m_2, d_2\rangle \in G_P$, we have $\langle a_2\cdot a_1, b_1, c_2, m_2 \cdot d_2 \cdot m_1, d_1\rangle \in G_P$.

\end{itemize}
\end{definition}
\begin{definition}
    The \textbf{total frontier} (often referred to as merely ``the frontier'') of a partial grid $\langle a, b, c, m, d\rangle \in G_P$ is $c\cdot m\cdot d$.
\end{definition}
    
We continue with three short lemmas about partial grids, all shown by structural induction. While simple in theory, the number of cases grows rapidly.

\begin{lemma} \label{lem:sides_pos_neg} \leanlink{https://github.com/hannahfechtner/braids_better/blob/098309403cda787eaf1fa954e07ff0ca3b2927ee/BraidProject/PartialGrid/Basic.lean#L125-L159}
    Consider some partial grid $\langle a, b, c, m, d\rangle\in G_P$. The left side ($a$) and right frontier ($d$) must be negative augmented braid words; the top side ($b$) and bottom frontier ($c$) must be positive augmented braid words. 
\end{lemma}

\begin{lemma} \label{lem:extend_frontier}\leanlink{https://github.com/hannahfechtner/braids_better/blob/098309403cda787eaf1fa954e07ff0ca3b2927ee/BraidProject/PartialGrid/Basic.lean#L187-L245} (Extending a partial grid's spine)
    Consider a partial grid $\langle a, b, c, m, d\rangle$. For any non-empty \textit{negative} augmented grid word $a'$ there exists a partial grid $\langle a'\cdot a, b, \varepsilon, a'\cdot c\cdot m, d\rangle$. Visually, 
\begin{center}    
\begin{tikzcd}
{} \arrow[dd, "a"'] \arrow[rr, "b"] &                     & {} \arrow[d, "d"] \\
                                    &                     & {}                \\
{} \arrow[r, "c"'] \arrow[d, "a'"'] & {} \arrow[ru, "m"', dash, squiggly] &                   \\
{}                                  &                     &                  
\end{tikzcd} 
\end{center}
For any non-empty \textit{positive} augmented grid word $b'$, there exists a partial grid $\langle a, b\cdot b', c, m \cdot d \cdot b', \varepsilon \rangle$. Visually,
\begin{center}
\begin{tikzcd}
{} \arrow[dd, "a"'] \arrow[rr, "b"] &                     & {} \arrow[d, "d"] \arrow[r, "b'"] & {} \\
                                    &                     & {}                                &    \\
{} \arrow[r, "c"']                  & {} \arrow[ru, "m"', dash, squiggly ] &                                   &   
\end{tikzcd}
\end{center}
\end{lemma}

\begin{lemma} \label{lem:frontierup} \leanlink{https://github.com/hannahfechtner/braids_better/blob/098309403cda787eaf1fa954e07ff0ca3b2927ee/BraidProject/PartialGrid/Basic.lean#L309-L391}
    The middle frontier of any partial grid is either empty or of the form $\overline{\alpha}m'\beta$ for some letters $\alpha, \beta\in S^{\mathds{1}}$ and some $m'\in BW^{\mathds{1}}$. Note that this is only true because we keep track of empty arrows in partial grids using $BW^{\mathds{1}}$.
\end{lemma}

The following two theorems are not used in the proof of Theorem \ref{thm:big_one}, but are basic properties of partial grids, so we present them now.

First, we can reflect a partial grid along its top-left to bottom-right diagonal.
\begin{theorem}\label{thm:mirror_partial_grid} (Reflecting) \leanlink{https://github.com/hannahfechtner/braids_better/blob/098309403cda787eaf1fa954e07ff0ca3b2927ee/BraidProject/PartialGrid/Basic.lean#L51-L107}
    Given a partial grid $\langle a, b, c, m, d\rangle$, there is a partial grid $\langle \overline b, \overline a, \overline d, \overline m, \overline c,\rangle$
\begin{center}
\begin{tikzcd}
{} \arrow[dd, "a"'] \arrow[rr, "b"] \arrow[rrdd, dashed, dash, cyan, shorten <=-17 pt, shorten >=-8pt] &                     & {} \arrow[d, "d"] \\
                                                         &                     & {}                \\
{} \arrow[r, "c"']                                       & {} \arrow[ru, "m"'{yshift=-3pt}, dash, squiggly] & {}               
\end{tikzcd}
$\implies$
\begin{tikzcd}
{} \arrow[dd, "\overline{b}"'] \arrow[rr, "\overline a"] &                     & {} \arrow[d, "\overline c"] \\
                                    &                     & {}                \\
{} \arrow[r, "\overline d"']                  & {} \arrow[ru, "\overline m"', dash, squiggly] &                  
\end{tikzcd}
\end{center}
\end{theorem}

This theorem shortens proofs significantly; for example, the second half of Lemma \ref{lem:extend_frontier} becomes an immediate consequence of the first half. We also use this mirroring trick in the following theorem about splitting partial grids: 

\begin{theorem}\label{thm:splittability} (Vertical Splitting) \leanlink{https://github.com/hannahfechtner/braids_better/blob/098309403cda787eaf1fa954e07ff0ca3b2927ee/BraidProject/PartialGrid/Splittability.lean#L11-L396}
    Given $p = \langle a, b_1b_2, c, m, d\rangle\in G_P$ with nonempty $b_1, b_2$, we have either
    \begin{itemize}
        \item $p$ splits into two smaller partial grids $p_1 = \langle a, b_1, c_1, m_1, d_1\rangle\in G_P$ and $p_2 = \langle d_1, b_2, c_2, m_2, d\rangle\in G_P$, with $c\cdot m = c_1 \cdot m_1 \cdot c_2 \cdot m_2$

\begin{tikzcd}
{} \arrow[dd, "a"'] \arrow[r, "b_1"] 
& {} \arrow[r, "b_2"] 
& {} \arrow[d, "d"] \\
& p & {} \\
{} \arrow[r, "c"', preaction={draw=magenta, draw opacity=.35, line width=6pt, -, shorten >=-4pt}] 
& {} \arrow[ru, "m"', preaction={draw=magenta, draw opacity=.35, line width=6pt, shorten <=-4pt}, dash, squiggly] 
&
\end{tikzcd} splits into
\begin{tikzcd}
{} \arrow[ddd, "a"'] \arrow[rr, "b_1"] 
&                     
& {} \arrow[rr, "b_2"] \arrow[dd, "d_1"] 
&                     
& {} \arrow[d, "d"] \\
& {p_1} 
&                                    
& {p_2} 
& {}                \\
&                     
& {} \arrow[r, "c_2"', preaction={draw=magenta, draw opacity=.35, line width=6pt, -, shorten >=-4pt, shorten <= -4pt}]                 
& {} \arrow[ru, "m_2"', dash, squiggly, preaction={draw=magenta, draw opacity=.35, line width=6pt, shorten <=-4pt}] 
&                   \\
{} \arrow[r, "c_1"', preaction={draw=magenta, draw opacity=.35, line width=6pt, -, shorten >=-4pt}]                   
& {} \arrow[ru, "m_1"', dash, squiggly, preaction={draw=magenta, draw opacity=.35, line width=6pt, shorten <=-4pt, shorten >= -7pt}] 
&                                    
&                     
&                  
\end{tikzcd}

(note that $c$ is not necessarily equal to $c_1$, as $m_1$ could be empty, making $c=c_1c_2$. $d$ may be either empty or nonempty)

\item $p$ splits into a smaller partial grid $p_1 = \langle a, b_1, c, m_1, \varepsilon\rangle\in G_P$ and a tail \begin{tikzcd}
{} \arrow[r, "b_2"] & {}
\end{tikzcd}, with $m =m_1\cdot b_2$. In this case, we must have $d=\varepsilon$, and so $m\cdot d = m_1\cdot b_2$

        \begin{tikzcd}
{} \arrow[dd, "a"'] \arrow[r, "b_1"] 
& {} \arrow[r, "b_2"] 
& {} \arrow[d, "d",preaction={draw=ProcessBlue, draw opacity=.35, line width=6pt, -, shorten <=-8pt}] \\
& p & {} \\
{} \arrow[r, "c"'] 
& {} \arrow[ru, "m"', preaction={draw=ProcessBlue, draw opacity=.35, line width=6pt, shorten >=-7pt}, dash, squiggly] 
&
\end{tikzcd}
splits into
\begin{tikzcd}
{} \arrow[dd, "a"'] \arrow[rr, "b_1"] &                        & {} \arrow[r, "b_2", preaction={draw=ProcessBlue, draw opacity=.35, line width=6pt, shorten <=-6pt, -}] & {} \\
                                     &     p_1                   &                     &    \\
{} \arrow[r, "c"']                    & {} \arrow[ruu, "m_1"', dash, squiggly, preaction={draw=ProcessBlue, draw opacity=.35, line width=6pt, shorten >= -7pt}] &                     &   
\end{tikzcd}
\end{itemize}
The proof is by induction on $p$.

\vspace{.6em}
(Horizontal Splitting) \leanlink{https://github.com/hannahfechtner/braids_better/blob/098309403cda787eaf1fa954e07ff0ca3b2927ee/BraidProject/PartialGrid/Splittability.lean#L360-L396} We have the analogous property for a partial grid $p' = \langle a_1a_2, b, c, m, d\rangle$ with nonempty $a_1, a_2$. The proof reflects $p'$, splits it vertically, and reflects the result(s).
\end{theorem}

\subsection{Building the Partial Grid}

We next build a partial grid from an $R_{\mathsf{grid}}$--derivation. We may only do so when the derivation begins with a word in negative-positive form, as that will correspond to the spine of a partial grid. 
\label{rwpgrid}
\begin{theorem}\label{thm:partial_grid_of_rw} \leanlink{https://github.com/hannahfechtner/braids_better/blob/098309403cda787eaf1fa954e07ff0ca3b2927ee/BraidProject/PartialGrid/Build.lean#L10-L26}
    For nonempty $a,b\in BW^{\mathds{1}, +}$, if $\overline{a}b \; \xrightarrow[R_{\mathrm{grid}}]{*} c$, then there is a partial grid with spine $(a, b)$ and frontier $c$.
\end{theorem}

We begin by defining a function called \lstinline{add_cell} which adds a single cell to a partial grid. 

\begin{definition} \leanlink{https://github.com/hannahfechtner/braids_better/blob/098309403cda787eaf1fa954e07ff0ca3b2927ee/BraidProject/PartialGrid/AddCell.lean#L150-L427}
    \lstinline{add_cell} is a function taking in some $(x,y)\in R_{\mathsf{grid}}$ and a partial grid $P_1 =\langle a, b, c, m_1\cdot x\cdot m_2,d\rangle$, and returning a partial grid $P_2=\langle a, b, c', m',d'\rangle$ such that $c'\cdot m'\cdot d'=c\cdot m_1\cdot y \cdot m_2\cdot d$, $c$ is a prefix of $c'$, and $d$ is a suffix of $d'$. The final two specifications are used only internally to the definition. By contrast, the first of the specifications (that $c'\cdot m'\cdot d'=c\cdot m_1\cdot y \cdot m_2\cdot d$) is both necessary for the definition and used externally for other theorems later in the project.

    We define \lstinline{add_cell} by recursion on $P_1$. 
    \begin{itemize}
        \item \textbf{a cell} : This case would lead to a contradiction. Partial grids from the cell case have an empty middle frontier. As $(x, y)\in R_{\mathsf{grid}}$, $x$ is not the empty string, so $P_1$ has a nonempty middle frontier. 
        \item \textbf{a spine} : Here $c=d=\varepsilon$ and $m_1\cdot x \cdot m_2=ab$. 
    First, note that as $(x, y)\in R_{\mathsf{grid}}$, there exist letters $\alpha, \beta, \gamma, \delta\in S^{\mathds{1}}$ such that $x = \overline{\alpha}\beta$ and $y=\gamma\overline{\delta}$. Thus $m_1\cdot \overline{\alpha}\beta\cdot m_2=ab$. As $a$ is negative and $b$ is positive (by Lemma \ref{lem:sides_pos_neg}), we must have $a=m_1\overline{\alpha}$ and $b=\beta m_2$. 
        As $(\overline{\alpha}\beta, \gamma\overline{\delta})\in R_{\mathsf{grid}}$, we have $P_3=\langle\overline{\alpha}, \beta, \gamma, \varepsilon, \overline{\delta}\rangle\in G_P$. 
        Consider the case where neither $m_1$ nor $m_2$ is the empty string (we omit the other, simpler cases).
        By Lemma \ref{lem:extend_frontier}, we may extend the left side of $P_3$ to obtain a partial grid $P_4=\langle m_1\overline{\alpha}, \beta, \varepsilon, m_1\gamma, \overline{\delta}\rangle$, and then the right side of $P_4$ to obtain a partial grid $P_5=\langle m_1\overline{\alpha}, \beta m_2, \varepsilon, m_1\gamma \overline{\delta}m_2, \varepsilon \rangle$. Hence, we have a partial grid $\langle a, b, c, m_1\cdot y\cdot m_2, d\rangle$ as desired. The condition on the equality of the total frontier immediately holds. As $c=d=\varepsilon$ is both a prefix and a suffix of any word, those two conditions are also fulfilled.
        \begin{center}
\begin{tikzcd}
{} \arrow[r, "\beta "] \arrow[d, "\alpha"'] & {} \arrow[d, "\delta"] \arrow[r, "m_2"] & {} \\
{} \arrow[r, "\gamma"'] \arrow[d, "m_1"']   & {}                                      &    \\
{}                                          &                                         &   
\end{tikzcd}
\end{center}
        \item \textbf{horizontal concatenation to a partial grid with non-empty middle frontier} : 
        

In this case, we have $b=b_1b_2$, partial grids $P_2=\langle a, b_1, c, n_1, d_1\rangle$ and $P_3=\langle d_1, b_2, c_2, n_2, d\rangle$, with $ n_1 \cdot c_2 \cdot n_2 = m_1\cdot x\cdot m_2$. By properties of middle frontiers (Lemma \ref{lem:frontierup}), we know that either $x$ is an infix of $n_1$ or of $n_2$. Consider the first case, where we have $n_1=n_3xn_4$, $m_1=n_3$ and $m_2 = n_4\cdot c_2\cdot n_2$; the other case is analogous. 

\begin{center}
\begin{tikzcd}
{} \arrow[ddddd, "a"'] \arrow[rrrr, "b_1"] &       &       &      & {} \arrow[dd, "d_1"] \arrow[rr, "b_2"] &      & {} \arrow[d, "d"] \\
       &        &         &       &       &  P_3    & {}       \\
     &      &    P_2   &   & {} \arrow[r, "c_2"']   & {} \arrow[ru, dash, squiggly, "n_2"] &     \\
     &         &       & {} \arrow[ru,dash, squiggly, "n_4"'] &        &       &           \\
     &          & {} \arrow[ru, Mulberry, "y"{pos=1, below, yshift=-22pt}, dash, to path={ -| (\tikztotarget) \tikztonodes }]
  \arrow[ru, "x"{pos=0, above, yshift=25pt}, dash, to path={ |- (\tikztotarget) \tikztonodes }] &       &             &         &            \\
{} \arrow[r, "c"']     & {} \arrow[ru, dash, squiggly,  "n_3"'] &      &      &       &     &                    
\end{tikzcd}
\end{center}
Here, recursively call \lstinline{add_cell} on $P_2$ to obtain a partial grid $P_4 = \langle a, b_1, c', n', d_1'\rangle$ with $c'\cdot n'\cdot d_1'=c\cdot n_3\cdot y\cdot n_4 \cdot d_1$, $c$ a prefix of $c'$ and $d_1$ a suffix of $d_1'$. Since $d_1$ is a suffix of $d_1'$, we know that there exists some $d_2$ such that $d_1'=d_1d_2$. We consider the case in which $d_2$ and $n'$ are nonempty (the other cases are simpler). 

\begin{center}
\begin{tikzcd}
{} \arrow[dddd, "a"'] \arrow[rr, "b_1"] &                      & {} \arrow[dd, "d_1"] \\
                                        &                      &                      \\
                                        & P_4                  & {} \arrow[d, "d_2"]  \\
                                        &                      & {}                   \\
{} \arrow[r, "c'"']                     & {} \arrow[ru, "n'"'] &                     
\end{tikzcd}
\end{center}
We extend the left side of $P_3$ to obtain a partial grid $P_5=\langle d_1', b_2, \varepsilon, d_2c_2n_2, d\rangle$. 
\begin{center}
\begin{tikzcd}
{} \arrow[dd, "d_1"'] \arrow[rr, "b_2"] &                       & {} \arrow[d, "d"] \\
                                        &                       & {}                \\
{} \arrow[r, "c_2"'] \arrow[d, "d_2"']  & {} \arrow[ru, "n_2"'] &                   \\
{}                                      &                       &                  
\end{tikzcd}
\end{center}
Next append $P_4$ and $P_5$ to obtain the desired partial grid $ P_6 = \langle a, b_1b_2, c', n'd_2c_2n_2, d\rangle$. 

\begin{center}
\begin{tikzcd}
{} \arrow[dddd, "a"'] \arrow[rr, "b_1"] &                      & {} \arrow[dd, color = gray, "d_1"] \arrow[rr, "b_2"] &                       & {} \arrow[d, "d"] \\
                                        &                      &                                        &                       & {}                \\
                                        &                   & {} \arrow[d, "d_2"] \arrow[r, "c_2"']  & {} \arrow[ru, "n_2"'] &                   \\
                                        &                      & {}                                     &                       &                   \\
{} \arrow[r, "c'"']                     & {} \arrow[ru, "n'"'] &                                        &                       &                  
\end{tikzcd}
\end{center}

To show that the first specification holds, we must show $c'\cdot n'\cdot d_2\cdot c_2\cdot n_2\cdot d = c\cdot m_1\cdot y\cdot m_2 \cdot d$. Recall that our inductive hypothesis gave us $c'\cdot n'\cdot d_1'=c\cdot n_3\cdot y\cdot n_4 \cdot d_1$. As $d_1'=d_2\cdot d_1$, we have
$c'\cdot n'\cdot d_2=c\cdot n_3\cdot y\cdot n_4$. Appending $c_2\cdot n_2\cdot d$ to both sides, we have
$c'\cdot n'\cdot d_2\cdot c_2\cdot n_2\cdot d=c\cdot n_3\cdot y\cdot n_4\cdot c_2\cdot n_2\cdot d$. In the case we have selected, we noted that $m_1=n_3$ and $m_2 = n_4\cdot c_2\cdot n_2$. Thus, 
$c'\cdot n'\cdot d_2\cdot c_2\cdot n_2\cdot d=c\cdot m_1\cdot y\cdot m_2 \cdot d$, as desired.
For the final two conditions, the inductive hypothesis states that $c$ is a prefix of $c'$, and $d$ is immediately a suffix of $d$.
        \item The other horizontal case and the two vertical cases proceed similarly        
    \end{itemize}
\end{definition}

\noindent We are now prepared to prove Theorem \ref{thm:partial_grid_of_rw}. 
\begin{proof} \leanlink{https://github.com/hannahfechtner/braids_better/blob/098309403cda787eaf1fa954e07ff0ca3b2927ee/BraidProject/PartialGrid/Build.lean#L14-L26}
        Consider nonempty, positive $a,b\in BW^{\mathds{1}}$. Assume $\overline{a}b \; \xrightarrow[R_{\mathrm{grid}}]{*} c$. Thus 
        $\overline{a}b \underset{R_{\mathsf{grid}}}{\rightarrowtriangle} c$. We proceed by induction on this derivation.

        \begin{itemize}
            \item \textbf{reflexivity :} In this case, we have $c = \overline{a}b$. Then we use the spinal partial grid given by $a$ and $b$, $\langle \overline{a}, b, \varepsilon, \overline{a}b, \varepsilon\rangle$. 
            \item \textbf{transitive one-step reduction :} We assume that $\overline{a} b \underset{R_{\mathsf{rev}}}{\rightarrowtriangle} d x e$, $R_{\mathsf{rev}}(x, y)$, and $c = d y e$. Our inductive hypothesis states that there is a partial grid $P_1:=\langle\overline{a}, b, c_1, m_1, d_1\rangle$ with $c_1m_1d_1=dxe$. We note that $x$ must be an infix of $m_1$, as $x$ is in negative-positive form, while $c_1$ is positive and $d_1$ negative. By the same reasoning, there exist $m_2, m_3\in BW^{\mathds{1}}$ such that $m_1=m_2xm_3$, $d=c_1m_2$, and $e=m_3d_1$. 
            We call \lstinline{add_cell} with $(x,y)$ and $P_1$ to obtain $P_2 :=\langle\overline{a}, b, c_1, m_2ym_3, d_1\rangle$. Then $P_2$'s frontier is $c_1m_2ym_3d_1=dye=c$ as desired.
        \end{itemize}
\end{proof}
\subsection{Empty Middle Frontiers}

Finally, we connect partial grids to grids by proving that a partial grid with empty middle frontier is in fact a grid. 

\begin{theorem}\label{thm:grid_of_pg} \leanlink{https://github.com/hannahfechtner/braids_better/blob/098309403cda787eaf1fa954e07ff0ca3b2927ee/BraidProject/PartialGrid/ToGrid.lean#L82-L84}
    Given $\langle \overline{a}, b, c, \varepsilon, \overline{d}\rangle \in G_P$, we have $\langle s^{\mathds{1}}(a), s^{\mathds{1}}(b), s^{\mathds{1}}(c), s^{\mathds{1}}(d)\rangle\in G$. 
\end{theorem}

\noindent The proof is by induction on the given partial grid. 

We next note that any partial grid whose total frontier is in positive-negative form has an empty middle frontier : 

\begin{lemma}\label{lem:square_fontier} \leanlink{https://github.com/hannahfechtner/braids_better/blob/098309403cda787eaf1fa954e07ff0ca3b2927ee/BraidProject/PartialGrid/Basic.lean#L393-L408}
    If a partial grid has total frontier $c\overline{d}$ for some $c, d\in BW^{\mathds{1},+}$, then it has bottom frontier $c$, middle frontier $\varepsilon$, and right frontier $\overline{d}$. (This is a consequence of Lemma \ref{lem:frontierup})
\end{lemma}

\subsection{Conclusion}
We are now prepared to give a proof of the remaining half of Theorem \ref{thm:big_one}, restated as follows : 
\begin{theorem} \label{big_one_one_dir} \leanlink{https://github.com/hannahfechtner/braids_better/blob/098309403cda787eaf1fa954e07ff0ca3b2927ee/BraidProject/SemiThueToGrid.lean#L10-L71}
    For any $a, b, c, d\in BW^+$, if $\overline{a}b\xrightarrow[R_{\mathrm{rev}}]{*}c\overline{d}$, then we have a grid $\langle a, b, c, d\rangle$.
\end{theorem}

\begin{proof}
Consider $a, b, c, d\in BW^+$. Assume $\overline{a}b\xrightarrow[R_{\mathrm{rev}}]{*}c\overline{d}$. By Theorem \ref{thm:general_new_rw}, we have $e, f\in BW^{\mathds{1}, +}$ such that $i^{\mathds{}1}(\overline{a})i^{\mathds{1}}(b)\xrightarrow[R_{\mathrm{grid}}]{*}e\overline{f}$ with $s^{\mathds{1}}(e)=c$ and $s^{\mathds{1}}(f)=d$. By Theorem \ref{thm:partial_grid_of_rw}, we have a partial grid $P_1$ with spine $(i^{\mathds{1}}(\overline{a}),\,i^{\mathds{1}}(b))$ and total frontier $e\overline{f}$. By Lemma \ref{lem:square_fontier}, we have $P_1=\langle i^{\mathds{1}}(\overline{a}),i^{\mathds{1}}(b),e,\varepsilon, \overline{f}\rangle$. Finally, by Theorem \ref{thm:grid_of_pg}, we have a grid $\Big\langle s^{\mathds{1}}\big(i^{\mathds{1}}(a)\big), s^{\mathds{1}}\big(i^{\mathds{1}}(b)\big), s^{\mathds{1}}(e), s^{\mathds{1}}(f))\Big\rangle = \langle a, b, c, d\rangle$ as required.  

\end{proof}

Theorem \ref{thm:big_one} can be used to give new, visual proofs of previous results. For example, note that properties of cancellativity and common multiples (the Ore conditions) are used for existence and uniqueness of grids, and hence for the proof of Theorem \ref{thm:big_one}. Thus, grids essentially capture Ore localization for a small subset of braid words (for those pairs of braid words which are the source and target of some grid).  We thus have a slick proof of the following : 

\begin{theorem} \leanlink{https://github.com/hannahfechtner/braids_better/blob/098309403cda787eaf1fa954e07ff0ca3b2927ee/BraidProject/SemiThueToGrid.lean#L73-L77}
    For any $a, b, c, d\in BW^+$, if $\overline{a}b\xrightarrow[R_{\mathrm{rev}}]{*}c\overline{d}$, then $ac\equiv^+bd$.
\end{theorem} 

\begin{proof}
    Assume $\overline{a}b\xrightarrow[R_{\mathrm{rev}}]{*}c\overline{d}$. Then by Theorem \ref{thm:big_one}, we have $\langle a, b, c, d\rangle \in G$. Hence by \ref{thm:grid_equiv}, $ac\equiv^+bd$. 
\end{proof}

Alternatively, we could have noted that since $\overline{a}b\xrightarrow[R_{\mathrm{rev}}]{*}c\overline{d}$, $\overline{a}b\equiv c\overline{d}$. In other notation, $a^{-1}b\equiv cd^{-1}$, and thus $ac \equiv bd$. As $ac$ and $bd$ are both positive braid words, by Ore localization (Theorem \ref{thm:braid_group_to_monoid_equiv_of_pos}), we know $ac\equiv^+bd$ as desired. 

\section{Metrics}

This chapter introduces metrics on grids, partial grids, and Semi-Thue derivations. In the Lean formalization, doing so requires us to distinguish between two ways in which an inductively defined mathematical object may be given : as a proposition asserting that the object exists, or as data recording how the object was constructed. A proposition in Lean has type \lstinline{Prop}, while mathematical data lives in \lstinline{Type}. To disambiguate the two, we begin by giving simple examples of \lstinline{Prop} and \lstinline{Type}.

Suppose that we have a predicate  $p : \mathbb N\to\mathsf{Prop}$. The proposition ``There exists an $n$ such that $p(n)$ holds'' is written as
\lstinline{∃ n : ℕ, p n}, which
has type \lstinline{Prop}. A proof of this proposition establishes that some natural number satisfying $p$ exists. We can discuss an arbitrary such number using \lstinline{Classical.choose} to give it a name, and \lstinline{Classical.choose_spec} to assert that $p$ holds of said named number. This selection uses the axiom of choice; \lstinline{Classical.choose} establishes that a suitable $n$ exists, but does not provide a computable procedure for finding it. Thus, a definition whose value depends on \lstinline{Classical.choose} is marked \lstinline{noncomputable}.

The corresponding \lstinline{Type}-valued concept is the subtype
\lstinline|{n : ℕ // p n}|.
An element of this subtype, say \lstinline|x : {n : ℕ // p n}|, consists of a particular natural number (\lstinline{x.val}) with which we may compute, together with a proof that \lstinline{p (x.val)} holds (called \lstinline{x.prop}). For example, a computable function can return \lstinline{x.val + 2}. We can prove specifications about said computable function using \lstinline{x.prop}. This subtype can also be written as the Sigma type \lstinline{Σ (n : ℕ), p(n)}, again of type \lstinline{Type}. Often in this formalization, we ignore the style guidelines to use subtypes, and rather favor sigma types. In the Lean infoview, a Sigma type (unlike a subtype) will show the type of $n$ explicitly, which is important when its type is complicated and requires rewriting to another form.

Now let us discuss the case of inductive definitions. If an inductive definition takes values in \lstinline{Prop}, Lean remembers only that the corresponding proposition has been proved; the particular constructor tree used to prove it is not computational data. If instead the inductive definition takes values in \lstinline{Type}, its constructors form data. Accordingly, we will pass from several of our earlier \lstinline{Prop}-valued relations to closely related \lstinline{Type}-valued versions whenever their construction sequence is needed to define a metric.

Why not simply make all our structures \lstinline{Type}-valued from the beginning? Working with unneeded data is unwieldy. Thus, the Mathlib library is designed to be \lstinline{Prop}-valued, so we cannot benefit from its extensive codebase when working with data (for example, we had to implement our own small data-carrying library of facts about prefixes and suffixes of lists). Similarly, many automated tactics, such as \lstinline{aesop} and \lstinline{grind} do not work on \lstinline{Type}-valued objects. Thus, we work in \lstinline{Prop} as much as possible, and only pass to \lstinline{Type} when the data is absolutely necessary. 

\subsection{Grids}

Currently, grids are \lstinline{Prop}-valued in Lean. Saying that $\langle a, b, c, d\rangle\in G$ means that there exists some method of using the constructors of $G$ to form a grid with frame $\langle a, b, c, d\rangle$. Note that there can be multiple ways of arriving at the same frame : 

\noindent 
\begin{tikzcd}
{} \arrow[rr, "b_1b_2"] \arrow[d, "a1"'] &  & {} \arrow[d, "d_1"] \\
{} \arrow[d, "a2"'] \arrow[rr, "x"']     &  & {} \arrow[d, "d_2"] \\
{} \arrow[rr, "c_1c_2"']                 &  & {}                 
\end{tikzcd}
and
\begin{tikzcd}
{} \arrow[dd, "a_1a_2"'] \arrow[r, "b_1"] & {} \arrow[r, "b_2"] \arrow[dd, "y"] & {} \arrow[dd, "d_1d_2"] \\
                                          &                                     &                         \\
{} \arrow[r, "c_1"']                      & {} \arrow[r, "c_2"']                & {}                     
\end{tikzcd}
are both
\begin{tikzcd}
{} \arrow[dd, "a_1a_2"'] \arrow[rr, "b_1b_2"] &  & {} \arrow[dd, "d_1d_2"] \\
                                              &  &                         \\
{} \arrow[rr, "c_1c_2"']                      &  & {}                     
\end{tikzcd} although they are constructed differently. In order to give a size metric on grids, we will need to track the data of how a given grid was constructed. It is entirely possible to give a definition of size which gives differing values for different construction sequences leading to the same frame. Our metric will provably give the same length to any pair of grids with matching sources and targets, but that is not the case for every metric. (In fact, I am not able to prove such a specification for the eventual size metric I will give on partial grids.)

We thus enhance our definition of grids as follows :

\begin{definition} \leanlink{https://github.com/hannahfechtner/braids_better/blob/098309403cda787eaf1fa954e07ff0ca3b2927ee/BraidProject/GridData/Basic.lean#L10-L20}

We define the set of grids with constructions $G^D$ inductively
\footnote{
    Most precisely, we define the set of constructions $K$ and grids with constructions $G^D$ mutually inductively as follows:
    
    The set of constructions $K$ is defined as the smallest set satisfying
    \begin{itemize}
        \item $\mathsf{trivial \,cell} \in K$
        \item $\mathsf{nontrivial \,cell} \in K$
        \item $\mathsf{horizontal \, (} g_1, g_2 \mathsf{)}\in K$ for $g_1, g_2\in G^D$ with appropriate bounds on the frames of $g_1$ and $g_2$ (that is, the fourth element of $g_1$ equals the first element of $g_2$)
        \item $\mathsf{vertical \, (} g_1, g_2 \mathsf{)}\in K$ for $g_1, g_2\in G^D$ with appropriate bounds on the frames of $g_1$ and $g_2$ (that is, the third element of $g_1$ equals the second element of $g_2$)
    \end{itemize}
     \quad The set of grids with constructions $G^D$ is defined as the smallest set satisfying
        \begin{itemize}
\item The trivial cells : 
\begin{itemize}
    \item $\langle\varepsilon, \varepsilon, \varepsilon, \varepsilon, \mathsf{trivial \, cell}\rangle \in G^D$
    \item For any $i\in \mathbb{N}$, 
$\langle\sigma_i, \varepsilon, \varepsilon, \sigma_i, \mathsf{trivial \, cell}\rangle \in G^D$
\item For any $i\in \mathbb{N}$, 
$\langle\varepsilon, \sigma_i, \sigma_i, \varepsilon, \mathsf{trivial \, cell}\rangle \in G^D$
\end{itemize}

\item The nontrivial cells : 
\begin{itemize}
\item For any $i\in \mathbb{N}$, $\langle\sigma_i,\sigma_i,\varepsilon,\varepsilon, \mathsf{nontrivial \, cell}\rangle\in G^D$ 
\item For $i,j\in \mathbb{N}$ with $|i-j|=1$, $\langle\sigma_j,\sigma_i,\sigma_i\sigma_j,\sigma_j\sigma_i, \mathsf{nontrivial \, cell} \rangle\in G^D$.
\item For $i,j\in\mathbb{N}$ with $|i-j|>1$, $\langle\sigma_j,\sigma_i,\sigma_i,\sigma_j, \mathsf{nontrivial \, cell}\rangle\in G^D$
\end{itemize}
\item Horizontal abutting: Given $a, b, c, d, e, f, g\in BW^{\mathds{1}, +}$ and $\kappa_1, \kappa_2 \in K$, if $g_1=\langle a, b, c,d, \kappa_1 \rangle$ and $g_2 =\langle d,e,f,g,\kappa_2 \rangle$ are elements of  $G^D$, then $\langle a,be,cf,g, \mathsf{horizontal}(g_1, g_2)\rangle\in G^D$.

\item Vertical abutting: Given $a, b, c, d, e, f, g\in BW^{\mathds{1}, +}$ and $\kappa_1, \kappa_2 \in K$, if  $g_1=\langle a, b, c,d, \kappa_1\rangle$ and $\langle e,c,f,g, \kappa_2\rangle$ are elements of $G^D$, then $\langle ae,b,f,dg, \mathsf{vertical}(g_1, g_2)\rangle\in G^D$.

\end{itemize}} :
    \begin{itemize}
\item The trivial cells\footnote{In Lean, each of the six cell cases is listed individually; here we combine them into two groups, trivial and nontrivial, to condense the presentation} : 
\begin{itemize}
    \item $\langle\varepsilon, \varepsilon, \varepsilon, \varepsilon, \mathsf{trivial \, cell}\rangle \in G^D$
    \item For any $i\in \mathbb{N}$, 
$\langle\sigma_i, \varepsilon, \varepsilon, \sigma_i, \mathsf{trivial \, cell}\rangle \in G^D$
\item For any $i\in \mathbb{N}$, 
$\langle\varepsilon, \sigma_i, \sigma_i, \varepsilon, \mathsf{trivial \, cell}\rangle \in G^D$
\end{itemize}

\item The nontrivial cells : 
\begin{itemize}
\item For any $i\in \mathbb{N}$, $\langle\sigma_i,\sigma_i,\varepsilon,\varepsilon, \mathsf{nontrivial \, cell}\rangle\in G^D$ 
\item For $i,j\in \mathbb{N}$ with $|i-j|=1$, $\langle\sigma_j,\sigma_i,\sigma_i\sigma_j,\sigma_j\sigma_i, \mathsf{nontrivial \, cell} \rangle\in G^D$.
\item For $i,j\in\mathbb{N}$ with $|i-j|>1$, $\langle\sigma_j,\sigma_i,\sigma_i,\sigma_j, \mathsf{nontrivial \, cell}\rangle\in G^D$
\end{itemize}
\item Horizontal abutting: If $g_1=\langle a, b, c,d, \kappa_1 \rangle$ and $g_2 =\langle d,e,f,g,\kappa_2 \rangle$ are elements of  $G^D$, then $\big\langle a,be,cf,g, \mathsf{horizontal}(g_1, g_2)\big\rangle\in G^D$.

Visually, 
\begin{tikzcd}
{} \arrow[r, "b"] \arrow[d, "a"'] & {} \arrow[d, "d"] \arrow[r, "e"] & {} \arrow[d, "g"] \\
{} \arrow[r, "c"']                & {} \arrow[r, "f"']               & {}               
\end{tikzcd}\item
Vertical abutting: If $g_1=\langle a, b, c,d, \kappa_1\rangle$ and $\langle e,c,f,g, \kappa_2\rangle$ are elements of $G^D$, then $\big\langle ae,b,f,dg, \mathsf{vertical}(g_1, g_2)\big\rangle\in G^D$.

Visually, 
\begin{tikzcd}
{} \arrow[r, "b"] \arrow[d, "a"']  & {} \arrow[d, "d"] \\
{} \arrow[r, "c"'] \arrow[d, "e"'] & {} \arrow[d, "g"] \\
{} \arrow[r, "f"']                 & {}               
\end{tikzcd}

\end{itemize}
\end{definition}

In Lean, this is handled entirely for us by switching out \lstinline{Prop} (which asserts the proposition that we have some grid) to \lstinline{Type}, which preserves the data of the grid's construction. 

\leanlink{https://github.com/hannahfechtner/braids_better/blob/098309403cda787eaf1fa954e07ff0ca3b2927ee/BraidProject/GridData/Basic.lean#L12}\begin{lstlisting} 
inductive GridData : FreeMonoid ℕ → FreeMonoid ℕ → FreeMonoid ℕ → FreeMonoid ℕ → Type
\end{lstlisting}
It is immediate to delete the construction from a grid-with-construction and get a grid. 

\begin{theorem} \label{thm:grid_of_grid_data}\leanlink{https://github.com/hannahfechtner/braids_better/blob/098309403cda787eaf1fa954e07ff0ca3b2927ee/BraidProject/GridData/Basic.lean#L64-L73}
    If $\langle a, b, c, d, \kappa\rangle\in G^D$, then $\langle a, b, c, d\rangle\in G$.
\end{theorem}

\begin{lstlisting}
def GridData.to_grid (h : GridData a b c d) : grid a b c d
\end{lstlisting}

Casting grids-with-constructions down to regular grids helps us convert many theorems (here, we show Theorem \ref{thm:grid_equiv} and the Remark below Definition \ref{def:diag_length}) as follows:

\leanlink{https://github.com/hannahfechtner/braids_better/blob/098309403cda787eaf1fa954e07ff0ca3b2927ee/BraidProject/GridData/Basic.lean#L98-L101}
\begin{lstlisting}
theorem braid_eq (h : GridData a b c d) : BraidMonoidInf.mk (a * c) = BraidMonoidInf.mk (b * d) := Grid.braid_eq_of_grid (to_grid h) 
\end{lstlisting}

\leanlink{https://github.com/hannahfechtner/braids_better/blob/098309403cda787eaf1fa954e07ff0ca3b2927ee/BraidProject/GridData/Basic.lean#L109-L110}
\begin{lstlisting}
theorem diag_length_eq (h : GridData a b c d) : a.length + c.length = b.length + d.length := Grid.diag_length_eq (to_grid h)
\end{lstlisting}

Conversely, given a grid, we know that there must be at least one construction sequence giving that frame :
\begin{theorem} \label{thm:grid_data_of_grid}\leanlink{https://github.com/hannahfechtner/braids_better/blob/098309403cda787eaf1fa954e07ff0ca3b2927ee/BraidProject/GridData/Basic.lean#L50-L62}
    If $\langle a, b, c, d\rangle\in G$, then there exists some $\kappa$ such that $\langle a, b, c, d, \kappa \rangle\in G^D$.
\end{theorem}
Note that the existential quantifier on $\kappa$ asserts that the set of construction sequences generating such a frame is nonempty. In Lean, \lstinline{Nonempty} has the same meaning : 

\leanlink{https://github.com/hannahfechtner/braids_better/blob/098309403cda787eaf1fa954e07ff0ca3b2927ee/BraidProject/GridData/Basic.lean#L50}\begin{lstlisting}
def GridData.of_grid (h : grid a b c d) : Nonempty (GridData a b c d)
\end{lstlisting}

Thus, any theorems about grids can be translated to grids-with-constructions.

We began with \lstinline{Prop}-valued grids because of the Mathlib support, and only move to the  \lstinline{Type}-valued grids (grids-with-constructions) when necessary. That occurs when we define a metric on grids-with-constructions. The metric is intended to represent the number of nontrivial cells occurring in the grid. 

\begin{definition} \leanlink{https://github.com/hannahfechtner/braids_better/blob/098309403cda787eaf1fa954e07ff0ca3b2927ee/BraidProject/GridData/Length.lean#L8-L18}
    Define a metric on grids-with-constructions $\ell_G : G^D \to \mathbb{N}$ recursively as follows : 
    \begin{itemize}
        \item $\ell_G\big(\langle a, b, c, d, \mathsf{trivial \, cell}\rangle\big) = 0$
        \item $\ell_G\big(\langle a, b, c, d, \mathsf{nontrivial \, cell}\rangle\big) = 1$
        \item $\ell_G\big(\langle a, b, c, d, \mathsf{horizontal}(g_1, g_2)\rangle\big) = \ell_G(g_1)+\ell_G(g_2)$
        \item $\ell_G\big(\langle a, b, c, d, \mathsf{vertical}(g_1, g_2)\rangle\big) = \ell_G(g_1)+\ell_G(g_2)$
    \end{itemize}
\end{definition}

We expand a number of theorems to discuss the grid size. For example, when splitting a grid, we show that the sum of the sizes of  two new, smaller grids equals the size of the original grid. We run through the full proof of splittability given in Lemma \ref{lem:grid_splittability} again, slightly augmented with  facts about size, to obtain the following lemma : 
\begin{lemma}(Update to Lemma \ref{lem:grid_splittability}) \leanlink{https://github.com/hannahfechtner/braids_better/blob/098309403cda787eaf1fa954e07ff0ca3b2927ee/BraidProject/GridData/Length.lean#L156-L257}
Given $G_0=\langle a_1a_2,b,c,d, \kappa_0\rangle \in G$, there exist grids $G_1=\langle a_1,b,m,d_1, \kappa_1\rangle$ and $G_2=\langle a_2,m,c,d_2, \kappa_2\rangle$ with $d=d_1d_2$ \textcolor{magenta}{and $\ell_G(G_0)=\ell_G(G_1)+\ell_G(G_2)$}. 

\noindent Similarly, given $G_3=\langle a,b_1b_2,c,d, \kappa_3\rangle \in G$, there exist grids $G_4=\langle a,b_1,c_1,m,\kappa_4 \rangle$ and $\langle m,b_2,c_2,d, \kappa_5\rangle$ with $c=c_1c_2$ \textcolor{magenta}{and $\ell_G(G_3)=\ell_G(G_4)+\ell_G(G_5)$}. 

\end{lemma}

Our main theorem is that the spine of a grid suffices to determine its length. We begin with a lemma discussing certain small cases : 

\begin{lemma}\label{lem:determinative_spine_length} [Determinative Spines -- Length] \mbox{} \leanlink{https://github.com/hannahfechtner/braids_better/blob/098309403cda787eaf1fa954e07ff0ca3b2927ee/BraidProject/GridData/Length.lean#L28-L153}
For $a,b,c,d \in BW^+$ and $i,j \in \mathbb{N}$,
\begin{itemize}
    \item If $\langle \varepsilon, b, c, d\rangle\in G$, then $\ell_G\big(\langle \varepsilon, b, c, d\rangle\big)=0$.
    \item If $\langle a, \varepsilon, c, d\rangle\in G$, then $\ell_G\big(\langle a, \varepsilon, c, d\rangle\big)=0$.
    \item If $\langle a, a, c, d\rangle\in G$, then $\ell_G\big(\langle a, a, c, d\rangle\big) = |a|$. Here, $|a|$ is the length of the word $a$.
    \item If $\langle \sigma_i, \sigma_j, c, d\rangle\in G$ with $|i-j|=1$, then $\ell_G\big(\langle \sigma_i, \sigma_j, c, d\rangle\big)=1$.
    \item If $\langle \sigma_i, \sigma_j, c, d\rangle\in G$ with $|i-j|>1$, then $\ell_G\big(\langle \sigma_i, \sigma_j, c, d\rangle\big)=1$.
\end{itemize}
\end{lemma}

\noindent Now, thanks to uniqueness of grid spines (Theorem \ref{thm:grid_uniqueness}), we may show the full theorem.
\begin{theorem}\label{thm:same_type_same_length_grid} \leanlink{https://github.com/hannahfechtner/braids_better/blob/098309403cda787eaf1fa954e07ff0ca3b2927ee/BraidProject/GridData/Length.lean#L261-L293}
    Given $g_1=\langle a,b,c,d, \kappa_1\rangle \in G^D$ and $g_2=\langle a,b,e, f, \kappa_2\rangle \in G^D$, $\ell_G(g_1)=\ell_G(g_2)$.
\end{theorem}
\begin{proof} (Sketch)
    We proceed by induction on $g_1$. In the base cases where $g_1$ is a single cell, Lemma \ref{lem:determinative_spine_length} allows us to conclude that $g_2$ has the same size as $g_1$. We present the horizontal inductive case; the vertical inductive case is analogous. 

    Assume that $g_1=\big\langle a,b,c,d, \mathsf{horizontal}(g_3, g_4)\big\rangle$, with $b=b_1b_2$, $c =c_1c_2$, and 
    $g_3=\langle a, b_1, c_1,m, \kappa_1 \rangle$ and $g_4=\langle m,b_2,c_2,d,\kappa_2\rangle$ to which the inductive hypothesis applies. Thus, $\ell_G(g_1)=\ell_G(g_3)+\ell_G(g_4)$.
    Then we split $g_2$ to obtain $g_5=\langle a, b_1, c_5,d_5, \kappa_1 \rangle$ and $g_6=\langle d_5,b_2,c_6,f,\kappa_2\rangle$ with $\ell_G(g_2)=\ell_G(g_5)+\ell_G(g_6)$. Uniqueness of (data-carrying) grid spines (Theorem \ref{thm:grid_uniqueness} combined with Theorem \ref{thm:grid_of_grid_data}) applied to $g_3$ and $g_5$ gives us $c_5=c_1$ and $d_5=m$, and so we may apply our induction hypothesis to those two to obtain $\ell_G(g_3)=\ell_G(g_5)$. Since $d_5=m$, $g_6=\langle m,b_2,c_6,f,\kappa_2\rangle$, and another application of Theorem \ref{thm:grid_uniqueness} to $g_4$ and $g_6$ gives us $c_6=c_2$. Finally, another application of our inductive hypothesis gives us $\ell_G(g_4)=\ell_G(g_6)$, which is sufficient to conclude.  
\end{proof}

We use the fact that the spine of a grid suffices to determine its size to inspire the following metric $\ell_B$ on pairs of braid monoid words.

\begin{definition}\label{def:pair_length}\leanlink{https://github.com/hannahfechtner/braids_better/blob/098309403cda787eaf1fa954e07ff0ca3b2927ee/BraidProject/GridData/Length.lean#L22-L24}
    Given $a,b\in BW^+$, we know by
    Theorem \ref{thm:grid_existence} and Theorem \ref{thm:grid_data_of_grid}
    that there exist $c,d\in BW^+$ and a construction
    $\kappa$ such that
    $g=\langle a,b,c,d,\kappa\rangle\in G^D$.
    We define
    \[
        \ell_B(a,b):=\ell_G(g).
    \]
    This is well-defined independently of the choice of
    $c,d$, and $\kappa$ by
    Theorem~\ref{thm:same_type_same_length_grid}.
\end{definition}

In Lean, we take a slightly different approach to well-definedness.
We have noncomputable \lstinline{GridData.existence} (via \lstinline{Prop}-valued grid existence and Theorem \ref{thm:grid_data_of_grid}). We then use the axiom of choice to select such a grid to serve as $g$. We define \lstinline{ab_len}
to be the length of this selected grid. We do not use Theorem \ref{thm:same_type_same_length_grid} and do not need to prove that our metric is independent of $c, d,$ and $\kappa$. Eventually, when we use this metric, we will show that $\ell_B$ is bounded in the desired manner for any $c,d,\kappa$, which will suffice in place of well-definedness. Furthermore, there is no issue with $\ell_B$ being noncomputable -- we can show termination of a computable algorithm via noncomputable termination metrics.

\subsection{Partial Grids}

We similarly update the definition of partial grids to track the data of their construction. In fact, throughout our Lean formalization we have used \lstinline{Type}-valued partial grids. Unlike grids, where there were significant theorems not mentioning size (such as stability, Theorem \ref{thm:stability}), almost every theorem relating to partial grids discusses size. So we had little use for \lstinline{Prop}--valued partial grids, and hence did not implement them. Therefore, the reader should view the following definition of $G_P$ as the true definition, and look back upon Definition \ref{def:partial_grid} as a simplified version. Thus, we do not need the clunky notation $G_P^D$.

\begin{definition} (Update to Definition \ref{def:partial_grid})
     We define the set of partial grids $G_P$ as follows\footnote{Again, this would be most rigorously defined mutually inductively with a set of constructions $K_p$.} :
\begin {itemize}
\item \textbf{a trivial cell} For any $a,b,c,d\in BW^{\mathds{1}, +}$, if $(\overline{a}\,b, \;c\,\overline{d})\in R_{trivial}$, we have $\langle \overline{a},\; b,\; c, \; \varepsilon,\; \overline{d}\textcolor{magenta}{, \mathsf{trivial \,cell}}\rangle \in G_P$

\item \textbf{a nontrivial cell} For any $a,b,c,d\in BW^{\mathds{1}, +}$, if $(\overline{a}\,b, \;c\,\overline{d})\in R_{nontrivial}$, we have \newline $\langle \overline{a},\; b,\; c, \; \varepsilon,\; \overline{d}\textcolor{magenta}{, \mathsf{nontrivial \,cell}}\rangle \in G_P$

\item \textbf{a spine} For any $a,b\in \left\{BW^{\mathds{1},+}\,\backslash \, \{\varepsilon\}\right\}$, $\langle \overline{a},\; b,\; \varepsilon, \; \overline{a}\, b,\; \varepsilon\textcolor{magenta}{, \mathsf{spine}}\rangle \in G_P$ 

\item \textbf{horizontal concatenation to a partial grid with empty middle frontier}  :

Given partial grids $g_1=\langle a_1, b_1, c_1, \varepsilon, d_1\textcolor{magenta}{, \kappa_1}\rangle \in G_P$ and $g_2=\langle d_1, b_2, c_2, m_2, d_2\textcolor{magenta}{, \kappa_2}\rangle \in G_P$, we have $\big\langle a_1, b_1\cdot b_2, c_1\cdot c_2, m_2, d_2\textcolor{magenta}{, \mathsf{horizontal_{\varepsilon}}(g_1, g_2)}\big\rangle \in G_P$.

\item \textbf{horizontal concatenation to a partial grid with non-empty middle frontier} ($m_1$ is non-empty) : 

Given $g_1=\langle a_1, b_1, c_1, m_1, d_1\textcolor{magenta}{, \kappa_1}\rangle \in G_P$ and $g_2=\langle d_1, b_2, c_2, m_2, d_2\textcolor{magenta}{, \kappa_2}\rangle \in G_P$, we have $\big\langle a_1, b_1\cdot b_2, c_1, m_1 \cdot c_2 \cdot m_2, d_2\textcolor{magenta}{, \mathsf{horizontal}(g_1, g_2)}\big\rangle \in G_P$.

\item \textbf{vertical concatenation to a partial grid with empty middle frontier} :

Given $g_1=\langle a_1, b_1, c_1, \varepsilon, d_1\textcolor{magenta}{, \kappa_1}\rangle \in G_P$ and $\langle a_2, c_1, c_2, m_2, d_2\textcolor{magenta}{, \kappa_2}\rangle \in G_P$, we have $\big\langle a_2\cdot a_1, b_1, c_2, m_2,d_2\cdot d_1\textcolor{magenta}{, \mathsf{vertical_{\varepsilon}(g_1, g_2)}}\big\rangle \in G_P$.

\item \textbf{vertical concatenation to a partial grid with non-empty middle frontier} 

Given $\langle a_1, b_1, c_1, m_1, d_1\textcolor{magenta}{, \kappa _1}\rangle \in G_P$ and $\langle a_2, c_1, c_2, m_2, d_2\textcolor{magenta}{, \kappa_2}\rangle \in G_P$, we have $\big\langle a_2\cdot a_1, b_1, c_2, m_2 \cdot d_2 \cdot m_1, d_1\textcolor{magenta}{, \mathsf{vertical}(g_1, g_2)}\big\rangle \in G_P$.

\end{itemize}
\end{definition}

\begin{definition} \leanlink{https://github.com/hannahfechtner/braids_better/blob/098309403cda787eaf1fa954e07ff0ca3b2927ee/BraidProject/PartialGrid/Basic.lean#L34-L48}
    Define a metric on partial grids $\ell_P$ by the following : 
    \begin{itemize}
        \item $\ell_P\left(\langle a, b, c, m, d, \mathsf{trivial \, cell}\rangle\right) = 0$
        \item $\ell_P\left(\langle a, b, c, m, d, \mathsf{nontrivial \, cell}\rangle\right) = 1$
        \item $\ell_P\left(\langle a, b, c, m, d, \mathsf{spine}\rangle\right) = 0$
        \item $\ell_P\left(\langle a, b, c, m, d, \mathsf{horizontal_\varepsilon}(p_1, p_2)\rangle\right) = \ell_P(p_1)+\ell_P(p_2)$
        \item $\ell_P\left(\langle a, b, c, m, d, \mathsf{horizontal}(p_1, p_2)\rangle\right) = \ell_P(p_1)+\ell_P(p_2)$
        \item $\ell_P\left(\langle a, b, c, m, d, \mathsf{vertical_\varepsilon}(p_1, p_2)\rangle\right) = \ell_P(p_1)+\ell_P(p_2)$
        \item $\ell_P\left(\langle a, b, c, m, d, \mathsf{vertical}(p_1, p_2)\rangle\right) = \ell_P(p_1)+\ell_P(p_2)$
    \end{itemize}
\end{definition}

\noindent It is possible to extend a partial grid's spine without affecting its size.
\begin{lemma}\label{lem:extend_frontier_w_length} (Update to Theorem \ref{lem:extend_frontier})
    Consider a partial grid $p = \langle a, b, c, m, d, \kappa\rangle$. For any non-empty \textit{negative} augmented grid word $a'$ there exists a partial grid $p_1 =\langle a'\cdot a, b, \varepsilon, a'\cdot c\cdot m, d, \kappa_1\rangle$ \textcolor{magenta}{such that $\ell_{P}(p)=\ell_{P}(p_1)$}. For any non-empty \textit{positive} augmented grid word $b'$, there exists a partial grid $p_2 =\langle a, b\cdot b', c, m \cdot d \cdot b', \varepsilon , \kappa_2\rangle$ \textcolor{magenta}{such that $\ell_{P}(p)=\ell_{P}(p_2)$}.
\end{lemma}

Splittability also behaves as we might expect in regard to the size of the resulting partial grids.  

\begin{theorem}\label{thm:split_with_length} (Update to Theorem \ref{thm:splittability})
    Given $p = \langle a, b_1b_2, c, m, d\rangle\in G_P$ with nonempty $b_1, b_2$, we have either
    \begin{itemize}
        \item $p$ splits into two smaller partial grids $p_1 = \langle a, b_1, c_1, m_1, d_1\rangle\in G_P$ and $p_2 = \langle d_1, b_2, c_2, m_2, d\rangle\in G_P$, with $c\cdot m = c_1 \cdot m_1 \cdot c_2 \cdot m_2$\textcolor{magenta}{, and $\ell_P(p)=\ell_P(p_1)+\ell_P(p_2)$}.

\item $p$ splits into a smaller partial grid $p_1 = \langle a, b_1, c, m_1, \varepsilon\rangle\in G_P$ and a tail \begin{tikzcd}
{} \arrow[r, "b_2"] & {}
\end{tikzcd}, with $m =m_1\cdot b_2$ \textcolor{magenta}{and $\ell_P(p)=\ell_P(p_1)$}.
\end{itemize}
\end{theorem}

We update our \lstinline{add_cell} function to add in a fact about the size of the resulting partial grid : 

\begin{definition} \label{def:new_add_cell}
    \lstinline{add_cell} is a function taking in some $(x,y)\in R_{\mathsf{grid}}$ and a partial grid $p_1 =\langle a, b, c, m_1\cdot x\cdot m_2,d\textcolor{magenta}{,\kappa_1}\rangle$, and returning a partial grid $p_2=\langle a, b, c', m', d'\textcolor{magenta}{, \kappa_2}\rangle$ where $c'\cdot m'\cdot d'=c\cdot m_1\cdot y\cdot m_2\cdot d$, $c$ is a prefix of $c'$, $d$ is a suffix of $d'$\textcolor{magenta}{, and the following property holds : if $(x,y)\in R_{\mathsf{triv}}$, $\ell_P(p_2)=\ell_P(p_1)$; otherwise, we have $(x,y)\in R_{\mathsf{nontriv}}$, and $\ell_P(p_2)=\ell_P(p_1)+1$}. 
\end{definition}

Note that we could separate out the final property about length, and prove it separately as a theorem about the output of \lstinline{add_cell} (in contrast to, say, the prefix and suffix properties, which are required for the definition of \lstinline{add_cell}). However, like the definition of \lstinline{add_cell}, such a proof would have 17 cases\footnote{The cases for horizontal and vertical appending to a partial grid with non-empty frontier have 5 sub-cases each. The case for a single cell comprises 4 subcases. The other 3 inductive cases (horizontal and vertical appending to a partial grid with empty middle frontier and the spinal case) do not split further.}, so it is easier to fold the length property into the main definition.

Our main goal will be to show that given a grid and a partial grid sharing essentially the same spine (aside from any empty arrows), the size of the partial grid is at most that of the grid. Our proof requires us to strengthen the inductive statement with an auxiliary fact about prefixes and suffixes of the grid's target (part (1) below).

\begin{theorem} \label{thm:nested_partial_grids_and_grids} \leanlink{https://github.com/hannahfechtner/braids_better/blob/098309403cda787eaf1fa954e07ff0ca3b2927ee/BraidProject/PartialGrid/NestedFrame.lean#L11-L233}
Let $g = \langle a, b, c, d\rangle\in G$ and $p = \langle a_0, b_0, c_0, m_0, d_0\rangle \in G_P$. Recall that $s^{\mathds{1},+} : BW^{\mathds{1}, +}\to BW^+$ removes the signs on the letters, and deletes all instances of $\mathds{1}$ from a word. It is generally used to convert (positive) partial grid edges to grid edges. If $a = s^{\mathds{1},+}(\overline {a_0})$ and $b=s^{\mathds{1},+}({b_0})$ (that is, the spines agree, ignoring empty arrows), then
\begin{enumerate}
    \item $s^{\mathds{1},+}({c_0})$ is a prefix of $c$ and $s^{\mathds{1},+}(\overline {d_0})$ is a prefix of $d$
    \item $\ell_P(p)\leq\ell_G(g)$.
\end{enumerate}

Here is the image to keep in mind : 
\begin{center}
\begin{tikzcd}
{} \arrow[dd, "a"'] 
   \arrow[rr, "b"] 
   \arrow[rr, "b_0"', draw=none,
     preaction={draw=magenta, draw opacity=.35, line width=4pt, -},
     text=magenta]
   \arrow[dd, "a_0", draw=none,
     preaction={draw=magenta, draw opacity=.35, line width=4pt, -},
     text=magenta]
&
&
{} \arrow[dd, "d"] 
   \arrow[d, "d_0"', draw=none,
     preaction={draw=magenta, draw opacity=.35, line width=4pt, -},
     text=magenta]
\\
&
&
{} 
\\
{} \arrow[rr, "c"'] 
   \arrow[r, "c_0", draw=none,
     preaction={draw=magenta, draw opacity=.35, line width=4pt, -},
     text=magenta]
&
{} \arrow[ru, "m_0",draw=none,
     preaction={draw=magenta, draw opacity=.35, line width=4pt, shorten <=-4pt, shorten >=-4pt},
     text=magenta, dash, squiggly]
&
{}
\end{tikzcd}
\end{center}

\end{theorem}

We note that we could prove $(1)$ independently, and then use $(1)$ to later prove $(2)$. (We cannot proceed conversely; $(2)$ cannot be proven without $(1)$.) Due to the large number of cases in the proof, and the work to organize the hypotheses for each case, it is preferable to combine $(1)$ and $(2)$ as we have done, proving them simultaneously.
\begin{proof} (Sketch)
    We proceed by induction on $g$, and present one illustrative (sub-sub-)case. Consider the case where $g = \big\langle a, b, c, d, \mathsf{horizontal}(g_1, g_2)\big\rangle$ for grids $g_1 = \langle a, b_1, c_1, d_1, \kappa_1\rangle$ and $g_2 = \langle d_1, b_2, c_2, d, \kappa_2\rangle $ to which the inductive hypothesis applies, with $b=b_1b_2$ and $c=c_1c_2$. Of course, $\ell_G(g)=\ell_G(g_1)+\ell_G(g_2)$.
    As $s^{\mathds{1},+}( {b_0}) = b=b_1b_2$, we know there exist some $b_3, b_4\in BW^{1,+}$  such that $b_0=b_3b_4$ and $b_1=s^{\mathds{1},+}( {b_3})$ and $b_2=s^{\mathds{1},+}( {b_4})$. We may thus split  $p=\langle a_0, b_3b_4, c_0, m_0, d_0\rangle$ vertically. Consider the case in which $p$ splits in to two sub-grids $p_1=\langle a_0, b_3, c_3, m_3, d_3 \rangle$  and $p_2 = \langle d_3, b_4, c_4, m_4, d_0\rangle$ with $c_3m_3c_4m_4 = c_0m_0$ and $\ell_P(p)=\ell_P(p_1)+\ell_P(p_2)$ (the other case is simpler). Further take the case in which $m_3$ is non-empty, and thus $c_0=c_3$ and $m_0=m_3c_4m_4$. We apply our inductive hypothesis to $g_1$ and $p_1$ to conclude that $s^{\mathds{1}, +}(c_3)$ is a prefix of $c_1$ (and thus $s^{\mathds{1}, +}(c_0)$ is a prefix of $c=c_1c_2$ as desired), $s^{\mathds{1}, +}\left(\overline{d_3}\right)$ is a prefix of $d_1$, and $\ell_P(p_1)\leq \ell_G(g_1)$. Thus, there is some negative augmented braid word $d_5$ such that $d_1=s^{\mathds{1}, +}\left(\overline{d_3}\right)s^{\mathds{1}, +}\left(\overline{d_5}\right)$. We may extend the left frontier of $p_2$ to obtain a partial grid $p_3=\langle d_5 \cdot d_3, b_4, \varepsilon, d_5 \cdot c_4 \cdot m_4, d_0\rangle$ with $\ell_P(p_3)=\ell_P(p_2)$. Then we know $d_1 = s^{\mathds{1}, +}\left(\overline{d_5\cdot d_3}\right)$, and so we may apply our inductive hypothesis to $g_2$ and $p_3$ to conclude that $s^{\mathds{1},+}(\overline {d_0})$ is a prefix of $d$ as desired, and that $\ell_P(p_3)\leq\ell_G(g_2)$. Thus, our final condition, that $\ell_P(p)\leq \ell_G(g)$, is satisfied.
\end{proof}

\subsection{Semi-Thue Derivations}

We return to our two definitions of Semi-Thue systems, defining data-carrying versions so that we may define a metric on them. In Lean, this once again corresponds to switching the output type from \lstinline{Prop} to \lstinline{Type}.
\begin{definition} \leanlink{https://github.com/hannahfechtner/braids_better/blob/098309403cda787eaf1fa954e07ff0ca3b2927ee/BraidProject/DataCarrying/SemiThue.lean#L8-L11}

Consider a set $A$ and a binary relation $R$ on $A$-strings. Then $\xrightarrow[R]{*}$ is the smallest relation satisfying : 
\begin{itemize}

\item reflexivity : For all $a\in A^*$, $\langle a \xrightarrow[R]{*} a, \mathsf{refl}\rangle$

\item closure : For all $a, b\in A^*$, $a\xrightarrow[R]{} b$ implies  $\langle a\xrightarrow[R]{*} b, \mathsf{closure} (a, b)\rangle$

\item transitivity : For all $a, b, c \in A^*$, if  $s_0=\langle a\xrightarrow[R]{*} b, \kappa_0\rangle$ and  $s_1=\langle b\xrightarrow[R]{*} c, \kappa_1\rangle$, then  $\langle a \xrightarrow[R]{*} c, \mathsf{trans} (s_0, s_1)\rangle$
\end{itemize}
\end{definition}

\begin{definition} \leanlink{https://github.com/hannahfechtner/braids_better/blob/098309403cda787eaf1fa954e07ff0ca3b2927ee/BraidProject/DataCarrying/SemiThue.lean#L104-L108}
    We define a metric on Semi-Thue systems $\ell_{S_R}$, parametrized by the relation $R$ and a size metric $\ell_R : R\to \mathbb{N}$ on pairs of $A$-strings satisfying $R$. To avoid triple subscripts, we do not add the particular choice of $\ell_R$ to our notation, allowing it to be inferred from context. In the below definition, for readability we use an $\_$ as a placeholder, as no information from the first element of the pair is needed. 
    \begin{itemize}

    \item $\ell_S \big(\langle\, \text{\textunderscore}, \mathsf{refl})\rangle\big) = 0 $

\item $\ell_{S_R}\big(\langle\, \_, \mathsf{closure} (a, b)\rangle \big) = \ell_R(a, b)$

\item $\ell_{S_R}\big(\langle \,\_, \mathsf{trans} (s_0, s_1)\rangle \big) = \ell_{S_R}(s_0)+\ell_{S_R}(s_1)$
    \end{itemize}
\end{definition}

\begin{definition} \leanlink{https://github.com/hannahfechtner/braids_better/blob/098309403cda787eaf1fa954e07ff0ca3b2927ee/BraidProject/DataCarrying/SemiThue.lean#L20-L23}
    Given a set $A$ and a binary relation $R$ over $A$-strings, let $\underset{R}{\rightarrowtriangle}$ be the smallest relation satisfying the below two properties :
\begin{itemize}
\item reflexivity : For all $a\in A^*$, $\langle a \underset{R}{\rightarrowtriangle} a, \mathsf{refl}\rangle $

\item one-step transitivity :
For all $a, b_1,x,b_2,y \in A^*$, given $d_0= \langle a\underset{R}{\rightarrowtriangle} b_1xb_2, \kappa_0\rangle $ and  $xRy$, then  $\langle a\underset{R}{\rightarrowtriangle} b_1yb_2, \mathsf{step\, transitivity} (d_0, x, y)\rangle$
\end{itemize}
\end{definition} 

\begin{definition} \leanlink{https://github.com/hannahfechtner/braids_better/blob/098309403cda787eaf1fa954e07ff0ca3b2927ee/BraidProject/DataCarrying/SemiThue.lean#L220-L222} We define a metric on derivation-style Semi-Thue systems $\ell_{D_R}$, again parametrized by the relation $R$ and a size metric $\ell_R$. 
    \begin{itemize}
\item $\ell_{D_R}\big(\langle \_, \mathsf{refl}\rangle\big) = 0 $

\item $\ell_{D_{R}}\big(\langle\_, \mathsf{step\, transitivity} (d_0, x, y)\rangle\big) = \ell_{D_R}(d_0)+\ell_R(x, y)$
\end{itemize}
\end{definition}

\begin{theorem} \leanlink{https://github.com/hannahfechtner/braids_better/blob/098309403cda787eaf1fa954e07ff0ca3b2927ee/BraidProject/DataCarrying/SemiThue.lean#L245-L267}
    Given $s_0 = \langle a \xrightarrow[R]{*} b, \kappa_0\rangle$, there exists a derivation $d_0=\langle a \underset{R}{\rightarrowtriangle} b, \kappa_1\rangle$ with $\ell_{S_R}(s_0)=\ell_{D_R}(d_0)$. The converse of this statement also holds. 
\end{theorem}

\noindent With the basics set up, we now move on to discuss our specific relations. 

\begin{definition} \leanlink{https://github.com/hannahfechtner/braids_better/blob/098309403cda787eaf1fa954e07ff0ca3b2927ee/BraidProject/Relations.lean#L14}
    We define $\ell_{R_{\mathsf{rev}}}$ simply : For all $(x,y)\in R_{\mathsf{rev}}$, $\ell_{R_{\mathsf{rev}}}(x,y)=1$.
    
    $\ell_{R_{\mathsf{grid}}}$ is defined as follows : \leanlink{https://github.com/hannahfechtner/braids_better/blob/098309403cda787eaf1fa954e07ff0ca3b2927ee/BraidProject/Relations.lean#L54-L60}
    \begin{itemize}
    \item For all $(x,y)\in R_{\mathsf{triv}}$, $\ell_{R_{\mathsf{grid}}}(x,y)=0$
    \item For all $(x,y)\in R_{\mathsf{nontriv}}$, $\ell_{R_{\mathsf{grid}}}(x,y)=1$
\end{itemize}
\end{definition}

Next, we update our results from the previous chapter with facts about our metrics. The proofs proceed in an entirely analogous manner. Here are the main conclusions (of course, each of these requires a number of lemmas to be updated as well; we omit those statements and proofs) :
\begin{theorem} (Update to the backwards direction of Theorem \ref{thm:big_one})\leanlink{https://github.com/hannahfechtner/braids_better/blob/098309403cda787eaf1fa954e07ff0ca3b2927ee/BraidProject/SemiThueReversing.lean#L186-L204}
    Given a partial grid $p=\langle a, b, c, m, d, \kappa\rangle$, there exists a derivation $s = \langle \overline{a}b \xrightarrow[R_{\mathsf{rev}}]{*}cm\overline{d}, \kappa'\rangle$ \textcolor{magenta}{with $\ell_P(p)=\ell_{S_{R_{\mathsf{rev}}}}(s)$}.
\end{theorem}

\begin{theorem} (Update to Theorem \ref{thm:new_rw})
    Let $a,b,c,d\in BW^{\mathds{1}, +}$. Given $s_0 = \langle \overline{a}b \xrightarrow[R_{\mathsf{rev}}]{*} c\overline{d}, \kappa_0\rangle$, there exist $e,f \in BW^{\mathds{1}, +}$ and a derivation $s_1=\langle i^{\mathds{1}}(\overline{a}b) \xrightarrow[R_{\mathsf{grid}}]{*} e\overline{f}, \kappa_1\rangle$ with $s^{\mathds{1}}(e)=c, \, s^{\mathds{1}}(f)=d$ \textcolor{magenta}{and $\ell_{S_{R_{\mathsf{rev}}}}(s_0)=\ell_{S_{R_{\mathsf{grid}}}}(d_1)$}. 
\end{theorem}

\begin{theorem} \label{thm:new_partial_grid_of_rw} (Update to Theorem \ref{thm:partial_grid_of_rw})
    For nonempty $a,b\in BW^{\mathds{1}, +}$, given $s_0=\langle\overline{a}b \; \xrightarrow[R_{\mathrm{grid}}]{*} c, \kappa_0\rangle$, then there is a partial grid $p$ with spine $(a, b)$ and frontier $c$ \textcolor{magenta}{such that $\ell_{S_{R_{\mathsf{rev}}}}(s_0)=\ell_P(p)$}.
\end{theorem}

\subsection{Alternative Choices}\label{sec:alt_choices}

We have made a number of seemingly arbitrary decisions, mainly in the precise definition of partial grids, as well as that of their size. Here we present their shortfalls, and discuss alternative strategies with their own benefits and shortfalls. As the definitions and theorems here do not form part of the final project, the Lean formalizations linked may be a bit rougher, and we take a few more liberties in our presentation, deviating from the exact Lean formalization for the sake of exposition. 

Let us begin with the shortfalls of our current definitions. 
There is one tantalizing theorem we have not been able to prove using structural induction. It would be nice to know that the frame of a partial grid determines its length :

\begin{conjecture}\label{conj:same_type_same_length}
    Given partial grids $p_1=\langle a, b, c, m, d, \kappa_1\rangle$ and $p_2=\langle a, b, c, m, d, \kappa_2\rangle$, we have $\ell_{P}(p_1)=\ell_{P}(p_2)$
\end{conjecture}

The obvious strategy of inducting on one of either $p_1$ or $p_2$ fails. Let us play through one case of this proof attempt to see where the failure lies. Say we induct on $p_1$, and are in the case of horizontal appending (to a partial grid with non-empty middle frontier). 

\begin{center}
\begin{tikzcd}
{} \arrow[dd, "a"'] \arrow[rr, "b"] &                     & {} \arrow[d, "d"] \\
                                    &   p_1                  & {}                \\
{} \arrow[r, "c"']                  & {} \arrow[ru, "m"', dash, squiggly] &                  
\end{tikzcd}
is in fact
\begin{tikzcd}
{} \arrow[ddd, "a"'] \arrow[rr, "b_3"] 
&                     
& {} \arrow[rr, "b_4"] \arrow[dd, "d_3"] 
&                     
& {} \arrow[d, "d"] \\
& {p_3} 
&                                    
& {p_4} 
& {}                \\
&                     
& {} \arrow[r, "c_4"']                 
& {} \arrow[ru, "m_4"', dash, squiggly] 
&                   \\
{} \arrow[r, "c"']                   
& {} \arrow[ru, "m_3"', dash, squiggly] 
&                                    
&                     
&                  
\end{tikzcd}
\end{center}

That is, we know there are nonempty $b_3, b_4$ with $b=b_3b_4$ and nonempty $m_3$ and some $c_4, m_4$ with $m=m_3c_4m_4$, partial grids $p_3=\langle a, b_3, c, m_3, d_3, \kappa_3\rangle$ and $p_4=\langle d_3, b_4, c_4, m_4, d, \kappa_4\rangle$ to which the inductive hypothesis may be applied, with $p_1=\langle a, b_3\cdot b_4, c, m_3 \cdot c_4 \cdot m_4, d, \mathsf{horizontal}(p_3, p_4)\rangle$. We would next split $p_2$ horizontally using Theorem \ref{thm:split_with_length}, in hopes of applying our inductive hypothesis for $p_3$ and $p_4$ to the resultant grids. Consider the case in which $p_2$ splits into two partial grids $p_5,p_6$ as shown: 

\begin{center}

\begin{tikzcd}
{} \arrow[ddd, "a"'] \arrow[rr, "b_3"] 
&                     
& {} \arrow[rr, "b_4"] \arrow[dd, "d_5"] 
&                     
& {} \arrow[d, "d"] \\
& {p_5} 
&                                    
& {p_6} 
& {}                \\
&                     
& {} \arrow[r, "c_6"']                 
& {} \arrow[ru, "m_6"', dash, squiggly] 
&                   \\
{} \arrow[r, "c_5"']                   
& {} \arrow[ru, "m_5"', dash, squiggly] 
&                                    
&                     
&                  
\end{tikzcd}
\end{center}

Here, we have partial grids $p_5 = \langle a, b_3, c_5, m_5, d_5\rangle$ and $p_6 = \langle d_5, b_4, c_6, m_6, d\rangle$, with $c\cdot m = c_5 \cdot m_5 \cdot c_6 \cdot m_6$. We may further select the case where $m_5$ is nonempty, and so by properties of middle frontiers (Lemma \ref{lem:frontierup}), $c=c_5$ and $m=m_5c_6m_6$. We wish to apply the inductive hypothesis to $p_5$ and $p_3$. However, we do not know that their frames match. Their spines do match, as do their lower frontiers, but not their middle frontiers ($m_3$ and $m_5$) nor their right frontiers ($d_3$ and $d_5$). As $m_3\cdot c_4\cdot m_4=m_5\cdot c_6 \cdot m_6$,  we know that the shorter of $m_3$ and $m_5$ must be a prefix of the longer, but we do not know the crucial concluding condition : that $m_3$ and $m_5$ have the same length. We could also deduce, with significant extra work\footnote{To do so, we would need to define a new version of grids which track empty arrows, then note that by existence and uniqueness of such grids, $p_3$ and $p_5$ lie within the same grid. Finally, we'd need to update the prefix-suffix fact in Theorem \ref{thm:nested_partial_grids_and_grids} to discuss such empty-arrow-tracking grids.}, that either $d_3$ is a suffix of $d_5$ or vice versa.  Alas, these facts about the middle and right frontier are not enough to conclude. (The fact that two partial grids have the same spine and bottom frontier, middle frontiers one of which is a prefix of the other, and right frontiers one of which is a suffix of the other is not enough to conclude that they have matching frames). 

The situation is even worse for $p_4$ and $p_6$, as their spines do not necessarily match. Now, this could be amended in a number of cases. If $d_5$ were shorter than $d_3$, we could extend the spine of $p_6$ to match $d_3$, as we know that extending a partial grid's spine preserves its size (however, this would have consequences for its bottom and middle frontiers). 

If $d_5$ were longer than $d_3$, we might be able to use an alternate form of induction, some sort of strong induction (which we will discuss below), to extend the spine of $p_4$ and then apply our more generalized induction hypothesis. However, this still does not solve the problem of the mismatched frontiers between $p_4$ and $p_6$, so we remain at an impasse.

Thus, Conjecture \ref{conj:same_type_same_length} remains inconclusive. It is possible to complete the project without Conjecture \ref{conj:same_type_same_length}, but is certainly displeasing. Intuitively, partial grids with the same frame contain the same (multi)set of cells, and should have the same size. We wonder if perhaps it is our definition of a partial grid, with its complicated division of the frontier into bottom, middle, and right segments, that is to blame. Here is an alternative, simpler definition of a partial grid : 

\begin{definition} \leanlink{https://github.com/hannahfechtner/braids_better/blob/098309403cda787eaf1fa954e07ff0ca3b2927ee/BraidProject/AlternativeChoices/FrontierStyle.lean#L11-L28}
    The set of frontier-style partial grids $G_F$ is a subset of $(BW^{\mathds{1}})^3$ defined inductively by
    \begin{itemize}
        \item a spine : For any positive, nonempty $a, b\in BW^{\mathds{1}}$, we have $\langle \overline{a}, b, \overline{a}b, \mathsf{spine}\rangle\in G_f$. Visually,

        \begin{tikzcd}
{} \arrow[d, "a"'] \arrow[r, "b"] & {} \\
{}                                &   
\end{tikzcd}
        
        \item cell addition : Given $(x, y)\in R_{\mathsf{grid}}$ and $f_1 =\langle a, b, c_1xc_2, \kappa_1\rangle\in G_F$, we have 
        
        $\langle a, b, c_1yc_2, \mathsf{cell \, addition}(x, y,f_1)\rangle\in G_F$. Visually, 

\begin{tikzcd}
{} \arrow[ddd, "a"'] \arrow[rrr, "b"] &                       &                                                & {} \\
                                      & & {} \arrow[ru, "c_2"', dash, squiggly]              &    \\
                                      & {} \arrow[ru, Mulberry, "y"{pos=1, below, yshift=-22pt}, dash, to path={ -| (\tikztotarget) \tikztonodes }]
  \arrow[ru, "x"{pos=0, above, yshift=25pt}, dash, to path={ |- (\tikztotarget) \tikztonodes }]      &                                             &    \\
{} \arrow[ru, "c_1"', dash, squiggly]                 &                       &                                                & 
\end{tikzcd}

    \end{itemize}

\end{definition}

\begin{definition} \leanlink{https://github.com/hannahfechtner/braids_better/blob/098309403cda787eaf1fa954e07ff0ca3b2927ee/BraidProject/AlternativeChoices/FrontierStyle.lean#L32-L40}
    We define the length metric $\ell_{G_F} : G_F \to \mathbb{N}$ as follows :
    \begin{itemize}
        \item a spine : $\ell_{G_{F}}\big(\langle \_, \_, \_, \mathsf{spine}\rangle\big)=0$
        
        \item cell addition : $\ell_{G_{F}}\bigg(\big\langle \_ , \_, \_, \mathsf{cell \, addition}(x, y,f_1)\big\rangle\bigg) = \ell_{G_F}(f_0)+\ell_{R_{\mathsf{grid}}}(x, y)$.
        
    \end{itemize}
    
\end{definition}

Alas, we again encounter issues trying to prove Conjecture \ref{conj:same_type_same_length}, even with $G_F$ in place of $G_P$. The issue is that frontier-style partial grids do not have a unique destructor. We wish we had a theorem of the following form :

\begin{fallacy}\label{fal:destructor}
    Given $(x,y)\in R_{\mathsf{grid}}$ and a frontier-style partial grid $f_1=\langle a, b, c_1yc_2, \kappa_1\rangle$, there is a frontier-style partial grid $f_0=\langle a, b, c_1xc_2, \kappa_0\rangle$ with $\ell_F(f_1)=\ell_F(f_0)+\ell_{R_{\mathsf{grid}}}(x,y)$.
    \begin{center}
    \begin{tikzcd}
{} \arrow[ddd, "a"'] \arrow[rrr, "b"] &                       &                                                & {} \\
                                      & & {} \arrow[ru, "c_2"', dash, squiggly]              &    \\
                                      & {} \arrow[ru, "y"{pos=1, below, yshift=-22pt}, dash, to path={ -| (\tikztotarget) \tikztonodes }]
  \arrow[ru, Rhodamine, dash, "x"{pos=0, above, yshift=25pt}, to path={ |- (\tikztotarget) \tikztonodes }]      &                                             &    \\
{} \arrow[ru, "c_1"', dash, squiggly]                 &                       &                                                &  
\end{tikzcd}
\end{center}
\end{fallacy}
In fact, Fallacy \ref{fal:destructor} does not even hold for the base cases. Consider the very simple frontier-style partial grid $f_1 = \langle\overline{\mathds{1}}, \mathds{1}, \mathds{1}\overline{\mathds{1}}, \mathsf{cell\, addition}(\overline{\mathds{1}}\mathds{1}, \mathds{1}\overline{\mathds{1}}, \langle \overline{\mathds{1}}, \mathds{1}, \overline{\mathds{1}}\mathds{1}, \mathsf{spine}\rangle)$ 

\begin{center}
\begin{tikzcd}
{} \arrow[d, "\varepsilon"'] \arrow[r, "\varepsilon"] & {} \arrow[d, "\varepsilon"] \\
{} \arrow[r, "\varepsilon"']                         & {}                        
\end{tikzcd}
\end{center}

We know that, say, $(\overline{\sigma_3}\sigma_3, \mathds{1}\overline{\mathds{1}})\in R_{\mathsf{grid}}$ and so our conjecture states that we should have some frontier-style partial grid
$f_0 = \langle \overline{\mathds{1}}, \mathds{1}, \sigma_3\overline{\sigma_3}, \kappa_0\rangle$. Clearly, such a partial grid does not exist. We can amend Fallacy \ref{fal:destructor} to the following statement, which is both tedious to prove and nearly useless. Hence, it is not formalized in Lean, and listed merely as a fact. 

\begin{fact}
    Given $(x,y)\in R_{\mathsf{grid}}$ and a frontier-style partial grid $f_1=\langle a, b, c_1yc_2, \kappa_1\rangle$, there is a frontier-style partial grid $f_0=\langle a, b, c_1zc_2, \kappa_0\rangle$ for some $z \in BW^{\mathds{1}}$ with $(z, y)\in R_{\mathsf{grid}}$, and $\ell_F(f_1)=\ell_F(f_0)+\ell_{R_{\mathsf{grid}}}(z,y)$.
\end{fact}

Due to the lack of a unique destructor, we are unable to prove that the frame of a frontier-style partial grid suffices to determine its size.

\begin{conjecture}
    Given frontier-style partial grids $f_1=\langle a, b, c, \kappa_1\rangle$ and $f_2=\langle a, b, c, \kappa_2\rangle$, we know $\ell_F(f_1)=\ell_F(f_2)$.
\end{conjecture}

Again, it is instructive to see where the straightforward proof attempt fails. Say we proceed by induction on $f_1$. Consider the case where $c=c_1\cdot \mathds{1}\overline{\mathds{1}} \cdot c_2$, we have $f_0 = \langle a, b, c_1\cdot\overline{\sigma_3}\sigma_3\cdot c_2\rangle$ to which the inductive hypothesis applies, and $f_1=\big\langle a, b, c_1\cdot \mathds{1}\overline{\mathds{1}}\cdot c_2, \mathsf{cell \, appending } (\overline{\sigma_3}\sigma_3, \mathsf{1}\overline{\mathsf{1}}, f_0) \big\rangle$. We would like to in some manner destruct $f_2$ to apply the inductive hypothesis, but alas, we have no unique destructor. For all we know, there is a frontier-style partial grid $f_3=\langle a, b, c_1\cdot \overline{\sigma_{42}}\sigma_{42}\cdot c_2, \kappa_3\rangle $ and $f_2=\langle a, b, c_1\cdot \mathds{1}\overline{\mathds{1}}, \mathsf{cell \; appending}(\overline{\sigma_{42}}\sigma_{42}, \mathds{1}\overline{\mathds{1}}, f_3)\rangle$. 

So, frontier-style partial grids suffer the same shortcoming as partial grids --- we have not been able to prove that their frame determines their size. Do they improve on our original definition of partial grids in any substantial way? In some ways, yes. For example, \leanlink{https://github.com/hannahfechtner/braids_better/blob/098309403cda787eaf1fa954e07ff0ca3b2927ee/BraidProject/AlternativeChoices/FrontierStyle.lean#L50-L66} the \lstinline{add_cell} function would be nearly trivial to implement for frontier-style partial grids, as compared to the seventeen involved cases for our original definition.
In other ways, certainly not. I could not find a proof using frontier-style partial grids for Theorem \ref{thm:grid_of_pg} (if a partial grid's frontier is in negative-positive form, then it is in fact a grid). Thus, we would need to keep both structures and transfer theorems between them. \leanlink{https://github.com/hannahfechtner/braids_better/blob/098309403cda787eaf1fa954e07ff0ca3b2927ee/BraidProject/AlternativeChoices/FrontierStyle.lean#L90-L338} The work of setting up a (size-preserving) bijection between $G_P$ and $G_F$ is at least equal to the work needed to define \lstinline{add_cell} directly on $G_P$, and so we dispense with frontier-style partial grids for good. 

The enterprising reader is most welcome to formulate his own definition for partial grids. Here are three theorems such a definition must satisfy :

\begin{itemize}
    \item Theorem \ref{thm:grid_of_pg} : if a partial grid's frontier is in negative-positive form, then it is in fact a grid.
    \item Definition \ref{def:new_add_cell} : there exists a function with the same specification as \lstinline{add_cell}.
    \item Theorem \ref{thm:new_partial_grid_of_rw} : if $\overline{a}b$ rewrites to $c$, then there is a partial grid with spine $(a, b)$ and frontier $c$.
\end{itemize}

The ideal definition would also satisfy Conjecture \ref{conj:same_type_same_length}; that is, the frame of the partial grid determines its size.

This is not the only unwieldy issue with the partial grids we have defined. We do not have a rigorous definition of a sub-partial-grid, which could provide a method of strong induction on partial grids. Surprisingly, there is not a single proof in this project requiring strong induction on partial grids. Nevertheless, one might find it easier to use than structural induction in some cases, or perhaps give a new theorem, or a proof of a conjecture, requiring strong induction.

Here are a few simple ideas for inducting on a partial grid $\langle a, b, c, m, d\rangle$, none of which work :
\begin{itemize}
    \item We might try to induct on the length of the spine (that is, $|a|+|b|$). However, due to cells of the form \begin{tikzcd}
{} \arrow[rr, "\sigma_j"] \arrow[dd, "\sigma_i"'] &                           & {} \arrow[d, "\sigma_i"] \\
                                                  &                           & {} \arrow[d, "\sigma_j"] \\
{} \arrow[r, "\sigma_j"']                         & {} \arrow[r, "\sigma_i"'] & {}                      
\end{tikzcd}, this metric often increases when we split a partial grid. The issue persists even if we ignored the empty cells by using $|s^{\mathds{1}}(a)|+|s^{\mathds{1}}(b)|$ for our metric
\item We might try to induct on the diagonal length of the partial grid, $|a|+|c|$, or possibly $|b|+|d|$ (or even the sum of the two). However, this metric remains static when adding a cell that touches only the middle frontier. 
\item We might try to induct on the size of the partial grid using $\ell_P$. However, due to the presence of trivial cells, it is possible to split a partial grid and obtain one part with size zero, and one part with the original size.
\end{itemize}

The last option could be amended to be workable : we could change our metric to $\ell_P'$, which gives every cell size $1$ : 
\begin{definition}
    Define an alternative metric on partial grids $\ell_P'$ identically to $\ell_P$, except for the very first case, which we replace by
    \begin{itemize}
        \item $\ell_P \big(\langle a, b, c, m, d,\mathsf{trivial \, cell}\rangle \big) = 1$
    \end{itemize}
\end{definition}

This would almost work for strong induction, except trimming off tails (such as arises in the second case of splittability) does not decrease $\ell_P'$. It could be perhaps amended to use the length of the spine as a tie-breaker (formally, a lexicographic ordering on the pair $\langle \ell'_P, \mathsf{spine \, length\rangle}$). This new metric is less elegant, especially as $\ell'_P$ does not match nicely with $\ell_{S_{\mathsf{rev}}}$.

So, even the simplest methods for defining strong induction are reasonably complex (in order to derive strong induction from the just-mentioned lexicographic order, one would have to consider a grid $g$ and a smaller (by lexicographic ordering) grid $g_1$, and then repeatedly split $g_1$ and apply determinative spines, grid uniqueness, and grid appending in order to obtain grids to which the structural inductive hypothesis would apply). Strong induction is never used in this project, and so we do not formalize $\ell_P'$ nor the associated lexicographic ordering.

In sum, we have defined metrics on our main structures (Semi-Thue derivations, grids, and partial grids). We have been able to show a number of properties of each metric, properties which will let us prove termination of the subword reversing algorithm in the following chapter.

\section{The Reversing Procedure}

We are now prepared to implement Dehornoy's subword reversing algorithm, which determines whether a braid group element is equivalent to the identity element. (Of course, to check equivalence between elements $a$ and $b$, we need simply check if $a^{-1}b$ is equivalent to the identity.) The heart of this algorithm is a process called ``reversing'' which rewrites any braid group element to one in a certain standard form, preserving equivalence in the braid group.
As an intuitive overview, a string in $BW$ undergoes successive rewritings according to $R_{\mathsf{rev}}$ until no such rewrites remain to be performed. We begin with ``simple'' reversing, which takes in a braid word in positive-negative form. We then extend simple reversing to full reversing, which takes in any braid word. Finally, we implement the full subword reversing algorithm, which reverses a word, exchanges the positive and negative parts, and then reverses again.

\subsection{Monoid Solver}

We will first implement a solver for the braid monoid word problem. Much of the machinery will be expanded and re-used for the braid group word problem solver. We begin by defining simple reversing, in which we iteratively apply reversing rewrites to $\overline{a}b$ for $a,b\in BW^+$. (Note that many of our functions are given the superscript $^+$, as this simple version begins with two elements of the braid \textit{monoid}). Our eventual solver will include simple reversing as its main step.

To aid in exposition, we give terminology to the first element of any $R_{\mathsf{rev}}$ pair as follows : 

\begin{definition}
    For any $i,j\in\mathbb{N}$, we say that the braid word $\overline{\sigma_i}\sigma_j$ is an \textbf{open pair}. 
\end{definition}
As a first attempt, we might define said iterative rewriting
  $\mathsf{rev}_1^{+} : BW \to BW$
recursively as so: If there exist \( c_1, c_2 \in BW \) such that for some \( i, j \in \mathbb{N} \),
\( c = c_1 \cdot \overline{\sigma_i}\sigma_j \cdot c_2 \) where $c_1$ contains no open pair,
then
\[
  \mathsf{rev}_1^+(c)
  = \mathsf{rev}_1^+(c_1 \cdot \overline{\sigma_i}\sigma_j \cdot c_2)
  :=
  \begin{cases}
    \mathsf{rev}_1^+(c_1 \cdot c_2) & i = j, \\[6pt]
    \mathsf{rev}_1^+(c_1 \cdot \sigma_j \overline{\sigma_i} \cdot c_2) & |i - j| \ge 2, \\[6pt]
    \mathsf{rev}_1^+(c_1 \cdot \sigma_j \sigma_i \overline{\sigma_j} \overline{\sigma_i} \cdot c_2) & |i - j| = 1.
  \end{cases}
\]
Otherwise,
\[
  \mathsf{rev}_1^+(c) = c.
\]

The immediate problem is to show well-definedness of this function -- that is,
that the recursive calls in $\mathsf{rev}_1^{+}$ do not go on endlessly\footnote{The replacement may introduce new negative-positive pairs, so termination is not straightforward}.
This will require the fact that our original input is of the form $\overline{a}\, b$, and $\overline{a}\, b$ rewrites under $R_{\mathsf{rev}}$ to the output of $\mathsf{rev}_1^{+}$ : 
\[
  \overline{a}\, b
\underset{R_{\mathsf{rev}}}{\rightarrowtriangle}
  \mathsf{rev}_1^{+}\!\left(\overline{a}\, b\right).
\]
That is immediately circular -- how can we speak about the output of a function
before we know it terminates? Instead, we build this fact into the definition of the function by restricting
the domain and codomain of $\mathsf{rev}^{+}$.

\begin{definition} \leanlink{https://github.com/hannahfechtner/braids_better/blob/098309403cda787eaf1fa954e07ff0ca3b2927ee/BraidProject/Solver/Monoid.lean#L40-L41}
   Define the sets $T$ parametrized by $a,b \in BW^{+}$ as\footnote{We use the letter $T$ as any element of $T_{a,b}$ corresponds to a partial grid with a potentially jagged frontier, hence roughly in in the shape of a triangle.}
\[
  T_{a,b} = \{\, c \in BW \mid  \overline{i^{+}(a)}\, i^{+}(b) \underset{R_{\mathsf{rev}}}{\rightarrowtriangle} c
  ) \,\}.
\]
Here we are briefly explicit about $i^+$ for precision. Note that we are using the data-carrying definition of $\underset{R_{\mathsf{rev}}}{\rightarrowtriangle}$. Our formalization defines elements of $T_{a,b}$ as a pair $(x,y)$, where $x\in BW$ and $y:=\langle \overline{i^{+}(a)}\, i^{+}(b)\underset{R_{\mathsf{rev}}}{\rightarrowtriangle}
  x, \kappa\rangle$. We thus use this convenient notation in our writeup.

\end{definition}

Then we inductively define the family of functions
\[
  \mathsf{rev}^{+}_{a,b} : T_{a,b} \to T_{a,b}
\]
as follows; note that the only differences from the first definition are the name of the function and the domain and range (and hence the added derivation). 

\begin{definition} \leanlink{https://github.com/hannahfechtner/braids_better/blob/098309403cda787eaf1fa954e07ff0ca3b2927ee/BraidProject/Solver/Monoid.lean#L53-L81}
    Given $a, b\in BW^+$, $\mathsf{rev}^{+}_{a,b} : T_{a,b} \to T_{a,b}$ is defined by

    \textbf{(Recursive case)} If there exist $x_1, x_2 \in BW$ such that for some
$i,j \in \mathbb{N}$,
$x = x_1 \cdot \overline{\sigma_i}\sigma_j \cdot x_2$ where $x_i$ contains no open pair, then
\[
  \mathsf{rev}^{+}_{a,b}(x_1 \cdot \overline{\sigma_i}\sigma_j \cdot x_2, y)
  =
\begin{cases}
\begin{aligned}
&\mathsf{rev}^{+}_{a,b}\Big(
    x_1 \cdot x_2,\\
&\qquad\qquad
    \big\langle
    \overline{i^{+}(a)}\, i^{+}(b)
    \underset{R_{\mathsf{rev}}}{\rightarrowtriangle}
    x_1\cdot x_2,\\
&\qquad\qquad\qquad
    \mathsf{step\,transitivity}
    (y,\overline{\sigma_i}\overline{\sigma_j},\varepsilon)
    \big\rangle
\Big)
\end{aligned}
& i=j,
\\[10pt]
\begin{aligned}
&\mathsf{rev}^{+}_{a,b}\Big(
    x_1 \cdot \sigma_j \sigma_i
    \overline{\sigma_j}\,\overline{\sigma_i}\cdot x_2,\\
&\qquad\qquad
    \big\langle
    \overline{i^{+}(a)}\, i^{+}(b)
    \underset{R_{\mathsf{rev}}}{\rightarrowtriangle}
    x_1 \cdot \sigma_j \sigma_i
    \overline{\sigma_j}\,\overline{\sigma_i}\cdot x_2,\\
&\qquad\qquad\qquad
    \mathsf{step\,transitivity}\big(
    y,\overline{\sigma_i}\overline{\sigma_j},
    \sigma_j \sigma_i
    \overline{\sigma_j}\,\overline{\sigma_i}
    \big)
    \big\rangle
\Big)
\end{aligned}
& |i-j|=1,
\\[10pt]
\begin{aligned}
&\mathsf{rev}^{+}_{a,b}\Big(
    x_1 \cdot \sigma_j \overline{\sigma_i}\cdot x_2,\\
&\qquad \qquad
    \big\langle
    \overline{i^{+}(a)}\, i^{+}(b)
    \underset{R_{\mathsf{rev}}}{\rightarrowtriangle}
    x_1 \cdot \sigma_j \overline{\sigma_i}\cdot x_2,\\
&\qquad\qquad\qquad
    \mathsf{step\,transitivity}\big(
    y,\overline{\sigma_i}\overline{\sigma_j},
    \sigma_j\overline{\sigma_i}
    \big)
    \big\rangle
\Big)
\end{aligned}
& |i-j|>1.
\end{cases}
\]

\textbf{(Non-recursive case)} Otherwise,
\[
  \mathsf{rev}^{+}_{a,b}(x,y) = (x,y).
\]
\end{definition}

Verifying the domain and range for the different cases in the definition of $\mathsf{rev}^+_{a,b}$ is easily done. 

For our formalization, we write a simple recursive function, \lstinline{FindOpenPair}, \leanlink{https://github.com/hannahfechtner/braids_better/blob/098309403cda787eaf1fa954e07ff0ca3b2927ee/BraidProject/Solver/FindOpenPair.lean#L7-L19} to find the first occurrence of a substring of the form $\overline{\sigma_i}\,\sigma_j$. However, we cannot use the existing \lstinline{List.Infix} infrastructure in Mathlib as it is non-data-carrying. We thus develop a small repository of data-carrying proofs about basic List properties. \leanlink{https://github.com/hannahfechtner/braids_better/blob/098309403cda787eaf1fa954e07ff0ca3b2927ee/BraidProject/DataCarrying/List.lean}

Next, we must show termination of $\mathsf{rev}^+_{a,b}$. We first define a metric on elements of $T_{a,b}$ (what all that data-carrying work has been for!).

\begin{definition} \leanlink{https://github.com/hannahfechtner/braids_better/blob/098309403cda787eaf1fa954e07ff0ca3b2927ee/BraidProject/Solver/Monoid.lean#L47}
    We define the family of functions $\ell_{T_{a, b}} : T_{a,b} \to \mathbb{N}$ (each often written $\ell_T$ when $a$ and $b$ can be inferred). Given $(x,y)\in T_{a,b}$, we know $x\in BW$ and we have a derivation $y= x \underset{R_{\mathsf{rev}}}{\rightarrowtriangle}
  \overline{i^{+}(a)}\, i^{+}(b) $. Then we define
    $$\ell_T(x,y):=\ell_B(a,b)-\ell_S(y)$$
\end{definition}

\begin{theorem}\label{thm:metric_nontruncated} \leanlink{https://github.com/hannahfechtner/braids_better/blob/098309403cda787eaf1fa954e07ff0ca3b2927ee/BraidProject/Solver/Monoid.lean#L21-L38}
    For any $(x,y)\in T_{a,b}$, we have $\ell_S(y)\leq \ell_B(a,b)$.
\end{theorem}
\begin{proof}
    Consider first the case for nonempty $a,b$. Since $y= \overline{i^{+}(a)}\, i^{+}(b) \underset{R_{\mathsf{rev}}}{\rightarrowtriangle} x
   $, by Theorem \ref{thm:partial_grid_of_rw} we know there is a partial grid $p = \langle a, b, c, m, d, \kappa\rangle$ with $x = s^{\mathds{1}}(c m d)$ and $\ell_P(p)=\ell_S(y)$. We also know that $\ell_P(p)\leq\ell_B(a,b)$ by definition of $\ell_B$ (Definition \ref{def:pair_length}) and Theorem \ref{thm:nested_partial_grids_and_grids}. Hence, $\ell_S(y)\leq\ell_B(a,b)$ as desired.

   When either of $a$ or $b$ is empty, we know that $\overline{i^{+}(a)}\, i^{+}(b)$ contains no open pair, and thus we have $\ell_S(y)=0$, satisfying our desired inequality.
\end{proof}
\begin{theorem}\label{monoid_terminates} \leanlink{https://github.com/hannahfechtner/braids_better/blob/098309403cda787eaf1fa954e07ff0ca3b2927ee/BraidProject/Solver/Monoid.lean#L75-L81}
    The metric $\ell_T$ decreases in each recursive call in $\mathsf{rev}^+_{a,b}$.
\end{theorem}
\begin{proof}
    By Theorem \ref{thm:metric_nontruncated}, it suffices to show that $\ell_S(y)$ increases in each recursive call. This is immediate, as each recursive call generates a new derivation using step transitivity, which increases $\ell_S(y)$ by the size of the relation used, and $\ell_{R_{\mathsf{rev}}}$ is always 1. Hence, in each recursive call, $\ell_S(y)$ increases by 1, and thus $\ell_B(a,b)-\ell_S(y)$ decreases by 1. 
\end{proof}

\medskip

As $\ell_T$ is a natural-number metric, \leanlink{https://github.com/hannahfechtner/braids_better/blob/098309403cda787eaf1fa954e07ff0ca3b2927ee/BraidProject/Solver/Monoid.lean#L74}it is well-founded, and hence our iterative rewriting function $\mathsf{rev}^{+}$ terminates. We finally write a wrapper function which will be used to solve the word problem in the monoid.

\begin{definition} \leanlink{https://github.com/hannahfechtner/braids_better/blob/098309403cda787eaf1fa954e07ff0ca3b2927ee/BraidProject/Solver/Monoid.lean#L128-L129}
We define simple reversing of a monoid pair $r^{+} : BW^{+} \to BW^{+} \to BW$ by
\[
  r^{+}(a,b)
  = \mathsf{rev}^{+}_{a,b}\!\left(\overline{a}\, b, \Big\langle \overline{a}\, b \underset{R_{\mathsf{rev}}}{\rightarrowtriangle} \overline{a}\, b, \mathsf{refl}\Big\rangle\right)
\]

We define a monoid solver $\mathsf{solve^+} : (BW^+)^2\to \mathsf{Bool}$ by $\mathsf{solve^+} (a,b)= \mathsf{true}$ iff $r^+(a,b)=\varepsilon$.
\end{definition}

\noindent We amass a number of theorems about $r^+$ in order to show correctness of $\mathsf{solve^+}$.
\begin{lemma} \label{lem:reversing_rw_rev} \leanlink{https://github.com/hannahfechtner/braids_better/blob/098309403cda787eaf1fa954e07ff0ca3b2927ee/BraidProject/Solver/Monoid.lean#L136-L137}
$\overline{a}\, b
  \xrightarrow[R_{\mathrm{rev}}]{*}
  r^{+}(a,b) $ holds immediately.
\end{lemma}

We next show that the output of $r^+(a,b)$ is in positive-negative form.

\begin{lemma} \leanlink{https://github.com/hannahfechtner/braids_better/blob/098309403cda787eaf1fa954e07ff0ca3b2927ee/BraidProject/Solver/Monoid.lean#L83-L124}
    For $a,b\in BW^+$, $r^+(a,b)$ contains no open pair (substring of the form $\overline{\sigma_i}\sigma_j$).
\end{lemma}
\begin{proof}
    We show that for any $x\in BW$, $y=\langle \overline{a}b \underset{R_{\mathsf{rev}}}{\rightarrowtriangle}x, \kappa\rangle$  we must have that $\mathsf{rev}^+_{a,b}(x,y)$ contains no open pair. We proceed by induction on $(x,y)\in T_{a,b}$ using strong induction on $\ell_{T_{a,b}}$. We then split into the four cases in the definition of $\mathsf{rev}^+_{a,b}$. Let's look at one illustrative case. Say we have $x=x_1\cdot\overline{\sigma_i}\sigma_j\cdot x_2$ with $|i-j|>1$, and thus $\mathsf{rev}^+_{a,b}(x,y) = \mathsf{rev}^+_{a,b}(x_1\cdot\sigma_j\overline{\sigma_i}\cdot x_2,
    \big\langle
    \overline{i^{+}(a)}\, i^{+}(b)
    \underset{R_{\mathsf{rev}}}{\rightarrowtriangle}
    x_1 \cdot \sigma_j \overline{\sigma_i}\cdot x_2,
    \mathsf{step\,trans}\big(
    y,\overline{\sigma_i}\overline{\sigma_j},
    \sigma_j\overline{\sigma_i}
    \big)
    \big\rangle)$. 

    Note that $\ell_{T_{a,b}}(x,y)=1+\ell_{T_{a,b}}((x_1\cdot\sigma_j\overline{\sigma_i}\cdot x_2,
    \big\langle
    \overline{i^{+}(a)}\, i^{+}(b)
    \underset{R_{\mathsf{rev}}}{\rightarrowtriangle}
    x_1 \cdot \sigma_j \overline{\sigma_i}\cdot x_2,
    \mathsf{step\,trans}\big(
    y,\overline{\sigma_i}\overline{\sigma_j},
    \sigma_j\overline{\sigma_i}
    \big)
    \big\rangle))$, and so we apply our inductive hypothesis to conclude that $\mathsf{rev}^+_{a,b}(x_1\cdot\sigma_j\overline{\sigma_i}\cdot x_2,
    \big\langle
    \overline{i^{+}(a)}\, i^{+}(b)
    \underset{R_{\mathsf{rev}}}{\rightarrowtriangle}
    x_1 \cdot \sigma_j \overline{\sigma_i}\cdot x_2,
    \mathsf{step\,trans}\big(
    y,\overline{\sigma_i}\overline{\sigma_j},
    \sigma_j\overline{\sigma_i}
    \big)
    \big\rangle)$, and thus $\mathsf{rev}^+_{a,b}(x,y)$ contains no open pair. Finally, we conclude that $r^+(a,b)$ contains no open pair as desired.
\end{proof}

\begin{lemma} \leanlink{https://github.com/hannahfechtner/braids_better/blob/098309403cda787eaf1fa954e07ff0ca3b2927ee/BraidProject/Solver/FindOpenPair.lean#L154-L174}
If a word $w\in BW$ contains no open pairs, then it is in positive-negative form -- that is, there exist $c,d\in BW^+$ such that $w=c\overline d$. (Proof by induction on $w$.)
\end{lemma}

Finally, we put together the two previous lemmas to obtain our desired theorem :
\begin{theorem} \label{thm:reversing_pos_neg} \leanlink{https://github.com/hannahfechtner/braids_better/blob/098309403cda787eaf1fa954e07ff0ca3b2927ee/BraidProject/Solver/Monoid.lean#L131-L134}
    For any $a,b\in BW^+$, there exist $c,d \in BW^+$ such that $r^+(a,b)=c\overline{d}$.
\end{theorem}

We can now prove a correctness theorem for our monoid solver: 
\begin{theorem}\label{thm:monoid_correctness}
    $\mathsf{solve^+}(a,b)=\mathsf{true}$ if and only if $a\equiv^+ b$
\end{theorem}

\begin{proof}
    Forward direction: \leanlink{https://github.com/hannahfechtner/braids_better/blob/098309403cda787eaf1fa954e07ff0ca3b2927ee/BraidProject/Solver/MonoidCorrectness.lean#L9-L30}Assume $\mathsf{solve^+}(a,b)=\mathsf{true}$, so $r^+(a,b)=\varepsilon$. Then by Lemma \ref{lem:reversing_rw_rev}, $\overline{a}b \xrightarrow[R_{\mathrm{rev}}]{*}\varepsilon$. Thus by Theorem \ref{lem:reverse_equiv}, $\overline{a}b \equiv\varepsilon$, so $a\equiv b$ and by Ore localization \ref{thm:group_equiv_monoid_equiv_ore}, $a\equiv^+ b$

    Backward direction: \leanlink{https://github.com/hannahfechtner/braids_better/blob/098309403cda787eaf1fa954e07ff0ca3b2927ee/BraidProject/Solver/MonoidCorrectnessHardDirection.lean#L13-L37}Assume $a\equiv^+ b$. Then by Theorem \ref{thm:grid_one_one_iff_equiv}, there is a grid with source $(a,b)$ and target $(\varepsilon, \varepsilon)$. Further, by Theorem \ref{thm:reversing_pos_neg} we know that for some $c,d\in BW^+$, $r^+(a,b)=c\overline{d}$. Hence by Lemma \ref{lem:reversing_rw_rev}, $\overline{a}b \xrightarrow[R_{\mathrm{rev}}]{*} c\overline{d}$. By Theorem \ref{thm:big_one}, there is a grid with source $(a,b)$ and target $(c,d)$. By uniqueness of grid sinks (Theorem \ref{thm:grid_uniqueness}), $c=\varepsilon$ and $d=\varepsilon$. Thus $r^+(a,b)=\varepsilon$ and $\mathsf{solve^+}(a,b)=\mathsf{true}$.
\end{proof}

\leanlink{https://github.com/hannahfechtner/braids_better/blob/098309403cda787eaf1fa954e07ff0ca3b2927ee/BraidProject/Solver/MonoidCorrectness.lean#L33-L43}In fact, for the forward direction, there is a method which does not use Ore's theorem, but rather works through grids. Since $\overline{a}b \xrightarrow[R_{\mathrm{rev}}]{*}\varepsilon\,(=\varepsilon\overline{\varepsilon})$, by Theorem \ref{thm:big_one}, there is a grid from $(a,b)$ to $(\varepsilon, \varepsilon)$, and hence by Theorem \ref{thm:grid_equiv} $a\equiv^+b$.

We also have an alternative proof for the backward direction, which uses a sort of restricted confluence for Semi-Thue systems over $R_{\mathsf{rev}}$. We begin with two theorems about partial grids ``nested'' inside grids which share essentially the same spine (up to deletion of $\mathds{1}$ and $\overline{\mathds{1}}$), which will aid in the confluence proof :

\begin{theorem} \leanlink{https://github.com/hannahfechtner/braids_better/blob/098309403cda787eaf1fa954e07ff0ca3b2927ee/BraidProject/PartialGrid/FrontierToSink.lean#L11-L233}
    Let $g = \langle a, b, c, d\rangle\in G$ and $p = \langle a_0, b_0, c_0, m_0, d_0\rangle \in G_P$. Recall that $s^{\mathds{1}}$ deletes all instances of $\mathds{1}$ or $\overline{\mathds{1}}$ from a braid word. If $a = s^{\mathds{1},+}(\overline {a_0})$, $b=s^{\mathds{1},+}({b_0})$, and there exists $c_1, d_1 \in BW^+$ such that $c = s^{\mathds{1},+} (c_0)\cdot c_1$ and $d = s^{\mathds{1},+}(\overline {d_0})\cdot d_1$, then we have $s^{\mathds{1}}(m_0) \xrightarrow[R_{\mathrm{rev}}]{*} c_1 \cdot \overline{d_1}$.

Here is the image to keep in mind -- we claim, roughly, that the purple section rewrites under $R_{rev}$ to the green section, both read from the bottom-left to the top-right : 
\begin{center}
\begin{tikzcd}
{} \arrow[dd, "a"'{xshift=-4pt}] 
   \arrow[rr, "b"{yshift=4pt}] 
   \arrow[rr, "b_0"', draw=none,
     preaction={draw=magenta, draw opacity=.35, line width=4pt, -},
     text=magenta]
   \arrow[dd, "a_0", draw=none,
     preaction={draw=magenta, draw opacity=.35, line width=4pt, -},
     text=magenta]
&
&
{} \arrow[dd, "d"{xshift=4pt}] 
   \arrow[d, "d_0"', draw=none,
     preaction={draw=magenta, draw opacity=.35, line width=4pt, -},
     text=magenta]
\\
&
&
 {} \arrow[d, "d_1", draw=none,
     preaction={draw=green, draw opacity=.35, line width=2pt, -},
     text=ForestGreen]
\\
{} \arrow[rr, "c"'{yshift=-4pt}] 
   \arrow[r, "c_0", draw=none,
     preaction={draw=magenta, draw opacity=.35, line width=4pt, -},
     text=magenta]
&
{} \arrow[ru, "m_0",draw=none,
     preaction={draw=violet, draw opacity=.35, line width=4pt, shorten <=-4pt, shorten >=-4pt},
     text=violet, dash, squiggly]
   \arrow[r, "c_1"', draw=none,
     preaction={draw=green, draw opacity=.35, line width=2pt, -},
     text=ForestGreen]
&
{}
\end{tikzcd}
\end{center}
 
\end{theorem}

The proof follows the strategy of Theorem \ref{thm:nested_partial_grids_and_grids}, inducting on the grid and proceeding by extensive casework in each inductive case, splitting the partial grid as appropriate, and using prefix-suffix fact (1) from Theorem \ref{thm:nested_partial_grids_and_grids} to help apply the inductive hypothesis. 

From the above theorem, we can immediately deduce our second needed theorem, that the middle frontier of a partial grid reverses to the sink of the grid in which it is nested : 

\begin{theorem}\label{thm:middle_to_sink} \leanlink{https://github.com/hannahfechtner/braids_better/blob/098309403cda787eaf1fa954e07ff0ca3b2927ee/BraidProject/PartialGrid/FrontierToSink.lean#L243-L250}
Let $g = \langle a, b, c, d\rangle\in G$ and $p = \langle a_0, b_0, c_0, m_0, d_0\rangle \in G_P$. If $a = s^{\mathds{1},+}(\overline {a_0})$ and $b=s^{\mathds{1},+}({b_0})$, then we have $s^{\mathds{1}}(c_0m_0d_0) \xrightarrow[R_{\mathrm{rev}}]{*} c \overline d$.
\end{theorem}

Now we are prepared to state and prove restricted confluence. 

\begin{theorem} \label{thm:restricted_confluence} \leanlink{https://github.com/hannahfechtner/braids_better/blob/098309403cda787eaf1fa954e07ff0ca3b2927ee/BraidProject/PartialGrid/FrontierToSink.lean#L269-L292}
    Let $a,b\in BW^+$ and $w_1,w_2\in BW$. If $\overline{a}b\xrightarrow[R_{\mathrm{rev}}]{*} w_1$ and $\overline{a}b\xrightarrow[R_{\mathrm{rev}}]{*} w_2$, then there exists some $w\in BW$ such that $w_1\xrightarrow[R_{\mathrm{rev}}]{*} w$ and $w_2\xrightarrow[R_{\mathrm{rev}}]{*} w$.
\end{theorem}
\begin{proof}
    First consider the case of nonempty $a$ and $b$. We know there are $c,d\in BW^+$ such that $\langle a, b, c, d\rangle \in G$, by existence of grids (Theorem \ref{thm:grid_existence}). Let $w=c\overline{d}$. Since $\overline{a}b\xrightarrow[R_{\mathrm{rev}}]{*} w_1$, there is a partial grid with spine $(a,b)$ and total frontier $t_1$ such that $s^{\mathds{1}} = w_1$, so by Theorem \ref{thm:middle_to_sink}, $s^{\mathds{1}}(t_1) \xrightarrow[R_{\mathrm{rev}}]{*} c \overline d$, and thus $w_1\xrightarrow[R_{\mathrm{rev}}]{*} c \overline d$ as desired. We show that $w_2 \xrightarrow[R_{\mathrm{rev}}]{*} c \overline d$ in exactly the same manner.

    In the case where either $a$ or $b$ is empty, $\overline{a}b$ contains no open pairs, so $\overline{a}b=w_1$ and $\overline{a}b=w_2$. Thus, we may use $w_1$ for $w$, and both desired derivations are instances of reflexivity.
\end{proof}

\leanlink{https://github.com/hannahfechtner/braids_better/blob/098309403cda787eaf1fa954e07ff0ca3b2927ee/BraidProject/Solver/MonoidCorrectnessHardDirection.lean#L40-L54}Now we can return to our second proof of the reverse direction of the monoid solver's correctness from Theorem \ref{thm:monoid_correctness}. Assume that $a\equiv^+ b$. Then by Theorem \ref{thm:grid_one_one_iff_equiv}, there is a grid with source $(a,b)$ and target $(\varepsilon, \varepsilon)$. Thus, by the forward direction of \ref{thm:big_one}, $\overline{a}b \xrightarrow[R_{\mathrm{rev}}]{*} \varepsilon$. Since $\overline{a}b \xrightarrow[R_{\mathrm{rev}}]{*} r^+(a,b)$, by Theorem \ref{thm:restricted_confluence}, there exists some $c\in BW$ such that $\varepsilon \xrightarrow[R_{\mathrm{rev}}]{*} c$ (and thus $c = \varepsilon$) and $r^+(a,b)\xrightarrow[R_{\mathrm{rev}}]{*} c$. Since $r^+(a,b)$ is in positive-negative form, we must have $r^+(a,b)=c=\varepsilon$ as desired.

This second proof is not very different from the first, but substitutes one instance of grids for restricted confluence.

We thus have created both a solver for the braid monoid ($\mathsf{solve^+}$), as well as a method for reversing any word in negative-positive form ($\mathsf{rev^+}$).
\subsection{Group Solver}

We now work to implement a solver for the word problem in the braid group. This again largely relies on a method to iteratively apply reversing rewrites to a braid word, akin to $\mathsf{rev^+}$. However, we may no longer restrict the input word to one in negative-positive form.  We will still need to restrict the codomain of this new reversing function. 
 \begin{definition} \leanlink{https://github.com/hannahfechtner/braids_better/blob/098309403cda787eaf1fa954e07ff0ca3b2927ee/BraidProject/DataCarrying/SignedList.lean#L23}
     We define the set of braid words in positive-negative form as $BW^{PN}$ -- that is, the set of all $w\in BW$ of the form $c\overline{d}$ for $c, d\in BW^+$.
 \end{definition}

 \begin{definition} \leanlink{https://github.com/hannahfechtner/braids_better/blob/098309403cda787eaf1fa954e07ff0ca3b2927ee/BraidProject/Solver/Group.lean#L221-L271}
 We define the reversing function $\mathsf{rev} : BW \to BW^{PN}$ as follows.\footnote{Here, unlike in our definition of $\operatorname{rev^+}$, we do not restrict the initial input in any manner, so we will not need any wrapper function akin to $r^+$.} Consider some $w\in BW$. We know that there must exist some $u_0,v_0\in BW^+$, not both the empty string, and some $w'\in BW$ such that $w =  \overline{u_0}v_0w'$.\footnote{In Lean, this is implemented as \lstinline{separate_first_pair} \leanlinkfootnote{https://github.com/hannahfechtner/braids\_better/blob/098309403cda787eaf1fa954e07ff0ca3b2927ee/BraidProject/Solver/SeparatePair.lean\#L64-L67}} Note that $w'$ has length strictly less than $w$.

 Let $c,d \in BW^+$ be such that $r^+(u_0,v_0) = c\overline{d}$. Let our recursive call $\mathsf{rev}(w')$ return $e\overline f$. Then define $\mathsf{rev}(w)=c\cdot r^+(d,e) \cdot \overline{f}$. Note that this is indeed in positive-negative form.
 
\end{definition}

 We can show that $\mathsf{rev}$ preserves equivalence via induction and appeals to the behavior of $r^+$ : 

\begin{theorem}\label{thm:reverse_equiv} \leanlink{https://github.com/hannahfechtner/braids_better/blob/098309403cda787eaf1fa954e07ff0ca3b2927ee/BraidProject/Solver/Group.lean#L14-L219}
    For any $w\in BW$, we have $w \xrightarrow[R_{\mathrm{rev}}]{*} \mathsf{rev}(w)$.
\end{theorem}

For our Lean formalization, it is easier (but not strictly necessary) to bundle this fact into the definition of $\mathsf{rev}$. Thus, we define the codomain of $\mathsf{rev}$ as follows :

\leanlink{https://github.com/hannahfechtner/braids_better/blob/098309403cda787eaf1fa954e07ff0ca3b2927ee/BraidProject/Solver/Group.lean#L8-L11}
\begin{lstlisting}
structure ReverseResult (L : List (ℕ × Bool)) where
  result : List (ℕ × Bool)
  ordered : SignedList.PosNegData result
  derivation : SemiThueData reversing L result
\end{lstlisting}

That is, the codomain is the set of braid words (\lstinline{result}) in positive-negative order (\lstinline{ordered}) such that $w \xrightarrow[R_{\mathrm{rev}}]{*} \mathsf{rev}(w)$ (\lstinline{derivation}).

With reversing for words in $BW$ defined, we may move on to the word problem in the braid group. Our group solver will be a bit more complex than our monoid solver -- it does not suffice to check if the output of $\mathsf{rev}$ is the empty string.

\begin{definition} \leanlink{https://github.com/hannahfechtner/braids_better/blob/098309403cda787eaf1fa954e07ff0ca3b2927ee/BraidProject/Solver/Group.lean#L273-L275}
    We define a procedure $\mathsf{solve} : BW^2 \to \mathsf{Bool}$ for determining equality in the braid group. Given $w_1, w_2 \in BW$, let $\mathsf{rev}(\overline{w_1}w_2)=c\overline{d}$. Then, $\mathsf{solve}(w_1, w_2) = \mathsf{true}$ iff $\mathsf{solve^+}(d,c) = \mathsf{true}$
    (that is, iff $r^+(d,c)=\varepsilon$).
\end{definition}

Next we show correctness of this procedure. 

\begin{theorem} \leanlink{https://github.com/hannahfechtner/braids_better/blob/098309403cda787eaf1fa954e07ff0ca3b2927ee/BraidProject/Solver/GroupCorrectnessHardDirection.lean#L26-L28}
    For any $w_1, w_2\in BW$, $\operatorname{solve}(w_1,w_2)=\mathsf{true}$ if and only if $w_1\equiv w_2$.
\end{theorem}

\begin{proof}
    Forward direction \leanlink{https://github.com/hannahfechtner/braids_better/blob/098309403cda787eaf1fa954e07ff0ca3b2927ee/BraidProject/Solver/GroupCorrectness.lean#L37-L57}: Assume $\mathsf{solve}(w_1,w_2)=\mathsf{true}$. Letting $\mathsf{rev}(\overline{w_1}w_2)=c\overline{d}$, we know we have $\mathsf{solve^+}(d,c)=\mathsf{true}$. By correctness of $\mathsf{solve^+}$, we know that $c\equiv^+d$. Thus $c\equiv d$ and $c\overline{d}\equiv\varepsilon$. 

    Returning to the fact that $\mathsf{rev}(\overline{w_1}w_2)=c\overline{d}$, we know by Theorems \ref{thm:reverse_equiv} and \ref{thm:reversing_to_braid_group_equivalence} that $\overline{w_1}w_2\equiv c\overline{d}$, and so $\overline{w_1}w_2\equiv \varepsilon$ as desired.

    Reverse direction : \leanlink{https://github.com/hannahfechtner/braids_better/blob/098309403cda787eaf1fa954e07ff0ca3b2927ee/BraidProject/Solver/GroupCorrectnessHardDirection.lean#L7-L24}Assume $w_1\equiv w_2$. Then $\overline{w_1}w_2\equiv 1$. Let $\mathsf{rev}(\overline{w_1}w_2)=c\overline{d}$. Then we have $\overline{w_1}w_2\equiv c\overline{d}$ again by Theorems \ref{thm:reverse_equiv} and \ref{thm:reversing_to_braid_group_equivalence}, and so $c\overline{d}\equiv 1$. Thus, $c\equiv d$, and so by Ore localization, $c\equiv^+d$. By correctness of our monoid solver, we know $\mathsf{solve}^+(d,c)=\mathsf{true}$. As $\operatorname{solve}(w_1, w_2)=\mathsf{solve}^+(d,c)$, we have our desired conclusion.
\end{proof}

Et voila! Thus, we have a terminating and correct algorithm to solve the braid word problem! 

\noindent Given two words in the braid group, we can run \lstinline{braid_group_word_solver}. For example,

\leanlink{https://github.com/hannahfechtner/braids_better/blob/098309403cda787eaf1fa954e07ff0ca3b2927ee/BraidProject/Solver/Examples.lean#L6-L7}
\begin{lstlisting}
#eval braid_word_solver [(4, false), (1, true), (2, true), (1, true)]
  [(2, true), (1, true), (2, true), (4, false)]
\end{lstlisting}
returns true, as the two words represent equal elements in the group.

\leanlink{https://github.com/hannahfechtner/braids_better/blob/098309403cda787eaf1fa954e07ff0ca3b2927ee/BraidProject/Solver/GroupCorrectnessHardDirection.lean#L70-L88}Or given two elements in the group $B$, we can use two applications of \lstinline{Quotient.exists_rep} to build \lstinline{braid_solver}. For example, we can run

\leanlink{https://github.com/hannahfechtner/braids_better/blob/098309403cda787eaf1fa954e07ff0ca3b2927ee/BraidProject/Solver/Examples.lean#L13}
\begin{lstlisting}
#eval braid_solver (σi 1 * σi 2 * σi 1) (σi 2 * σi 1 * σi 2 * (σi 3)⁻¹* (σi 1 * σi 3)⁻¹)
\end{lstlisting}
This comes with the specification that

\leanlink{https://github.com/hannahfechtner/braids_better/blob/098309403cda787eaf1fa954e07ff0ca3b2927ee/BraidProject/Solver/GroupCorrectnessHardDirection.lean#L90-L93}\begin{lstlisting}
theorem braid_solver_correct {a b : BraidGroupInf} : braid_solver a b ↔ a = b 
\end{lstlisting}

Results have been instantaneous using \lstinline{#eval} on small cases (approximately .33 seconds using the \lstinline{#time} command). Comparing two elements of $B_8$ each with length $500$ took about 17 seconds. The runtime of the code could certainly be improved. 
A number of constructive definitions carry unnecessary data -- data that makes proofs easier (by showing as an object is defined that a given data-carrying property holds, instead of later proving the property by structural induction), but is technically unnecessary. While \lstinline{Prop}-valued properties are ignored by the compiler, data-carrying properties, such as \lstinline{SignedList.PosNegData} or \lstinline{SemiThueData} are computed, and so slow down the calculation. The main contribution of this project is the mathematical theory, and so we have not spent significant time attempting to speed up code runtime (especially as it would come at the cost of much more work to set up the proofs).

\subsection{Finite Braid Group Solver}

We now return at last to the case of the finite braid groups. Note that we have been able to avoid the work of implementing Ore localization and grids for finite braid groups -- we haven't even shown simple monoid properties like common multiples nor cancellativity.

The solver is incredibly simple : we cast the finite braid group words to infinite braid group words, and run our braid group solver. In this section, unlike Section \ref{sec:finite_def}, we are careful to keep track of the following casts:

\begin{definition} \leanlink{https://github.com/hannahfechtner/braids_better/blob/098309403cda787eaf1fa954e07ff0ca3b2927ee/BraidProject/SignedListBounded.lean#L37-L38}
    Let $i_n$ be the inclusion map from $BW_n$ to $BW$. Let $s_n$ map appropriately-bounded elements of $BW$ to $BW_n$. \leanlink{https://github.com/hannahfechtner/braids_better/blob/098309403cda787eaf1fa954e07ff0ca3b2927ee/BraidProject/Additions/FreeMonoid.lean#L170-L171}(By abuse of notation, we also let both of these functions be used with $BW^+$ and $BW^+_n$ in place of $BW$ and $BW_n$ respectively).
\end{definition}

\begin{definition} \leanlink{https://github.com/hannahfechtner/braids_better/blob/098309403cda787eaf1fa954e07ff0ca3b2927ee/BraidProject/Solver/Finite.lean#L151-L152}
    We define the finite braid group word solver as $\mathsf{solve_n}(a,b):=\mathsf{solve}\big(i_n(a),i_n(b)\big)$ for $a,b\in BW_n$.
\end{definition}
Termination is thus immediate. Correctness requires one main lemma about reversing derivations, shown by induction :  

\begin{lemma}\label{lem:rev_bounded}\leanlink{https://github.com/hannahfechtner/braids_better/blob/098309403cda787eaf1fa954e07ff0ca3b2927ee/BraidProject/Solver/Finite.lean#L46-L87}
    Given $c\in BW$, if the index of the largest generator occurring in $c$ is strictly less than $n-1$, and $c \xrightarrow[R_{\mathrm{rev}}]{*} d$, then the index of the largest generator occurring in $d$ is also strictly less than $n-1$, and moreover, $s_n(c)\equiv_n s_n(d)$.
\end{lemma}

\begin{theorem} \leanlink{https://github.com/hannahfechtner/braids_better/blob/098309403cda787eaf1fa954e07ff0ca3b2927ee/BraidProject/Solver/Finite.lean#L154-L165}
    Given $a,b\in BW_n$, $\mathsf{solve_n}(a,b)=\mathsf{true}$ iff $a\equiv_nb$.
\end{theorem}

\begin{proof}
    Forward direction : \leanlink{https://github.com/hannahfechtner/braids_better/blob/098309403cda787eaf1fa954e07ff0ca3b2927ee/BraidProject/Solver/Finite.lean#L94-L138}Assume $\mathsf{solve_n}(a,b)=\mathsf{true}$. Then $\mathsf{solve}\big(i_n(a),i_n(b)\big)=\mathsf{true}$. Letting $\mathsf{rev}\big(\overline{(i_n(a)}, i_n(b)\big)=c\overline{d}$, we know we have $\mathsf{solve^+}(c,d)=\mathsf{true}$. By correctness of $\mathsf{solve^+}$, we know we have $c\equiv^+d$. Returning to the fact that $\mathsf{rev}\big(\overline{(i_n(a)}, i_n(b)\big)=c\overline{d}$, we know that $\overline{(i_n(a)}, i_n(b)\xrightarrow[R_{\mathrm{rev}}]{*} c\overline{d}$. We know that the index of the largest generator in $\overline{(i_n(a)}, i_n(b)$ must be strictly less than $n-1$, so by Lemma \ref{lem:rev_bounded}, the same holds for the word $c\overline{d}$, and hence $c$ and $d$ themselves. Thus Theorem \ref{thm:braid_monoid_inf_to_fin} allows us to conclude $s_n(c)\equiv^+_n s_n (d)$  from $c\equiv d$. We obtain $s_n(c\overline{d})\equiv\varepsilon$.  By the second part of Lemma \ref{lem:rev_bounded}, we know from $\overline{(i_n(a)}, i_n(b)\xrightarrow[R_{\mathrm{rev}}]{*} c\overline{d}$ that
 $s_n(\overline{(i_n(a)}, i_n(b)) \equiv_n s_n(c\overline{d})$. As $s_n(c\overline{d})\equiv\varepsilon$, we have $\overline{a} b \equiv_n \varepsilon$, so we have $a\equiv_nb$ as desired.

    Reverse direction : Assume $a\equiv_nb$. Then $i_n(a)\equiv i_n(b)$ and so $\mathsf{solve}(i_n(a),i_n(b))=\mathsf{true}$. Hence $\mathsf{solve_n}(a,b)=\mathsf{true}$ as desired.
\end{proof}

As a consequence of the correctness of this algorithm, we may now prove one elusive fact from the end of Section \ref{sec:finite_def} :

\begin{theorem}\leanlink{https://github.com/hannahfechtner/braids_better/blob/098309403cda787eaf1fa954e07ff0ca3b2927ee/BraidProject/Solver/Finite.lean#L167-L172}
    If $a,b\in BW$ have no generator with index $n-1$ or larger, and $a\equiv b$, then $s_n(a)\equiv_n s_n(b)$.
\end{theorem}

\begin{proof}
    Assume $a\equiv b$. By correctness of $\mathsf{solve}$, we know $\mathsf{solve}(a,b)=\mathsf{true}$. We know $a=i_n(s_n(a))$ and the analogous property holds for $b$, so $\mathsf{solve}(i_n(s_n(a)),i_n(s_n(b)))=\mathsf{true}$. By correctness of $\mathsf{solve_n}$, we have $s_n(a)\equiv_n s_n(b)$ as desired.
\end{proof}
\section{Future Work}

The most tempting work to be completed is a proof of Conjecture \ref{conj:same_type_same_length} (any two partial grids with the same frame have the same size). On the topic of grids, another improvement would be to generalize their definition to one parametrized by a set of relations. Further, we could classify which relations lead to grid existence and/or uniqueness, following Dehornoy's many papers which provide criteria such as right complemented presentations, the cube criterion, homogeneity, etc \cite{dehornoy_cancellativity_2018,  dehornoy_complete_2003, dehornoy_completeness_2000, dehornoy_garside_2019_fixed}. Such work could eventually lead to a generalization of the word problem solvers presented here to Garside groups and others. 

There are also other algorithms to solve the word problem which rely on the results developed in this project. For example, Garside's method of universal denominators (as mentioned in the introduction) could now be implemented fairly easily. We have also completed initial work necessary for more complex, popular, and efficient-in-practice methods such as Dehornoy's handle reduction. A formalization of this algorithm could also be useful to other projects, such as Hurtig's generalization of braid group for machine knitting, which uses the normal form from handle reduction as part of his algorithm to determine equivalence in a certain braided monoidal category \cite{hurtig_polynomial-time_2025}. (There is a plan to complete a Lean formalization of the mathematics in this paper).\footnote{\url{https://github.com/nhurtig/lean-category}} Other projects are working in the same vein of category-theoretic braids, such as Marco's formalization project.\footnote{\url{https://marcodavid.net/ainfinity/}} It would thus be potentially useful to identify the braid group with its category-theoretic generalization, as discussed in the introduction.

A larger goal is to formalize topological definitions of the braid group (for example, one under homotopy equivalence and another under ambient isotopy), and prove that the algebraic definition matches the topological ones. (There is some emerging, unpublished activity from AI theorem provers on a related problem using the configuration-space definition of braid groups.)

Beyond formalization of proofs, I am also working on a Lean widget to visualize braids. Given a formalized algorithm, the widget can show a braid transforming as the algorithm is applied to it (such as building its fraction form). In the reverse order, moving from imagery to Lean code, the strings in the widget will be able to be clicked and dragged in order to generate proofs.

\section{Acknowledgements}
I would like to thank Jeremy Avigad for his thoughtful feedback and careful edits throughout the project. I am grateful to Patrick Massot for his advice on the formalization of Ore localization for presented monoids. I thank Jim McCann, Jenny Lin, and Nat Hurtig for their brainstorming contributions concerning the definitions of partial grids, and Bhavik Mehta for his enthusiastic, spirited discussions of alternative definitions of partial grids, which led to the inclusion of Section \ref{sec:alt_choices}. I also thank the Mathlib community (especially Yaël Dillies, Johann Commelin, and numerous other maintainers) for their suggestions and assistance, particularly regarding the definitions of presented monoids and braid groups. Lastly, this project has been greatly bolstered by the congenial atmosphere of the Hoskinson Center for Formal Mathematics built by all the students and visitors who have gathered together there to talk Lean and support one another.

The first version of the codebase was completed without AI assistance. During refinements, Github’s Copilot’s line-completion tool was occasionally used (every output was immediately checked by hand). Claude Opus 5.5 was asked to search for duplicate code, summarize the contents of files, and suggest alternatives for theorem and function names.

The prose of the paper used no AI tools. ChatGPT 5 and 6 assisted in Latex formatting. Claude Opus 5.5 was sometimes used to search the codebase for formalized statements given a Latex statement (it provided the name of the theorem and the file in which it was found; the Github links were copied by me). 

This work was supported in part by NSF Grant DMS-2434614 and DARPA ExpMath Grant HR0011262E028.
\newpage

\printbibliography[heading=bibintoc,title={References}]

\newpage \appendix
\section{Maps}

Here, we summarize the many maps used between braid structures in this project. First, we give a diagram of functions between the various infinite braid structures. In the paper, some are often left implicit, such as $i^+$, but all are carefully noted in the Lean formalization.
\medskip
\[\begin{tikzcd}
	{BW_{pos}} & BW && {BW^{\mathds{1}}} & {BW^{\mathds{1}}_{pos}} \\
	&& B & \\
	&& {B^+} \\
	&& {BW^+}
	\arrow["\subset"{description}, draw=white, from=1-1, to=1-2]
	\arrow["{s^+}"{description}, curve={height=30pt}, from=1-1, to=4-3]
	\arrow["{i^\mathds{1}}"{description}, shift left=3, from=1-2, to=1-4]
	\arrow["q"', from=1-2, to=2-3]
	\arrow["{s^\mathds{1}}"{description}, from=1-4, to=1-2]
	\arrow["{q^\mathds{1}}", from=1-4, to=2-3]
	\arrow["\supset"{description}, draw=white, from=1-5, to=1-4]
	\arrow["{s^{\mathds{1}, +}}"{description}, curve={height=-30pt}, from=1-5, to=4-3]
	\arrow["\iota"', from=3-3, to=2-3]
	\arrow["{i^+}"{description}, curve={height=-24pt}, from=4-3, to=1-2]
	\arrow["{i^{\mathds{1},+}}"{description}, curve={height=24pt}, from=4-3, to=1-4]
	\arrow["{q^+}"', from=4-3, to=3-3]
\end{tikzcd}\]

\begin{longtable}{@{}llll@{}}
\label{tab:maps}\\

\toprule
\textbf{Map} &
\textbf{Domain} &
\textbf{Codomain} &
\textbf{Lean name} \\
\midrule
\endhead

$q$ & $BW$ & $B$ & \lstinline{BraidGroupInf.mk} \\

$q^{\mathds{1}}$ & $BW^{\mathds{1}}$ & $B$ &
\lstinline{BraidGroupInf.mk ∘ SignedOptionList.toSignedList} \\

$q^+$ & $BW^+$ & $B^+$ & \lstinline{BraidMonoidInf.mk} \\

$i^{\mathds{1}}$ & $BW$ & $BW^{\mathds{1}}$ &
\lstinline{SignedList.toSignedOptionList} \\

$s^{\mathds{1}}$ & $BW^{\mathds{1}}$ & $BW$ &
\lstinline{SignedOptionList.toSignedList} \\

$i^+$ & $BW^+$ & $BW$ & \lstinline{to_horizontal_edge_no_epsilon} \\
&             &       & \hspace{1em}(\lstinline{to_vertical_edge_no_epsilon} is the composition of \\
& & & \hspace{1.25em}the overbar operator $\overline{\hspace{.25em}\cdot\hspace{.25em}}$ with \lstinline{to_horizontal_edge_no_epsilon}) \\

$i^{\mathds{1},+}$ & $BW^+$ & $BW^{\mathds{1}}$ &
\lstinline{to_horizontal_edge} \\
\multicolumn{4}{@{}l@{}}{%
     \quad (note that $i^{\mathds{1},+}(\varepsilon):=\mathds{1}$; in Lean, we have \lstinline{to_horizontal_edge [] := [(none, true]})
} \\
&             &       & \hspace{1em}(\lstinline{to_vertical_edge} is the composition of \\
& & & \hspace{1.25em}the overbar operator $\overline{\hspace{.25em}\cdot\hspace{.25em}}$ with \lstinline{to_horizontal_edge}) \\

$s^+$ & $BW_{pos}$ & $BW^+$ & \texttt{SignedList.toList} \\

$s^{\mathds{1},+}$ & $BW^{\mathds{1}}_{pos}$ &
$BW^+$ & \texttt{SignedOptionList.toList} \\

$\iota$ & $B^+$ & $B$ & \lstinline{BraidMonoidInf.toBraidGroupInf} \\
\addlinespace[0.5em]
\multicolumn{4}{@{}l@{}}{%
    \emph{Note:} In Lean, the two inclusions are denoted using
    \lstinline{SignedList.is_true}.
} \\
\bottomrule
\end{longtable}

\newpage
We next consider functions between the finite braid structures, as well as maps to and from infinite braid structures. We define $BW_{\text{bounded}_n}$ as the set of infinite braid group words for which the index of the largest generator is strictly less than $n-1$. 
\[\begin{tikzcd}
	{BW^+_{\text{bounded}_n}} && {BW^+} & {BW_{\text{bounded}_n}} && BW \\
	\\
	& {BW^+_n} & {} && {BW_n} \\
	& {B_n^+} &&& {B_n} \\
	& {B^+} &&& B
	\arrow["\subset"{description}, draw=none, from=1-1, to=1-3]
	\arrow["{s_n^+}", curve={height=12pt}, from=1-1, to=3-2]
	\arrow["\subset"{description}, draw=none, from=1-4, to=1-6]
	\arrow["{s_n}", curve={height=12pt}, from=1-4, to=3-5]
	\arrow["{i_n^+}", curve={height=12pt}, from=3-2, to=1-3]
	\arrow["{f_n}", from=3-2, to=3-5]
	\arrow["{q_n^+}"', from=3-2, to=4-2]
	\arrow["{i_{n}}", curve={height=12pt}, from=3-5, to=1-6]
	\arrow["{q_n}", from=3-5, to=4-5]
	\arrow["{\iota_n}"', from=4-2, to=4-5]
	\arrow["{j_n^+}"', from=4-2, to=5-2]
	\arrow["{j_n}", from=4-5, to=5-5]
\end{tikzcd}\]

\begin{longtable}{@{}llll@{}}

\label{tab:finite_infinite_maps}\\

\toprule
\textbf{Map} &
\textbf{Domain} &
\textbf{Codomain} &
\textbf{Lean name} \\
\midrule
\endhead

$s_n^+$ & $BW^+_{\text{bounded}_n}$ & $BW_n^+$ &
\lstinline{FreeMonoid.mapNatToFin} \\

$s_n$ & $BW_{\text{bounded}_n}$ & $BW_n$ &
\lstinline{SignedList.mapNatToFin} \\

$i_n^+$ & $BW_n^+$ & $BW^+$ &
\lstinline{List.map (fun x => x.1)} \\

$i_n$ & $BW_n$ & $BW$ &
\lstinline{List.map (fun x => (x.1.val, x.2))} \\

$f_n$ & $BW_n^+$ & $BW_n$ &
\lstinline{to_horizontal_edge_no_epsilon} \\

$q_n^+$ & $BW_n^+$ & $B_n^+$ &
\lstinline{BraidMonoidFin.mk n} \\

$q_n$ & $BW_n$ & $B_n$ &
\lstinline{BraidGroupFin.mk n} \\

$\iota_n$ & $B_n^+$ & $B_n$ &
\lstinline|BraidMonoidFin.toBraidGroupFin n| \\

$j_n^+$ & $B_n^+$ & $B^+$ &
\lstinline|BraidMonoidFin.toBraidMonoidInf n| \\

$j_n$ & $B_n$ & $B$ &
\lstinline|BraidGroupFin.toBraidGroupInf n| \\

\addlinespace[0.5em]
\multicolumn{4}{@{}l@{}}{%
    \emph{Note:}
    To implement the inclusion maps in Lean, we use \lstinline{SignedList.bounded} to define $BW_{\text{bounded}_n}$,
} \\
\multicolumn{4}{@{}l@{}}{%
    \qquad \quad and hypotheses of the form \texttt{$\forall$}\lstinline{ x ∈ w, x < n.pred} to discuss $BW^+_{\text{bounded}_n}$.
} \\
\bottomrule
\end{longtable}

\newpage 
\section{Dependency Graph}
Here is the dependency graph for the project: \leanlink{https://github.com/hannahfechtner/braids_better/blob/098309403cda787eaf1fa954e07ff0ca3b2927ee/braid_only.png}

\begin{center}
\includegraphics[width=.8\paperwidth, keepaspectratio]{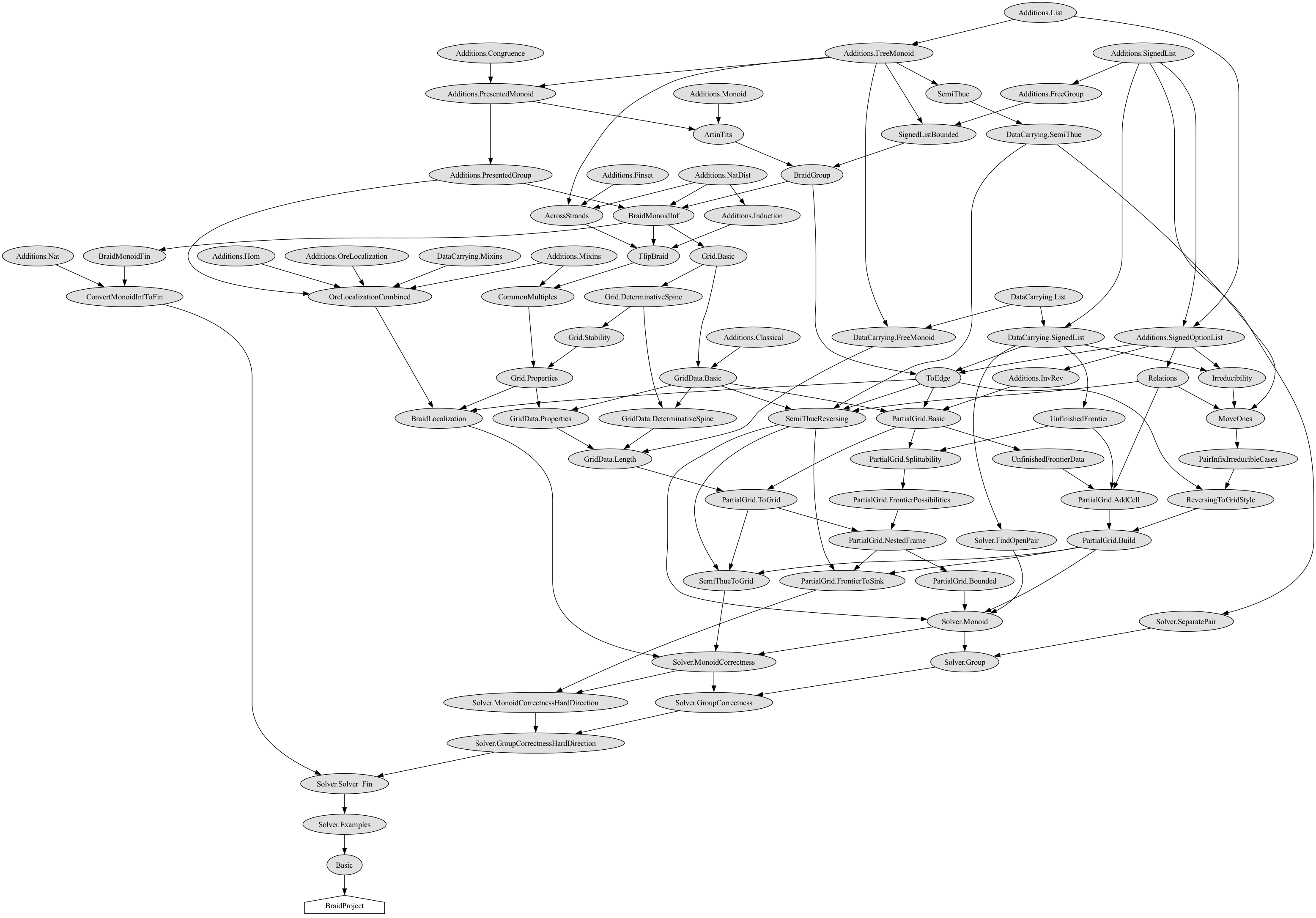}
\end{center}

\noindent The repository currently contains about 15,000 lines of code.

\end{document}